\documentclass[11pt]{article}
\renewcommand{\theequation}{\arabic{section}.\arabic{equation}}

\newtheorem{example}{Example}[section]

\newtheorem{note}{Note}[section]
\usepackage{bm}
\usepackage{amssymb}
\usepackage{amsmath}
\usepackage{epsfig}

\begin{document}
\begin{center}
{\LARGE\bf Nonparametric estimates of low bias}\\[1ex]
by\\[1ex]
Christopher S. Withers\\
Applied Mathematics Group\\
Industrial Research Limited\\
Lower Hutt, NEW ZEALAND\\[2ex]
Saralees Nadarajah\\
School of Mathematics\\
University of Manchester\\
Manchester M13 9PL, UK
\end{center}
\vspace{1.5cm}
{\bf Abstract:}~~We consider the problem of estimating an arbitrary
smooth functional of $k \geq  1 $ distribution functions (d.f.s.) in
terms of random samples from them.
The natural estimate replaces the
d.f.s by their empirical d.f.s.
Its bias is generally $\sim n^{-1}$,
where $n$ is the minimum sample size, with a {\it $p$th order} iterative estimate of bias $ \sim n^{-p}$ for any $p$.
For $p \leq 4$, we give an explicit estimate in terms of the first $2p - 2$ von
Mises derivatives of the functional evaluated at the empirical d.f.s.
These may be used to obtain {\it unbiased} estimates, where these exist and
are of known form in terms of the sample sizes; our form for such unbiased
estimates is much simpler than that obtained using polykays and tables of the
symmetric functions.
Examples include functions of a mean vector (such as the ratio of two means and the
inverse of a mean), standard deviation, correlation, return times and exceedances.
These $p$th order estimates require only $\sim n $ calculations.
This is in sharp contrast with computationally intensive bias reduction methods such as the
$p$th order bootstrap and jackknife, which require $\sim n^p $ calculations.

\noindent
{\bf AMS 2000 Subject Classifications:}~~Primary 62G05; Secondary 62G30.

\noindent
{\bf Keywords and Phrases:}~~Bias reduction; Correlation; Exceedances;
Multisample; Multivariate; Nonparametric; Ratio of means;
Return times; Standard deviation; von Mises derivatives.

\section{Introduction and summary}
\setcounter{equation}{0}

Let $T (F)$ be any {\it smooth functional} of one or more unknown
distributions $F$ based on random samples from them.
Bias reduction of estimates of $T (F)$, say $T (\widehat{F})$, has been a subject of considerable interest.
Traditionally bias reduction has been based on well known resampling methods like bootstrapping
and jackknifing in nonparametric settings, see Gray and Schucany (1972) and Efron (1982).
However, these methods may not be effective in complex situations when the sampling distribution
of the statistic changes too abruptly with the parameter, or when this distribution is very skewed and has heavy tails.
Also the robustness properties of $F$ may not be preserved for $T (F)$ for all $T (\cdot)$.

Recently, various analytical methods have been developed for bias reduction in parametric settings.
Withers (1987) developed methods for bias reduction based on Taylor series expansions.
Sen (1988) established asymptotic normality of
$\sqrt{n} \{ T (\widehat{F}) - T (F) \}$ as $n \rightarrow \infty$
under suitable regularity conditions.
Cabrera and Fernholz (1999, 2004) defined a {\it target estimator}:
for a given $T$ and a parametric family of distributions it is defined by setting
the expected value of the statistic equal to the observed value.
Cabrera and Fernholz (1999, 2004) established under suitable regularity conditions that
the target estimator has smaller bias and mean squared error than the original estimator.
See also Fernholz (2001).

Suppose we have $k \geq  1$ independent samples of sizes $n_1, \ldots, n_k$
from distribution functions (d.f.s) $F = (F_1, \ldots, F_k)$ on $R^{s_1}, \ldots, R^{s_k}$.
Let $\widehat{F} = ( \widehat{F}_1,\ldots, \widehat{F}_k)$ denote their
sample d.f.s and let $n$ be the minimum sample size.
The problem we consider in this paper is that of finding an estimate of low bias
for an arbitrary smooth functional $T(F)$.
The natural estimate $T \widehat{(F)}$  generally has bias  $\sim n^{-1}$, that is, $ O(n^{-1})$ as $n \rightarrow \infty$.

This paper has already been cited as an unpublished technical report in Withers and Nadarajah (2008).
The estimates proposed here have been compared to alternatives.
We showed in particular that our estimates consistently outperform bootstrapping, jackknifing and those due to
Sen (1988)  and Cabrera and Fernholz (1999, 2004).
We also provided computer programs in MAPLE for implementation of the proposed estimates.

The emphasis of this paper is to describe how to find estimates of low bias for $T (F)$.
Because of the material in Withers and Nadarajah (2008), the emphasis here will not be on numerical illustrations or applications.

For the reader's convenience, in Section 2, we repeat the definition of functional
derivatives and rules for obtaining them given in Withers (1988).
In Section 3, we have a formal asymptotic expansion of the form
\begin{eqnarray}
E T \widehat{(F)} = \sum_{r=0}^{\infty} n ^{-r} C_r,
\label{1.1}
\end{eqnarray}
where $C_0 = T(F)$.
The coefficient of $n^{-r}$ in $E T (\widehat{F})$, $C_r(F, T) = C_r$ may be written in terms of the (functional or von Mises) derivatives
of $T ( \widehat{F})$ of order $ \leq 2 r$, and is given there  explicitly for $r \leq  4$.

From (\ref{1.1}) if a functional $T_{(n)}(F)$ can be expanded as
\begin{eqnarray*}
T_{(n)} = \sum^\infty_{i=0} n^{-i} T_i(F)
\end{eqnarray*}
then
\begin{eqnarray*}
ET_{(n)}\left(\widehat{F}\right)= \sum^\infty_{i=0}n^{-i} C_i\left( {\bf T} \right),
\end{eqnarray*}
where
\begin{eqnarray*}
C_i({\bf T}) = \sum^i_{j=0} C_j  \left(F,T_{i-j}\right).
\end{eqnarray*}
Defining $T_i $ iteratively by $T_0 = T$ and
\begin{eqnarray}
T_i (F) = -\sum_{j=1}^i C_j \left(F, T_{i -j}\right)
\label{1.2}
\end{eqnarray}
for $i \geq 1$ it follows that for $p \geq 1$
\begin{eqnarray}
T_{np} (F) = \sum_{i=0}^{p-1} n^{-i} T_i (F)
\label{1.3}
\end{eqnarray}
satisfies a formal expansion of the form
\begin{eqnarray*}
E T_{np} \widehat{(F)} =  T (F) + O \left(n^{-p}\right).
\end{eqnarray*}
So, $T_{np} (\widehat{F})$ is {\it a $p$th order estimate} in the sense that it has bias $O (n^{-p})$.
This result was given for the case $k=1$, $p=2$ using
a different approach in an unpublished technical report by Jaeckel (1972).

Note that $T_i (\widehat{F})$ given
by (\ref{1.2}) is the coefficient of $n^{-i}$ in the expansion in powers of
$n^{-1}$ of the unbiased estimate (UE) of $T (F)$, if an UE exists.

Section 4  gives $T_i (F)$ explicitly in terms of the first $2i$ derivatives of
$T(F)$ for $i \leq  3$.
So, $T_{n4}\widehat{(F)}$ is an explicit
estimate of bias $ O( n^{-4})$.
Note 4.1 shows how to obtain from (\ref{1.3}) an estimate of bias $ O (n^{-p})$ of the form $S_{np}
(\widehat{F})$, where
\begin{eqnarray*}
S_{np} (F) = \sum_{i=0}^{p-1} S_i (F) / \left\{ (n-1) \cdots (n-i) \right\}.
\end{eqnarray*}
This estimate is {\it unbiased} for one sample if $T(F)$
is a polynomial in $ F$ (such as a moment or cumulant) of degree up to $p$.

Section 5 gives examples and makes comparisons with the UEs
of central moments and cumulants given by James (1958) and by Fisher (1929).
Our method is demonstrated to be give much simpler results for UEs of products
of moments than the polykay system of Wishart (1952) as expounded in
Section 12.22 of Stuart and Ord (1987) using tables of the symmetric functions.

Examples 5.1 to 5.3 estimate an arbitrary
function of the vector ${\bm \mu} (F)$, the mean of one multivariate distribution.
Example 5.2   specializes to $T (F) = {\bf a}^{\prime} {\bm \mu} (F) / {\bf b}^{\prime} {\bm \mu} (F)$,
where ${\bf a}$, ${\bf b}$ are given $s_1$-vectors, in particular for the ratio of means of a bivariate sample,
\begin{eqnarray*}
T (F) = \mu_1 (F) / \mu_2 (F).
\end{eqnarray*}
Examples 5.4  and 5.5 estimate an arbitrary function of the means of $k$
univariate distributions; in particular it
considers the case of two univariate samples $(k = 2, s_1 = s_2 = 1)$ with
\begin{eqnarray*}
T (F) = \mu \left(F_1\right) / \mu \left(F_2\right).
\end{eqnarray*}
Example 5.6 gives an explicit expression for the general derivative of the
$r$th central moment $\mu_r$.
Together with the chain rule of Appendix A this
enables one to obtain a $p$th order estimate of any smooth function of
moments.
In particular, we give fourth order estimates for {\it any central moment} and
UEs for $\mu_r$ for $r \leq 7$.

Examples 5.7 to 5.11 extend this to an arbitrary product of moments.
An alternative matrix method for obtaining  UEs of
products of moments is given there.
This involves obtaining
simultaneously the UEs of all moment products of a given degree.
Examples 5.12 to 5.15 give fourth order estimates of the standard deviation and functions of it.
Example 5.16 gives third order estimates of the ratio of the mean to the standard deviation.

Examples 5.17 to 5.21 give applications to return times and
exceedances.
Examples 5.22 and 5.23 illustrate how to obtain UEs for multivariate moments and cumulants from univariate
analogs.
Finally, Examples 5.24 and 5.25 give second order estimates for the correlation and its square.

The method can also be used to estimate with reduced bias any
cumulant of $T(\widehat{F})$.
This is illustrated in Section 6 which gives
a third order estimate for the covariance of
any estimate of the form ${\bf T} \widehat{(F)}$, where now ${\bf T}$ may
be a vector.
For example, by Example 5.1, if $k = 1 $ and $T(F)$ is any
function of ${\bm \mu} (F)$ (such as $\mu_1 (F) /  \mu_2 (F)$) if $s_1 = 2)$,
this estimate is a function of the mean and covariance of $F$
only, whereas $C_1$ depends also on the third moment.

Section 7 shows how to estimate the covariance of an estimate of bias.

There are, of course, other $p$th order estimates of $T(F)$, but they are all
{\it computationally intensive}, requiring $O( n^p) $ calculations
(except in special cases), whereas {\it our
method requires only $O( n) $ calculations} for fixed $p$.
The main examples are, firstly, the $(p-1)$th iterated bootstrap, $\widehat{\theta}_{p-1}$ of equation (1.35) of Hall (1992)
in which $(-1)^{i+1}$ should be inserted in the right hand side; and, secondly, the
$p$th order jackknife $\widehat{\theta}^{p-1}$ of equation (4.17) of Schucany {\it et al.} (1971), a ratio of
$p\times p$ determinants.
To see that this requires $O( n^p)$ calculations note that $t_p$ of their equation (4.19) requires $O( n^p)$ calculations.

The techniques given here can also be applied to quantify their biases.
Note that if $A$ and $B$ are two $p$th order estimates of $T(F)$ then $A - B = O_p(n^{-p})$.

Appendix A gives a very useful chain rule for obtaining the derivatives of a function of a functional.
Appendix B gives some results used to obtain $\{T_i\}$ of (\ref{1.3}).
Appendix C shows how to estimate the number of simulated samples
needed to estimate the bias to within a given relative error.

Tiit (1988) by an entirely different method obtained
an expansion of the form (\ref{1.1}) for
\begin{eqnarray*}
m ({v})=T (F) =\prod_{i=1}^{s} E X^{v_{i}},
\end{eqnarray*}
where $ X \sim F$, and so also for $\mu_{r} (F)$.
For these cases he constructs estimates of bias $O (n^{-p})$ given $p \geq 1$.
He shows for $T(F) = m ({v})$ that the UE $T_{n \infty} (\widehat{F})$ converges if
$E | X |^{h} < \infty$, where $ h=\sum_{i=1}^{s} v_{i}$ and $n - 1 >$ the number of partitions of $h$.
His expression on page 12, Theorem 4, is incorrect.
He gives
\begin{eqnarray*}
var \ \widehat{m} ({v} ) = n^{-1} V + O \left(n^{-2}\right),
\end{eqnarray*}
where
\begin{eqnarray*}
V = m
\left({v}\right)^{2} \left( A - s^{2} \right), \mbox{ and } A = \sum_{i=1}^{s} m_{2v_{i}} m_{v_{i}}^{2}.
\end{eqnarray*}
Here, $A$ should be
\begin{eqnarray*}
\sum_{i,j=1}^{s} m_{v_{i}+v_{j}} m_{v_{i}}^{-1} m_{v_{j}}^{-1}.
\end{eqnarray*}
For the case $T(F) = \mu^{3}$ his Table 2 illustrates through simulations for
$F=U (0,1)$ and $n=5, 10$ how the bias of $T_{np} (\widehat{F}) $ falls to zero as $p$ increases.

\section{Functional partial derivatives and notation}
\setcounter{equation}{0}

Let $ {\cal F}_{s}$ denote the space of d.f.s on $R^{s}$.
Let ${\bf x}, {\bf y}, {\bf x}_1, \ldots, {\bf x}_r$ be points in $R^{s}$,
$F \in {\cal F}_{s}$ and $T: \ \mathcal{F}_{s} \rightarrow R$.
In Withers (1983) the $r$-th order functional derivative of $T (F)$ at $ ({\bf x}_1, \ldots, {\bf x}_r)$
\begin{eqnarray*}
T_{{\bf x}_1, \ldots, {\bf x}_r}= T_{F} \left( {\bf x}_1, \ldots, {\bf x}_r \right),
\end{eqnarray*}
was defined.
It is characterized by the formal
functional Taylor series expansion:  for $ G$ in  $\mathcal{F}_{s}$,
\begin{eqnarray}
T (G) - T(F) \approx \sum_{r=1}^{\infty} \int^r T_{F} \left( {\bf x}_1, \ldots, {\bf x}_r \right) \prod_{j=1}^r d
\left( G \left({\bf x}_j\right) - F \left({\bf x}_j\right)\right) / r !,
\label{2.1}
\end{eqnarray}
where $\int^r$ denote $r$ integral signs, and the constraints
$T_{{\bf x}_1 \cdots {\bf x}_r}$ is symmetric in its $r$ arguments, and
\begin{eqnarray*}
\int T_{{\bf x}_1 \cdots {\bf x}_r} dF \left({\bf x}_1\right) = 0.
\end{eqnarray*}
These imply $F ({\bf x}_j)$ in (\ref{2.1}) can be replaced by zero.
In particular, it was shown that, for $0 \leq \varepsilon \leq 1$,
\begin{eqnarray*}
T_{x} = \partial T\left(F + \varepsilon \left(\delta_{x} - F\right)\right) / \partial \varepsilon
\end{eqnarray*}
at $\varepsilon = 0$, where $\delta_{x}$ is the d.f putting mass 1 at $x$,
that is $\delta_{x} (y) = I (x \leq y) = 1$ if  $x \leq y$ and $0$ otherwise.
For example, $T (F) = F(y)$ has first derivative $T_{x} = T_{F} (x) = \delta_{x} (y) - F(y)=F(y)_{x}$, say.

Also, $T_{{\bf x}_1\cdots {\bf x}_r} = 0$ if $T (F)$ is a `polynomial in $F$' of degree
less than $r$ (for example, a moment or cumulant of $F$ of
order less than $r$), so that the Taylor series in (\ref{2.1}) consists of
only $r-1$ terms.
Note that $T(F)$ is a polynomial in $F$ of degree $m$
if for any $G$ in $\mathcal{F}_{s}$, $T (F + \varepsilon (G - F))$ is a
polynomial in $\varepsilon$ of degree $m$.

Suppose now that $F = (F_1, \ldots, F_k)$ consists of $k$
distributions on $R^{s_1}, \ldots, R^{s_k}$ and that $T (F)$ is
a real functional of $F$.
Then the {\it functional partial derivative} of $T (F)$ at
$\left (  \begin{array}{ccc}
a_1&\cdots&a_r\\
{\bf x}_1&\cdots&{\bf x}_r\end{array} \right )$ is defined by
\begin{eqnarray*}
T_{{\bf x}_1 \cdots {\bf x}_r}^{a_1\cdots a_r} = T_{F} \left (
\begin{array}{ccc}
a_1& \cdots &a_r\\
{\bf x}_1& \cdots& {\bf x}_r \end{array}\right ),
\end{eqnarray*}
where ${\bf x}_i$ in $R^{s_{a_i}}$ and $a_i$ in $\{1, 2, \ldots, k\}$,
and is obtained by treating the lower order functional partial derivatives
and $T(F)$ as functionals of $F_a$ alone for $a=a_1, \ldots, a_r$.
For example, $T^{a\cdots a}_{{\bf x}_1 \cdots {\bf x}_r}$ is the ordinary
functional derivative of $S(F_a) = T (F)$ at $({\bf x}_1 \cdots {\bf x}_r)$, and
$T^{a \cdots ab \cdots b}_{{\bf x}_1 \cdots {\bf x}_r {\bf y}_1 \cdots {\bf y}_s}$
is the ordinary functional derivative of $S(F_b) = T^{a \cdots a}_{{\bf x}_1 \cdots {\bf x}_r}$
at $({\bf y}_1 \cdots {\bf y}_{s})$.

Just as \ $\partial^2 f(x,y) / \partial x \partial y = \partial^2 f (x,y)/ \partial y \partial x$ under mild conditions,
swapping columns
of $T \begin{array}c
a_1 \cdots a_r\\{\bf x}_1 \cdots {\bf x}_r\end{array}$ (for example,
${}^{a_1}_{{\bf x}_1}$ and ${}^{a_2}_{{\bf x}_2}$) will not alter its value.
So, $T^{a \cdots ab \cdots b}_{{\bf x}_1 \cdots {\bf x}_r {\bf y}_1 \cdots {\bf y}_s}$
is also the ordinary functional derivative of $S(F_a) = T^{b \cdots b}_{{\bf y}_1 \cdots {\bf y}_r}$ at $({\bf x}_1 \cdots {\bf x}_s)$.

The partial derivatives may also be characterized by the formal
functional Taylor series expansion: for $G=(G_1, \ldots, G_k)
\in \mathcal{F}_{s_1} \times \cdots \times \mathcal{F}_{s_k}$,
\begin{eqnarray}
T(G) - T(F) \approx \sum_{r=1}^{\infty} \int^r T_{F} \left (
\begin{array}c a_1 \cdots a_r\\{\bf x}_1 \cdots {\bf x}_r\end{array}
\right) \prod_{j=1}^r d \left( G_{a_j} \left({\bf x}_j\right) - F_{a_j} \left({\bf x}_j\right)\right) / r !
\label{2.6}
\end{eqnarray}
with summation of the repeated subscripts $a_1 \cdots a_r$ over
their range $1 \cdots p$ implicit, together with the constraints
\begin{eqnarray}
T^{a_1 \cdots a_r}_{{\bf x}_1 \cdots {\bf x}_r}
\mbox{ is not altered by swapping columns,}
\nonumber
\end{eqnarray}
and
\begin{eqnarray*}
\int T ^{a_1 \cdots a_r}_{{\bf x}_1 \cdots {\bf x}_r} d F_{a_1} \left({\bf x}_1\right) = 0.
\end{eqnarray*}
These imply $F_{a_j} ({\bf x}_j)$ in (\ref{2.6}) can be replaced by zero.
The partial derivatives may be calculated using
\begin{eqnarray}
T^a_{\bf x} = S_{\bf x} \mbox{ for } S\left(F_a\right) = T (F),
\label{2.9a}
\end{eqnarray}
and
\begin{eqnarray}
T \begin{array}ca_1\cdots a_{r+1}\\{\bf x}_1\cdots {\bf x}_{r+1}\end{array} =
\left ( T\,\begin{array}ca_1 \cdots a_r\\{\bf x}_1 \cdots {\bf x}_r\end{array}
\right ){}^{a_{r+1}}_{{\bf x}_{r+1}}
+ \sum_{i=1 }^r \delta_{a_i a_{r+1}} T \left\langle
{}^{a_1 \cdots a_{r+1}}_{{\bf x}_1 \cdots {\bf x}_{r+1}}
\right\rangle_i,
\label{2.9b}
\end{eqnarray}
where $\delta_{ij}=1$ or $0$ for $i=j$ or $i \neq j$ and $\langle \rangle_i$ means `drop the $i$-th column'.
For example, $T_{\bf xy}^{ab} + \delta_{ab} T_{\bf y}^b$.
The proof of (\ref{2.9a})-(\ref{2.9b}) is as for equation (2.6) of Withers (1983).

\section{Expansions for bias}
\setcounter{equation}{0}

Perhaps the easiest method to obtain expressions for the bias coefficients
$ \{C_r \}$ of (\ref{1.1}) and the bias reduction coefficients
$\{ T_i (F) \}$ of (\ref{1.3}) is from their parametric analogs, given in equation (A.1) and
Appendix  D (for $i \leq  3 $) of Withers (1987).
The method is to identify $( \theta, \widehat{\theta}, t, \sum )$ with $( F, \widehat{F}, T, \int)$,
where the integral is with
respect to the appropriate distribution $F_i$.
This method was used
in Withers (1988) to derive non-parametric confidence intervals of level $ 1 - \alpha + O(n^{-j/2})$ from their parametric analogs.
It is convenient to set
\begin{eqnarray}
T ( a^i b^j \dots ) = \int \cdots \int T_{F} \left( \begin{array}{cc}
a^i b^j \\ x^iy^j \end{array} \cdots \right)
dF_a (x) dF_b (y) \cdots,
\label{3.1}
\end{eqnarray}
where $x^i$ denotes a string of $i$ $x$'s (not a product) and
similarly, for $a^i$.
In the notation of Withers (1988) this is $[1^i 2^j \cdots]_{ab \cdots}$.
Setting
\begin{eqnarray}
\lambda_a = n/n_a \mbox{ with } n = \min n_i,
\label{3.2}
\end{eqnarray}
the above approach yields
\begin{eqnarray}
&&
C_1 = |2|/2,
\
C_2  =  |3| /6 + \left |2^2 \right | /8,
\label{3.3.1}
\\
&&
C_3 = |4| /24 + |23| /12 + \left |2^3 \right | /48,
\label{3.3.2}
\\
&&
C_4 = |5| /120 + |24|/48+\left|3^2\right|/72+\left |2^23\right |/48+\left |2^4\right |/384,
\label{3.3.3}
\end{eqnarray}
where
\begin{eqnarray}
&&
|2| = \sum \lambda_{a} T \left(a^{2}\right),
\nonumber
\\
&&
|3| = \sum \lambda_{a}^{2} T \left( a^{3} \right),
\nonumber
\\
&&
\left |2^{2} \right| = \sum \lambda_{a_{1}} \lambda_{a_{2}} T \left(a_{1}^{2} a_{2}^{2}\right),
\nonumber
\\
&&
|4| = \sum \lambda_{a}^{3} \left\{ T \left(a^{4}\right) - 3 T \left(a^{2} a^{2}\right) \right\},
\nonumber
\\
&&
|23| = \sum \lambda_{a} \lambda_{b}^{2} T \left(a^{2} b^{3} \right),
\nonumber
\\
&&
\left |2^{3} \right | = \sum \lambda_{a_{1}} \lambda_{a_{2}} \lambda_{a_{3}} T \left( a_{1}^{2} \ a_{2}^{2}  \  a_{3} ^{2}\right),
\nonumber
\\
&&
|5| = \sum \lambda_{a}^{4} \left\{ T \left(a^{5}\right) - 10 T \left(a^{2} a^{3}\right) \right\},
\nonumber
\\
&&
|24| = \sum \lambda_{a} \lambda_{b}^{3} \left\{ T \left( a^{2} b^{4} \right) - 3 T \left(a^{2} b^{2} b^{2}\right) \right\},
\nonumber
\\
&&
\left |3^{2} \right | = \sum \lambda_{a_{1}}^{2} \lambda_{a_{2}}^{2} \ T \left(a_{1}^{3} a_{2}^{3}\right),
\nonumber
\\
&&
\left | 2^{2}3 \right | = \sum \lambda_{a_{1}} \lambda_{a_{2}} \lambda_{b}^{2} T \left( a_{1}^{2} \ a_{2}^{2} \ b^{3} \right),
\nonumber
\\
&&
\left |2^{4}\right | = \sum \lambda_{a_{1}} \lambda_{a_{2}} \lambda_{a_{3}}
\lambda_{a_{4}} T \left( a_{1}^{2} \ a_{2}^{2} \ a_{3}^{2} \ a_{4}^{2} \right).
\nonumber
\end{eqnarray}
For example, if $k = 1$ (one sample) then
\begin{eqnarray}
C_1 = T \left(1^2\right)/2,
\
C_2 = T \left(1^3\right)/6 + T\left(1^2 1^2\right)/8, \ldots .
\label{3.6}
\end{eqnarray}
More generally,
\begin{eqnarray}
&&
\left |A^i \right | = \sum \lambda^{A-1}_{a_1} \cdots \lambda^{A-1}_{a_i} \left |A^i \right |_{a_1 \cdots a_i},
\nonumber
\\
&&
\left |A^i B^j \right | =  \sum \lambda^{A-1}_{a_1} \cdots
\lambda^{A-1}_{a_i}\lambda^{B-1}_{b_1} \cdots  \lambda^{B-1}_{b_j} \left | A^i B^j \right |_{a_1 \cdots a_i b_1 \cdots b_j}
\label{3.3.4}
\end{eqnarray}
with each $a_1, \ldots, b_j$ summed over $1, \ldots, k$,
\begin{eqnarray*}
&&
\left |A^i B^j \right |_{a_1 \cdots a_i b_1 \ldots b_j} = T \left(a^{A}_1 \cdots a^{A}_i b^{B}_1 \cdots b^{B}_j\right) \mbox{ if } A   \mbox{ and }  B=2 \mbox{ or } 3,
\\
&&
|4|_a = T \left(a^4\right) - 3T \left(a^2 a^2\right),
\\
&&
|5|_a = T \left(a^{5}\right) - 10T \left(a^2 a^3\right),
\\
&&
|24|_{ab} = T\left(a^2b^4 \right) -3T \left(a^2b^2b^2\right).
\end{eqnarray*}
For example,
\begin{eqnarray*}
\left |A^2 \right | = \sum \lambda^{A-1}_{a_1} \lambda^{A-1}_{a_2} \left |A^2 \right |_{a_1a_2},
\end{eqnarray*}
and
\begin{eqnarray*}
\left |A^2 \right |_{a_1a_2}
&=&
T \left( a_1^{A} a_2^{A} \right) \mbox{ if } A=2 \mbox{ or } 3
\\
&=&
\int \int T_{F} \left( \begin{array}{cc} a_1^{A} & a_2^{A}\\
x^{A}& y^{A}  \end{array}\right)dF_{a_1}(x) dF_{a_2}(y),
\end{eqnarray*}
so for the one sample case ($k=1$),
\begin{eqnarray*}
&&
\left |A^i \right | = T \left(1^{A} \cdots 1^{A}\right) \mbox{ if } A=2 \mbox{ or } 3,
\\
&&
\left |A^i B^j \right | = T \left(1^{A} \cdots 1^{A} 1^{B} \cdots 1^{B}\right) \mbox{ if } A  \mbox{ and }  B=2 \mbox{ or } 3,
\\
&&
|4| = T\left(1^4\right)-3T\left(1^2 1^2\right),
\
|5| = T\left(1^5\right)-10T\left(1^2 1^3\right),
\
|24|=T\left(1^2 1^4\right)-3T\left(1^21^21^2\right).
\end{eqnarray*}

\begin{note}
The general term $C_r$ is given by equation (A.1) of Withers (1987), (\ref{3.2}), (\ref{3.3.4}), and
\begin{eqnarray*}
|ij \cdots|_{ab \cdots} = \int^i d^i \kappa^{\prime}_a \left( {\bf x}_1 \cdots {\bf x}_i \right)
\int^j d\kappa^{\prime}_b \left( {\bf y}_1 \cdots {\bf y}_j \right)
\cdots T_F \left( \begin{array}c a \cdots ab \cdots b \\ {\bf x}_1 \cdots {\bf x}_i {\bf y}_1 \cdots {\bf y}_j \end{array} \cdots \right ),
\end{eqnarray*}
where $\int^i d^i \kappa^{\prime}_a ({\bf x}_1 \cdots {\bf x}_i)$ is the Lebesgue-Stieltjes integral,
\begin{eqnarray*}
&&
{\bf x}_1 \wedge {\bf x}_2 \wedge \cdots = \min \left({\bf x}_1, {\bf x}_2, \ldots \right) \mbox{ taken componentwise },
\\
&&
f_{12 \cdots} = F_a \left( {\bf x}_1 \wedge {\bf x}_2 \wedge \cdots \right),
\\
&&
\kappa_a \left( {\bf x}_1 {\bf x}_2 \cdots \right) = \kappa \left(  {\bf Y}_1, {\bf Y}_2, \ldots \right),
\mbox{the joint cumulant at } {\bf Y}_j = I \left( {\bf X}_a \leq {\bf x}_j \right),
\\
&&
\kappa_a^{\prime} \left( {\bf x}_1 {\bf x}_2 \cdots \right) = \kappa_a \left( {\bf x}_1 {\bf x}_2 \cdots \right)
\mbox{ expressed as a function of } \left\{ f_{ij \cdots} \right\} \mbox{ at } f_i \equiv 0,
\end{eqnarray*}
and $I$ is the indicator function and $X_a \sim F_a$.
For example, using an obvious summation notation
\begin{eqnarray*}
&&
\kappa_a \left( {\bf x}_1 {\bf x}_2 \right) = f_{12} - f_1 f_2,
\\
&&
\kappa_a \left( {\bf x}_1 {\bf x}_2 {\bf x}_3 \right) = f_{123} - \sum^3 f_{12} f_3 +2f_1 f_2 f_3,
\\
&&
\kappa_a \left( {\bf x}_1 \cdots {\bf x}_4 \right) = f_{1 \cdots 4} -\sum^4 f_{123} f_4 - \sum^3 f_{12} f_{34},
\end{eqnarray*}
imply
\begin{eqnarray*}
&&
\kappa_a^{\prime} \left( {\bf x}_1 {\bf x}_2 \right) = f_{12},
\
\kappa_a^{\prime} \left( {\bf x}_1 {\bf x}_2 {\bf x}_3 \right) = f_{123},
\
\kappa_a^{\prime} \left( {\bf x}_1, \ldots, {\bf x}_4 \right) = f_{1 \cdots  4} - \sum^3 f_{12} f_{34}.
\end{eqnarray*}
\end{note}

\begin{note}
As a check if $k  = 1$, $(C_1, C_2) = (a_{11}, a_{12})$ on page 580 of Withers (1983).
\end{note}

\section{Estimates of bias $O(n^{-4})$}
\setcounter{equation}{0}

Here, we give expressions for $\{T_i, i \leq 3\}$ of (\ref{1.2}) and for
$\{S_i,i \leq 3\}$ of Note 4.1.
Estimates of bias $O(n^{-4})$ are then given by $T_{n4} \widehat{(F)}$ of (\ref{1.3}) and $S_{n4}(\widehat{F})$ of (\ref{4.799}), (\ref{1.6}).

From their parametric analogs in Appendix D of Withers (1987), we obtain (see Appendix B) in the notation of (\ref{3.3.4})
\begin{eqnarray}
T_1 (F) = -|2| /2,
\
T_2 (F) = |3| /3 + \left |2^2 \right | /8 - \sum \lambda^2_a T \left(a^2\right) /2,
\label{4.1}
\end{eqnarray}
and
\begin{eqnarray}
T_3 (F)
&=&
-\sum \lambda^3_a T \left(a^2\right) /2 +   \sum \lambda^3_a T \left(a^3\right) - \sum \lambda^3_a T \left(a^4\right) /4 +
\sum \lambda^3_a T \left(a^2 a^2 \right) /2
\nonumber
\\
&&
+ \sum \lambda_a \lambda^2_b T \left(a^2 b^2 \right) /4 - \sum   \lambda_a \lambda_b^2 T \left(a^2b^3\right) /6 -
\sum \lambda_a   \lambda_b \lambda_c T\left(a^2b^2c^2\right) /48.
\nonumber
\end{eqnarray}
For the one sample case $(k = 1)$, these reduce to
\begin{eqnarray}
&&
T_1 (F) = -T\left(1^2\right) /2,
\label{4.4}
\\
&&
T_2 (F) =  T \left(1^3\right) /3 + T \left(1^21^2\right)/8 - T\left(1^2\right)/2,
\label{4.599}
\\
&&
T_3 (F) = -T\left(1^2\right) /2 + T\left(1^3\right) - T\left(1^4\right) /4 + 3T\left(1^21^2 \right)/4 - T\left(1^21^3\right) /6
\nonumber
\\
&&
\qquad \qquad
- T\left(1^21^21^2\right) /48.
\label{4.699}
\end{eqnarray}

\begin{note}
Let $\{N_i(n), i \geq 0\}$ be given functions satisfying $N_i(n)/n^{-i} \rightarrow 1$.
Then (\ref{1.3}) may be rewritten as $S_{np}(F)+O(n^{-p})$, where
\begin{eqnarray}
S_{np} {(F)} = \sum_{i=0}^{p-1} N_i (n) S_i {(F)}.
\label{4.799}
\end{eqnarray}
So, $S_{np}\widehat{(F)}$ is a $p$th order estimate of $T(F)$.

Suppose now that it is known that there exists an UE and that
it has the form $S_{np}\widehat{(F)}$.
Then this gives a method of obtaining it.
For example, if $k = 1$ and $T(F)$ is a
polynomial of degree $p$ in $F$ (for example, a product of moments or cumulants of total
degree $p$), then the UE of $T(F)$ has the form (\ref{4.799}) with
\begin{eqnarray}
N_i (n) = 1/ (n-1)_i,
\label{1.6}
\end{eqnarray}
where $(r)_i = r!/(r-i)! = r(r-1) \cdots (r-i+1)$.
In this case, $\{ S_i \}$ are given in terms of $\{ T_i \}$ by equation (2.17.2) of Withers (1987):
\begin{eqnarray*}
S_0 = T,
\
S_1 = T_1,
\
S_2 = T_2 - T_1,
\
S_3 = T_3 - 3T_2 + 2T_1,
\ldots .
\end{eqnarray*}
If $k=1$ and we choose $N_i(n)$ as in (\ref{1.6}), then $S_i$ is generally a simpler expression than $T_i$:
\begin{eqnarray}
&&
S_0 (F) = T(F),
\
S_1 (F) = -T \left(1^2\right) /2,
\nonumber
\\
&&
S_2 (F) = T\left(1^3\right) /3 + T\left(1^21^2\right) /8,
\label{4.10}
\\
&&
S_3 (F) = -T \left(1^4\right) /4 + 3T\left(1^21^2\right)/8 - T\left(1^21^3\right)/6 - T\left(1^21^21^2\right)/48.
\label{4.11}
\end{eqnarray}
If $k \neq 1$,
\begin{eqnarray}
&&
S_0(F) = T(F),
\
S_1(F) = T_1(F) \mbox{ of (\ref{4.4})},
\nonumber
\\
&&
S_2(F) = T_2(F) - T_1(F) = |3|/3 + \left | 2^2 \right |/8  + \sum \left(\lambda_a - \lambda^2_a\right) T\left(a^2\right)/2,
\nonumber
\end{eqnarray}
and so on.
\end{note}

\begin{note}
For $p \geq 1$, set $e_{n,p} (T,F) = T_{np} (F)$ of (\ref{1.3}) and let
$\{U_i (F)\}$ be smooth.
Then a $p$th order estimate of
\begin{eqnarray*}
U_n (F) = \sum_{i=0}^{\infty} n^{-i} U_i (F)
\end{eqnarray*}
is
\begin{eqnarray}
U^{\star}_{(n)p} \widehat{(F)}=\sum_{i=0}^{p-1} n^{-i}e_{n,p-i} \left( U_{i}, \widehat{F} \right).
\label{1.8}
\end{eqnarray}
Let $\kappa_r ({\bf X})$ denote any $r$th order cumulant of ${\bf X}$, any $q \times 1$ random vector.
Then $n^{1-r} \kappa_r (T\widehat{(F)})$ can
be expanded in the form (\ref{1.8}); a method of obtaining $\{U_i\}$ is
illustrated in Section 6 for the case $r=2$.
\end{note}

\begin{note}
$ET \widehat{(F)}$ may be infinite or may not exist.
For example, this is
the case if $k = s = 1$, $T (F) = \mu (F)^{-I}$, $I \geq 1$ and $F$ has positive density at zero,
or $\dot{F} (x)$ approaches zero
too slowly as $ x \rightarrow  0$.
So, Quenouille (1956, page 356) is
wrong in giving $ \overline{X}^{-1}$ finite bias for $X \sim N$ (2,1).
In such cases, our
method may be salvaged provided we know an upper bound for $|T(F) |$,
say $|T(F)|< u < \infty $.
By large deviation theory $P (|T \widehat{(F)}| \geq  u) = O(\exp (-n \lambda))$,  where $\lambda  > 0$.
Typically, $ \widetilde{T}_{np} \widehat{(F)}$ is a $p$th order estimate of $T(F)$, where
\begin{eqnarray}
\widetilde{T}_{np} (F) = \left \{
\begin{array}{ll}
T_{np} (F), & \mbox{ if } |T (F)| < u,\\
c,  & \mbox{ otherwise }
\end{array} \right.
\label{1.11}
\end{eqnarray}
and $c$ is any finite constant, for example, $u$.

The estimates (\ref{4.799}) and (\ref{1.8}) can be adapted similarly, to give
$\widetilde{S}_{np} \widehat{(F)} $ and $\widetilde{U}_{np}^{\star} \widehat{(F)}$ say.
Similarly, if $U_{(n)}(F)$ is the formal expansion of $ n^{r-1} \kappa_r (T_{np} \widehat{(F)})$ then
\begin{eqnarray*}
U_{nq}^{\star} \widehat{(F)} I \left( \left| T \widehat{(F)} \right | < u \right)
\mbox{ is a $q$th order estimate of }
n^{r-1} \kappa_r \left( \widetilde{T}_{np} \widehat{(F)} \right)
\end{eqnarray*}
even if $\kappa_r (T \widehat{(F)} )$ is not finite.
For example, the variances in equations (10.17)-(10.20) of Kendall and Stuart (1977) are infinite if the
density at zero is positive.
\end{note}

\begin{note}
An alternative estimate of bias $O(n^{-p} )$ is
$T_{np}^{+} \widehat{(F)} = T_{nq} \widehat{(F)}$,
where $q  \leq p $ is the maximum integer such
that $ \{ n^{-i} T_i \widehat{(F)}, 0 \leq i \leq q \}$ decreases in
absolute value.
This may be useful if $T_{np} \widehat{(F)} $ diverges.
Note that $S_{np}^{+} (F)$ and $\widetilde{T}_{np}^{+} \widehat{(F)}$ may be defined
analogously from (\ref{4.799}) and (\ref{1.11}).
\end{note}

\section{Examples}
\setcounter{equation}{0}

\begin{example}
Suppose $k=1$, ${\bf X} \sim F$ on $R^s$ and $T(F)=g(\mu)$, where ${\bm \mu} = {\bm \mu} (F) =E {\bf X}$ has dimension $s_1=s$
and $g$ is a function with finite derivatives at ${\bm \mu}$.
By the chain rule (\ref{A6}) or (\ref{A7}) of Appendix A,
\begin{eqnarray*}
T_{F} \left( {\bf x}_1, \ldots, {\bf x}_r\right)= g_{j_1 \cdots j_r} \mu_{j_1 {\bf x}_1} \cdots  \mu_{j_r {\bf x}_r},
\end{eqnarray*}
where
\begin{eqnarray*}
\mu_{j {\bf x}} = \mu_{j F} ({\bf x}) = x_j - \mu_j,
\
g\cdots = g \cdots ({\bm \mu})
\end{eqnarray*}
are the partial derivatives of $g({\bm \mu})$ with respect to ${\bm \mu}$, and summation
of the repeated indices $j_1, \ldots, j_r$ over their range $1, \ldots, s$ is implicit.
So,
\begin{eqnarray*}
T \left( 1^{i_1}1^{i_2} \cdots \right) = g_{j_1 \cdots j_{i_1} k_1 \cdots k_{i_2} \cdots} \mu \left[ j_1 \cdots j_{i_1} \right]
\mu \left[ k_1 \cdots k_{i_2} \right] \cdots,
\end{eqnarray*}
where
\begin{eqnarray}
\mu \left[j_1 \cdots j_a \right] = \int \left( x_{j_1} - \mu_{j_1} \right) \cdots \left( x_{j_a} - \mu_{j_a} \right) dF({\bf x}),
\label{5.199}
\end{eqnarray}
the joint central moment.
So,
\begin{eqnarray*}
&&
T\left(1^2\right) = g_{ij} \mu [ij]=\sum_{i=1}^{s} g_{ii} \mu [ii]+  2 \sum_{1 \leq i < j \leq s} g_{ij} \mu [ij],
\\
&&
T \left(1^3\right) = g_{ijk}\mu[ijk],
\\
&&
T \left(1^4\right) = g_{ijkl} \mu [ijkl],
\\
&&
T\left(1^21^2\right) = {g_{j_1 j_2 k_1 k_2}} \mu \left[j_1j_2\right] \mu\left[k_1k_2\right],
\\
&&
T\left(1^21^3 \right) = g_{ijklm} \mu [ij] \mu [klm],
\\
&&
T\left(1^21^21^2\right) = g_{ijklmn} \mu [ij] \mu [kl] \mu [mn].
\end{eqnarray*}
So, by (\ref{4.4})-(\ref{4.699})
\begin{eqnarray*}
&&
T_1 (F) = -C_1 = -g_{ij} \mu [ij] /2,
\\
&&
T_2(F) = -g_{ij} \mu[ij] /2 + g_{ijk} \mu [ijk] /3 + g_{ijkl}\mu [ij] \mu [kl] /8,
\\
&&
T_3(F)  = -g_{ij} \mu[ij]/2 +g_{ijk} \mu[ijk] - g_{ijkl} \left\{\mu[ijkl] - 3 \mu[ij]\mu[kl] \right\}/4
\\
&&
\qquad \qquad
-g_{ijklm} \mu[ij] \mu [klm] /6 - g_{ijklmn} \mu[ij] \mu [kl] \mu [mn]/48.
\end{eqnarray*}
A $p$th order estimate of $T(F)$ is now given in terms of these by $T_{np}(\widehat {F})$ of (\ref{1.3}).
\end{example}

\begin{example}
Consider Example 5.1 with $g({\bm \mu}) = {\bm \alpha}^{\prime} {\bm \mu} / {\bm \beta}^{\prime} {\bm \mu} = N/D$, say,
where ${\bm \alpha}$, ${\bm \beta}$ are given $s$-vectors.
Its $i$th order partial derivative with respect to ${\bm \mu}$ is
\begin{eqnarray}
g_{j_{1} \cdots j_i} = (-1)^{i-1} (i-1)! D^{-i} \sum^i
\delta_{j_1} \beta_{j_2} \cdots \beta_{j_i},
\label{5.299}
\end{eqnarray}
where
\begin{eqnarray}
\delta_i = \alpha_i - \beta_i T(F)
\label{5.399}
\end{eqnarray}
and
\begin{eqnarray*}
\sum^m f_{i_1\cdots i_m} = f_{i_1\cdots i_m} +  f_{i_2\cdots i_mi_1} + \cdots + f_{i_mi_1\cdots i_{m-1}}.
\end{eqnarray*}
So,
\begin{eqnarray*}
&&
T\left(1^i\right) = (-1)^{i-1} i!D^{-i} \delta_{j_1} \beta_{j_2} \cdots   \beta_{j_i} \mu \left[ j_1 \cdots j_i \right],
\\
&&
T\left(1^2 1^2\right) = -4! D^{-4} \delta_{j_1} \beta_{j_2} \beta_{j_3}   \beta_{j_4} \mu \left[j_1 j_2\right] \mu \left[j_3 j_4\right],
\\
&&
T\left(1^21^3\right) = 4! D^{-5} \left\{ 2 \delta_{j_1} / \beta_{j_1}  + 3 \delta_{j_3} / \beta_{j_3} \right\} \beta_{j_1} \cdots
\beta_{j_5} \mu \left[j_1j_2\right] \mu \left[j_3 j_4 j_5\right],
\\
&&
T\left(1^21^21^2 \right) = -6!D^{-6} \delta_{j_1} \beta_{j_2} \cdots   \beta_{j_6} \mu \left[j_1j_2\right] \mu \left[j_3j_4\right] \mu \left[j_5j_6\right].
\end{eqnarray*}
In particular, for $g ({\bm \mu}) = \mu_1 / \mu_2$ (the ratio of means
for one bivariate sample),
\begin{eqnarray*}
&&
T\left(1^i \right) = (-1)^{i-1} i ! \mu_2^{-i} \left\{ \mu \left[12^{i-1}\right] - T(F) \mu \left[2^i\right] \right\},
\\
&&
T\left(1^21^2\right) = -4! \mu_2^{-4} \left\{ \mu [12] \mu\left[2^2\right] - T(F) \mu\left[2^2\right]^2 \right\},
\\
&&
T\left(1^21^3\right) = 4! \mu_2^{-5} \left\{ 2 \mu[12] \mu \left[2^3\right] + 3 \mu \left[2^2\right] \mu \left[12^2\right] - 5T(F) \mu \left[2^2\right] \mu \left[2^3\right] \right\},
\\
&&
T\left(1^21^21^2\right) = -6! \mu_2^{-6} \left\{ \mu [12] - T(F) \mu \left[2^2\right] \right\} \mu \left[2^2\right]^2,
\end{eqnarray*}
so
\begin{eqnarray*}
&&
S_1(F) = T_1(F) = -C_1 = \mu_2^{-2} \left\{ \mu [12] - T(F) \mu \left[2^2\right] \right\},
\\
&&
T_2(F) = 2 \mu_2^{-3} \left\{ \mu \left[12^2\right] - T(F) \mu \left[2^3\right] \right\} - T_1 (F) \left\{1 + 3 \mu_2^{-2} \mu \left[2^2\right] \right\},
\end{eqnarray*}
$S_2 (F)$ is the same as $T_2 (F)$ with `1 $+$' deleted,
\begin{eqnarray*}
T_3(F)
&=&
\mu_2^{-2} \left\{\mu [12] - T(F) \mu \left[2^2\right] \right\} \left\{ 1-18 \mu_2^{-2} \mu \left[2^2\right] -
8 \mu_2^{-3}  \mu \left[2^3\right] + 15 \mu_2^{-4} \mu \left[2^2\right]^2 \right\}
\\
&&
+6 \mu_2^{-3}\left\{ \mu \left[12^2\right] - T(F)\mu \left[2^3\right] \right\}
\left\{ 1 - 2\mu_2^{-2} \mu \left[2^2\right] \right\}
\\
&&
+6 \mu_2^{-4} \left\{ \mu \left[12^3\right] - T (F) \mu \left[2^4\right] \right\}
\end{eqnarray*}
and
\begin{eqnarray*}
S_3 (F)
&=&
\mu_2^{-2} \left\{ \mu[12] - T(F) \mu \left[2^2\right] \right\} \left\{ -9 \mu_2^{-2} \mu \left[2^2\right] -
8\mu_2^{-3} \mu \left[2^3\right] +15 \mu_2^{-4} \mu\left[2^2\right]^2 \right\}
\\
&&
-12 \mu_2^{-5} \left\{ \mu \left[12^2\right] - T(F) \mu \left[2^3\right] \right\} \mu \left[2^2\right] +
6 \mu_2^{-4} \left\{ \mu \left[12^3\right] - T(F) \mu \left[2^4\right] \right\}.
\end{eqnarray*}
\end{example}

\noindent
[Table 5.1 about here.]

\begin{example}
Consider Example 5.1 with $g({\bm \mu}) = ({\bm \alpha}^{\prime} {\bm \mu})^{p}= N^{p}$, say, where ${\bm \alpha}$ is a given $s$-vector.
The $i$th order partial derivative of $g({\bm \mu})$ with respect to ${\bm \mu}$ is
\begin{eqnarray*}
g_{j_1 \cdots j_i} = (p)_i N^{p-i} \alpha_{j_1} \cdots \alpha_{j_i}.
\end{eqnarray*}
Set
\begin{eqnarray*}
\alpha_{(i)} = N^{-i} \alpha_{j_1} \cdots \alpha_{j_i} \mu \left[j_1 \cdots j_i\right].
\end{eqnarray*}
Then
\begin{eqnarray*}
&&
T\left(1^i\right) = (p)_i N^{p} \alpha_{(i)},
\\
&&
T \left(1^2 1^2\right) = (p)_4 N^{p} \alpha^2_{(2)},
\\
&&
T \left(1^21^3\right) = (p)_5 N^{p} \alpha_{(2)} \alpha_{(3)},
\\
&&
T \left(1^21^21^2\right) = (p)_6 N^{p} \alpha^3_{(2)},
\\
&&
T_1 (F) = -C_1 = -(p)_2 N^{p} \alpha_{(2)}/ 2,
\\
&&
T_2 (F) = N^{p} \left\{ -(p)_2\alpha_{(2)}/2+(p)_3\alpha_{(3)}/3  + (p)_4\alpha^2_{(2)}/8\right\},
\\
&&
T_3 (F) = N^{p} \bigg\{ -(p)_2 \alpha_{(2)} /2 + (p)_3   \alpha_{(3)} - (p)_4 \left[  \alpha_{(4)} - 3 \alpha_{(2)}^2 \right] /4
\\
&&
\qquad \qquad
-(p)_5 \alpha_{(2)} \alpha_{(3)}/6-(p)_6\alpha^3_{(2)}/48  \bigg\}.
\end{eqnarray*}

In particular, for a univariate sample $(s = 1)$ with central moments $\{\mu_r\}$ and $g (\mu) = \mu^{p}$,
\begin{eqnarray}
S_1(F)
&=&
T_1 (F) = -(p)_2 \mu^{p-2} \mu_2 /2,
\nonumber
\\
T_2(F)
&=&
-(p)_2 \mu^{p-2} \mu_2 /2 + S_2 (F),
\nonumber
\\
S_2 (F)
&=&
(p)_3 \mu^{p-3} \mu_3 /3 + (p)_4 \mu^{p-4} \mu_2^2 /8,
\nonumber
\\
T_3(F)
&=&
-(p)_2 \mu^{p-2} \mu_2 /2 + (p)_3 \mu^{p-3} \mu_3 - (p)_4 \mu^{p-4} \left( \mu_4 - 3 \mu_2^2 \right) /4
\nonumber
\\
&&
-(p)_5 \mu^{p-5} \mu_3 \mu_2 /6 - (p)_6 \mu^{p-6} \mu_2^3 /48
\nonumber
\end{eqnarray}
and
\begin{eqnarray}
S_3(F) = -(p)_4 \mu^{p-4} \left(2 \mu_4 - 3 \mu_2^2\right) /8 - (p)_5 \mu^{p-5} \mu_3   \mu_2 /6 - (p)_6 \mu^{p-6} \mu_2^3 /48.
\nonumber
\end{eqnarray}
In particular, for $p$ a positive integer, by Note 4.1, an UE for $\mu^p$ is
\begin{eqnarray*}
\sum_{i=0}^{p-1} S_i \left(\widehat{F}\right)/(n-1)_i,
\end{eqnarray*}
where $S_0(F)=\mu^p$, and
\begin{eqnarray*}
&&
\mbox{for } p=2: S_1(F)=-\mu_2,
\\
&&
\mbox{for } p=3: S_1(F)=-3\mu\mu_2, S_2(F)=2\mu_3,
\\
&&
\mbox{for } p=4: S_1(F)=-6\mu^2\mu_2, S_2(F)=8\mu\mu_3+3\mu_2^2, S_3(F)=-6\mu_4+9\mu_2^2.
\end{eqnarray*}
These results may be checked by by solving the
system of equations given by Wishart (1952, page 5).
For $p=4$ the system has seven equations.
Alternatively, one may follow the method of Section 12.22 of Stuart and
Ord (1987) using their tables of the symmetric functions.
For example, after some labor one obtains for $p=4$ the UE $T_n(\widehat{F})$, where
\begin{eqnarray*}
(n-1)_3T_n(F)
&=&
\left(N^3-8n^2+23n-30\right)m_4 - n\left(n^2-7n+4\right)m_3m_1
\\
&&
-n\left(n^2-6n+6\right)m_2^2+n^2(n-9)m_2m_1^2+n^3m_1^4,
\end{eqnarray*}
where $m_i=EX^i$.
Clearly, our method gives a much simpler form.

For $ p=-1$, that is $T(F) = \mu^{-1}$, the above gives
\begin{eqnarray*}
S_{np} (F) = \sum_{i=0}^{p-1}{S_i(F)}/{(n-1)_i},
\end{eqnarray*}
where
\begin{eqnarray*}
&&
S_0 (F) = \mu^{-1},
\
S_1 (F) = -\mu^{-3} \mu_2,
\\
&&
S_2 (F) = -2 \mu^{-4} \mu_3 + 3 \mu^{-5} \mu_2^2,
\\
&&
S_3 (F) = -3 \mu^{-5} \left(2 \mu_4 -3 \mu_2^2\right) + 20 \mu^{-6} \mu_3 \mu_2 -15 \mu^{-7} \mu_2^3,
\end{eqnarray*}
so setting $\gamma_r = \mu_r \mu^{-r}$, $s_i = S_i(F)/T(F)$ is given by
\begin{eqnarray*}
&&
s_1 = -\gamma_2,
\\
&&
s_2 = -2\gamma_3 + 3\gamma_2^2,
\\
&&
s_3 = -3 \left(2\gamma_4-3\gamma_2^2\right) + 20\gamma_3\gamma_2 - 15\gamma_2^3.
\end{eqnarray*}
Some simulations estimating the
bias of $\widetilde{S}_{ni} \widehat{(F)}$ of (\ref{4.799}), (\ref{1.6}) and Note 4.3 with
$c=1/u=\mu/10$ for $ 1 \leq i \leq 4$, for $\mu^{-1}$, are given in Table 5.1.
For $n = 100$  and $p=2$, the estimates are poor: see Appendix C.
\end{example}

Example 5.1  estimated a smooth function of the mean of one multivariate
distribution.
We now estimate a smooth function of the means of $k$
univariate distributions.

\begin{example}
Suppose we have $k$ univariate samples (that is $s_1 = \cdots = s_k =  1)$ with $T (F) = g({\bm \mu})$,
where now ${\bm \mu} = ( \mu(F_1), \ldots, \mu(F_k) )$.
That is, $T(F)$ is a function of the means of $k$ univariate samples.
Then
\begin{eqnarray*}
T_{F}  \left( \begin{array}{ccc} a_1 \cdots a_r \\ x_1 \cdots x_r \end{array} \right ) = g_{a_1 \cdots a_r}
\mu_{a_1 x_1} \cdots \mu_{a_r x_r},
\end{eqnarray*}
where $g \cdots$ is the partial derivative with respect to ${\bm \mu}$ and
\begin{eqnarray*}
\mu_{ax} = \mu_{F_a}(x) = x - \mu \left(F_a\right) = x - \mu_a.
\end{eqnarray*}
So,
\begin{eqnarray*}
T \left(a^i b^j \cdots \right) = {g_{a^i b^j \cdots}} \mu_i [a]\mu_j [b] \cdots,
\end{eqnarray*}
where
\begin{eqnarray*}
\mu_i [a] = \mu_i \left(F_a\right) = \int \left(x-\mu_a\right)^i dF_a (x),
\end{eqnarray*}
the $i$th central moment of $F_a$.
So, for $\lambda_a$ of (\ref{3.2}),
\begin{eqnarray*}
C_1
&=&
\sum_a \lambda_a g_{aa} \mu_2 [a] /2,
\\
C_2
&=&
\sum_a \lambda_a^2 g_{aaa} \mu_3 [a] /6+
\sum_{ab}\lambda_a\lambda_bg_{aabb} \mu_2 [a] \mu_2 [b] /8,
\\
C_3
&=&
\sum \lambda_a^3 g_{aaaa} \left\{ \mu_4 [a] - 3 \mu_2[a]^2 \right\} /24
\\
&&
+ \sum \lambda_a \lambda_b^2 g_{aabbb}\mu_2 [a] \mu_3 [b] /12 + \sum \lambda_a \lambda_b \lambda_c g_{aabbcc} \mu[a]\mu_2 [b] \mu_2 [c] /48,
\\
T_1 (F)
&=&
-C_1,
\\
T_2(F)
&=&
\sum \lambda_a^2 g_{aaa} \mu_3 [a] /3 + \sum \lambda_a \lambda_b g_{aabb} \mu_1[a]\mu_2[b] /8 - \sum \lambda_a^2 g_{aa} \mu_2 [a] /2,
\\
T_3 (F)
&=&
-\sum \lambda_a^3 g_{aa} \mu_2 [a]/2 + \sum \lambda_a^3 g_{aaa} \mu_3 [a] -
\sum \lambda_a^3 g_{aaaa} \left\{ \mu_4 [a] /4+\mu_2 [a]^2 /2\right\}
\\
&&
+ \sum \lambda_a^2 \lambda_b g_{aabb} \mu_2 [a]\mu_2[b]/4  - \sum \lambda_a \lambda^2_b g_{aabbb}\mu_2 [a] \mu_3 [b] /6
\\
&&
- \sum \lambda_a \lambda_b \lambda_c g_{aabbcc} \mu_2[a]\mu_2[b] \mu_2 [c] /48.
\end{eqnarray*}
\end{example}

\begin{example}
Consider Example 5.4 with  $g({\bm \mu}) = {\bm \alpha}^{\prime} {\bm \mu}/{\bm \beta}^{\prime} {\bm \mu} = N/D$, say,
where ${\bm \alpha}$ and ${\bm \beta}$ are given $k$-vectors.
Set
\begin{eqnarray*}
&&
\gamma_a = \alpha_a /\beta_a - T(F),
\\
&&
A_{ikl} = D^{-kl} \sum_a \lambda_a^{i+kl-1} \beta_a^k\mu_k (a)^l \gamma_a,
\\
&&
B_{ikl} = \left\{ A_{ikl} \right\} \mbox{ at } \gamma_a \equiv 1,
\\
&&
A_k = A_{0k1},
\\
&&
B_k = B_{0k1}.
\end{eqnarray*}
Then, by (\ref{5.299}),
\begin{eqnarray*}
&&
C_1 = -A_2,
\\
&&
C_2 = A_3 - 6A_2 B_2,
\\
&&
C_3 = -A_4 + 3A_{022} + 6A_2 B_3 + 9 A_3 B_2 - 15A_2 B_2^2,
\\
&&
T_1(F) = A_2,
\\
&&
T_2(F) = 2A_3 - 3A_2 B_2 + A_{121},
\\
&&
T_3 (F) = A_{221}-9A_{131}-3A_3+6A_4-12A_{022}-3A_{121} B_2-3A_2 B_{121} - 8A_2B_3
\\
&&
\qquad \qquad
-12A_3 B_2 + 15A_2 B_2^2.
\end{eqnarray*}
In particular, for $g ({\bm \mu}) = \mu_1 / \mu_2$ (the ratio of means
for two univariate samples), setting $\nu_k = \mu_2^{-k} \mu_k[2]$,  we obtain
\begin{eqnarray*}
&&
C_1 = \lambda_2 \nu_2 \mu_1 / \mu_2,
\\
&&
C_2 = \lambda_2^2 \left(-\nu_3 + 6 \nu_2^2 \right) \mu_1 / \mu_2,
\\
&&
C_3 = \lambda_2^3 \left( \nu_4 - 3 \nu_2^2 - 15 \nu_2 \nu_3 + 15 \nu_2^3\right) \mu_1 / \mu_2,
\\
&&
T_1 (F) = -\lambda_2 \nu_2 \mu_1 / \mu_2,
\\
&&
T_2 (F) = \lambda_2^2 \left(-2 \nu_3 - \nu_2 + 3 \nu_2^2\right) \mu_1/\mu_2,
\\
&&
T_3 (F) = \lambda_2^3 \left(-6 \nu_4  -6 \nu_3 - \nu_2 - 15 \nu_2^3 + 20 \nu_3 \nu_2 + 18 \nu_2^2\right) \mu_1 / \mu_2.
\end{eqnarray*}
This may also be derived from (\ref{5.299}).
\end{example}

Central moments and functions of them may be viewed as functions of noncentral
moments and so dealt with using Examples 5.1 and 5.4.
However, it is much more
convenient to deal with them directly in terms of the derivatives of the
central moments.
We now give these.

\begin{example}
One univariate sample ( that is  $ k = s_1 = 1$) with $T(F) = \mu_r(F) = \mu_r$, the $r$th central moment of $X \sim F$.
Let $\mu = \mu (F)$ denote the mean of $F$.
Recall that $(r)_i = r!/(r - i)!$  and set $h_i = \mu_{x_i} = x_i - \mu$.
The general derivative of $\mu_r (F)$ is
\begin{eqnarray}
T_{x_1\cdots x_p} = \mu_{rF} \left(x_1 \cdots x_{p} \right) = (-1)^p\left\{ (r)_p\mu_{r-p} - (r)_{p-1} \sum_{i=1}^p
\left(h^{r-p}_i-\mu_{r-p+1} h_i^{-1}\right) \right\} \prod_{j=1}^{p} h_j.
\label{5.699}
\end{eqnarray}
For example,
\begin{eqnarray*}
&&
T_x = -r \mu_{r-1} \mu_x + \mu_x^r - \mu_r,
\\
&&
T_{xy} = (r)_2 \mu_{r-2} \mu_x \mu_y - r \sum_{xy}^2 \left( \mu_x^{r-1} - \mu_{r-1}\right) \mu_y,
\\
&&
T_{xyz} = -(r)_3 \mu_{r-3} \mu_x \mu_y \mu_z + (r)_2 \sum_{xyz}^3 \left(\mu_x^{r-2} - \mu_{r-2}\right) \mu_y \mu_z.
\end{eqnarray*}
These basic building blocks are written out more explicitly
up to $r=6$ in Appendix D.
Setting $ q = i_1 + i_2 + \cdots$, this gives
\begin{eqnarray}
\mu_r \left(1^{i_1}1^{i_2} \cdots\right)
&=&
\displaystyle
(-1)^{ {q}} \bigg[ (r)_{{q}}\mu_{r - {q} }\prod_{j=1}  \mu_{i_j}
\nonumber
\\
&&
\displaystyle
-(r)_{{q}-1}\sum_{I=1} i_I \left( \mu_{r - {q}+i_I} - \mu_{r -  {q}+1}  \mu_{i_I - 1}\right) \prod_{j \neq I}   \mu_{i_j} \bigg]
\label{5.799}
\\
&=&
\left\{ \begin{array}{ll}
\displaystyle
0, & \mbox{ if }  {q} > r,\\
\displaystyle
(-1)^{r-1} (r - 1)! \prod_{j=1} \mu_{i_j}, & \mbox{ if }  {q} = r.
\end{array}
\right.
\nonumber
\end{eqnarray}
For example,
\begin{eqnarray}
&&
\mu_r (1^2) = (r)_2 \mu_{r-2}\mu_2 - 2r\mu_r,
\label{4.4.1}
\\
&&
\mu_r \left(1^3\right) = -(r)_3 \mu_{r-3} \mu_3 + 3 (r)_2 \left(\mu_r - \mu_{r-2} \mu_2\right),
\label{4.4.2}
\\
&&
\mu_r \left(1^4\right) = (r)_4 \mu_{r-4} \mu_4 -4 (r)_3 \left(\mu_r - \mu_{r-3} \mu_3 \right),
\label{4.4.3}
\\
&&
\mu_r \left(1^21^2\right) = (r)_4\mu_{r-4} \mu_2^2 - 4 (r)_3 \mu_{r-2} \mu_2,
\label{4.4.4}
\\
&&
\mu_r \left(1^21^3\right) = -(r)_5 \mu_{r-5} \mu_3\mu_2 + (r)_4 \left( 2 \mu_{r-3}\mu_3 + 3 \mu_{r-2} \mu_2 - 3 \mu_{r-4}\mu_2^2 \right),
\label{4.4.5}
\\
&&
\mu_r \left(1^21^21^2\right) = (r)_6 \mu_{r-6}\mu_2^3 - 6 (r)_5 \mu_{r-4} \mu_2^2.
\label{4.4.6}
\end{eqnarray}
Substituting into the expressions of (\ref{3.3.1})-(\ref{3.3.3}) for the coefficient $C_i$
of $n^{-i} $ in the expansion of $E \mu_r \widehat{(F)}$ gives
\begin{eqnarray}
T_1 (F)
&=&
-C_1 = r \mu_r-(r)_2\mu_{r-2}\mu_2/2,
\label{5.999}
\\
C_2
&=&
(r)_2\mu_r/2-(r)_2 (r-1) \mu_{r-2}\mu_2/2 - (r)_3\mu_{r-3}\mu_3 /6 + (r)_4\mu_{r-4} \mu_2^2/8,
\label{5.1099}
\\
C_3
&=&
-(r)_3\mu_r /6 +\mu_{r-2}\mu_2(r)_3(r-1)/4+(r)_3 (r-2)\mu_{r-3}\mu_3/6
\nonumber
\\
&&
+(r)_4\mu_{r-4} \left( \mu_4-3(r-1)\mu_2^2 \right)/24-(r)_5\mu_{r-5}\mu_3\mu_2/12
\nonumber
\\
&&
+(r)_6\mu_{r-6}\mu_2^3/48,
\nonumber
\\
C_4
&=&
(r)_4\mu_r/24-(r)_4(r-7)\mu_{r-2}\mu_2/12-(r)_6\mu_{r-3}\mu_3/2
\nonumber
\\
&&
+\mu_{r-4}\left\{ -(r)_4(r-3) \ \mu_4/24+(r)_4(r^2-3r-8)\mu_2^2 /16 \right\}
\nonumber
\\
&&
+\mu_{r-5}\left\{ -(r)_5\mu_5/120 + (r)_6(r-2)\mu_3\mu_2/12 \right\}
\nonumber
\\
&&
+ (r)_6\mu_{r-6} \left( \mu_4\mu_2/48 +\mu_3^2/72-r \mu_2^3/48\right) - (r)_7\mu_{r-7}\mu_3\mu_2^2/48
\nonumber
\\
&&
+ (r)_8 \mu_{r-8} \mu_2^4 /384.
\nonumber
\end{eqnarray}
Substituting into the expressions of (\ref{4.599})-(\ref{4.699}) for the
coefficient $T_i \widehat{(F)}$ of $n^{-i}$ in the expansion for the UE of $\mu_r (F)$ gives
\begin{eqnarray*}
T_2(F) = r^2\mu_r - \left( r^3-r\right) \mu_{r-2}\mu_2/2 - (r)_3\mu_{r-3}\mu_3/3  + (r)_4\mu_{r-4}\mu_2^2/8,
\end{eqnarray*}
and
\begin{eqnarray}
T_3(F)
&=&
r^3\mu_r - \left(r^4-r\right)\mu_{r-2} \mu_2/2 -(r)_3(r+3)\mu_{r-3}\mu_3/3
\nonumber
\\
&&
+(r)_ 4\mu_{r-4}\left\{ -2\mu_4+(r+6)\mu_2^2 \right\} / 8 +(r)_5\mu_{r-5}\mu_3\mu_2/6 -(r)_6\mu_{r-6}\mu_2^3/48.
\nonumber
\end{eqnarray}
Similarly, from (\ref{4.10}) and (\ref{4.11}),
\begin{eqnarray*}
S_2(F) = (r)_2\mu_r - r^2(r-1)\mu_{r-2} \mu_2/2 - (r)_3\mu_{r-3} \mu_3/3 + (r)_4\mu_{r-4}\mu^2_2/8
\end{eqnarray*}
and
\begin{eqnarray*}
S_3(F)
&=&
(r)_3\mu_r - r(r)_3\mu_{r-2}\mu_2/2 -r(r)_3\mu_{r-3}\mu_3/3
\\
&&
-(r)_4 \mu_{r-4} \mu_4/4 + (4r-9)(r)_4\mu_{r-4}\mu^2_2/8 +(r)_5\mu_{r-5} \mu_3\mu_2/6 - (r)_6\mu_{r-6}\mu^3_2/48.
\end{eqnarray*}
Now from James (1958, page 6) the UE for $\mu_r$ has the form
\begin{eqnarray}
l_r = \left\{ \sum_{i=0}^{s} a_{ir} \widehat{(F)}n^{-i}\right \} / \prod_{i=1}^{r-1} (1-i/n)
\label{4.11a}
\end{eqnarray}
for $r=2s$ or $2s + 1$,
which can be recovered from $\{T_i, i \leq s \}$ as in Note 4.1.
So, the above $\{T_i, i \leq 3 \}$ provide UEs for $\mu_r$ for $r \leq 7$.
These were given for $r \leq 6$ on James (1958, page 6) and agree with our results.

For example, for $\mu_3$, $T(1^2)=-2\mu_2$, so $S_1(F)=3\mu_3$ and $T(1^3)=12\mu_3$, $T(1^21^2)=0$,
so $ S_2(F)=4 \mu_3$ and so the UE of $\mu_3 $ is
\begin{eqnarray*}
\mu_3\left( \widehat{F}\right) \left\{ 1+3 /(n-1)+4/(n-1)_2 \right\} \ = \mu_3\left(\widehat{F}\right)
\left\{ \left(1-n^{-1}\right)\left(1-2n^{-1}\right)\right\}^{-1}.
\end{eqnarray*}
For $r =7$, we obtain in this way $\{ a_{i7} = a_{i7} (F) \}$ of (\ref{4.11a}) as
\begin{eqnarray}
&&
a_{07} = \mu_7,
\
a_{17} =-7 \left(2 \mu_7 + 3 \mu_5\mu_2\right),
\nonumber
\\
&&
a_{27} = 7 \left(11\mu_7+39\mu_5\mu_2-10\mu_4\mu_3+15\mu_3\mu_2^2\right),
\nonumber
\\
&&
a_{37} =-7 \left(28 \mu_7 +192 \mu_5 \mu_2-80\mu_4\mu_3 +60\mu_3\mu_2^2\right).
\nonumber
\end{eqnarray}
\end{example}

\begin{example}
One univariate sample (that is $k=s_1=1$) with $T(F) = \prod_{j=2}^{\ q} \mu_j^{p_j}$ for $\{p_j\}$ arbitrary and $\{\mu_j \}$ as in Example 5.6.
Set $S_i({\bm \mu}) = \mu_i$ and $g({\bf S}) = \prod S_j^{p_j}$.
The ordinary partial derivatives of $g({\bf S})$ are
\begin{eqnarray*}
&&
g_i = p_i \mu_i^{-1}T(F),
\
g_{ij}=p_i \left( p_j - \delta_{ij} \right) \left(\mu_i\mu_j\right)^{-1}T(F),
\\
&&
g_{ijk} = p_i \left(p_j-\delta_{ij}\right) \left(p_k-\delta_{ik}-\delta_{jk}\right) \left(\mu_i\mu_j\mu_k\right)^{-1}T(F),
\end{eqnarray*}
and so on, where $\delta_{ij} = 1$ if $i=j$ and 0 otherwise.
Set
\begin{eqnarray*}
\left[{}^{a b}_{i j} \cdots \right] = \int \mu_{iF}\left(x^a\right) \mu_{jF} \left(x^b\right) \cdots dF (x).
\end{eqnarray*}
So, $[{}^a_i] = \mu_i(1^a)$ of  (\ref{5.799}) and by (\ref{5.699}), and
\begin{eqnarray*}
\left[ {}^{1 1}_{i j}\right] = ij \mu_{i-1} \mu_{j-1}  \mu_2-\sum_{ij}^2 i \mu_{i-1} \mu_{j+1} + \mu_{i+j}-\mu_i\mu_j,
\end{eqnarray*}
where $\sum^m_{i_1 \cdots i_m} f_{i_1 \cdots i_m}=\sum^m f_{i_1 \cdots i_m}$ is defined in Example 5.2.

By (\ref{A8}),
\begin{eqnarray}
-2T_1 (F) = 2C_1=T\left(1^2\right) = T(F) \left\{2<1,2>+<1,1>+<1^2> \right\},
\label{5.1999}
\end{eqnarray}
where
\begin{eqnarray}
&&
<1,2> = \sum_{i<j} p_ip_j\left[{}^{1 1}_{ij}\right] \mu_i^{-1} \mu_j^{-1},
\nonumber
\\
&&
<1,1> = \sum_i \left(p_i \right)_2 \left[{}^{1 1}_{i i}\right] \mu_i^{-2},
\nonumber
\\
&&
<1^2> = \sum_i p_i \left[ \stackrel2 i \right] \mu_i^{-1}.
\nonumber
\end{eqnarray}
Other terms are calculated similarly.
For example, $C_2$, $T_2 (F)$ and $S_2 (F) $ are given by (\ref{3.6}), (\ref{4.599}), and (\ref{4.10}) in
terms of $T(1^2)$, $T (1^3)$ and $T (1^2 1^2)$.
Also by (\ref{A9}) to (\ref{A11})
\begin{eqnarray}
T\left(1^3\right)
&=&
T(F)\bigg\{ \sum_{ijk} p_i \left(p_j-\delta_{ij}\right) \left(p_k- \delta_{ik}-\delta_{jk}\right)
\left( \mu_i\mu_j\mu_k \right)^{-1}\left[{}^{1 1 1}_{i j k}\right]
\nonumber
\\
&&
+3 \sum_{ij}p_i \left(p_j-\delta_{ij}\right) \left(\mu_i\mu_j\right)^{-1}\left[ {}^{21}_{i j}\right] +
\sum_i p_i\mu_i^{-1} \left[{}^3_i\right] \bigg\},
\label{5.2099}
\end{eqnarray}
and
\begin{eqnarray}
T\left(1^21^2\right)
&=&
T(F)\bigg\{ \sum_{ijkl}p_i\left(p_j-\delta_{ij}\right)\left(p_k- \delta_{ik} - \delta_{jk}\right)
\left(p_{l}-\delta_{il}-\delta_{jl} -\delta_{kl}\right)
\left( \mu_i\mu_j\mu_k\mu_{l} \right)^{-1}
\nonumber
\\
&&
\qquad \qquad
\times
\left [{}^{11}_{kl} \right]
\nonumber
\\
&&
+ \sum_{ijk} p_i\left(p_j-\delta_{ij}\right) \left(p_k-\delta_{i k}-\delta_{jk}\right)
\left( \mu_i\mu_j \mu_k \right)^{-1} G_{ijk}
\nonumber
\\
&&
+ \sum_{ij}p_i \left(p_j-\delta_{ij}\right) \left(\mu_i\mu_j\right)^{-1}
H_{ij} + \sum_i p_i \mu_i^{-1} \mu_i \left(1^21^2\right) \bigg\},
\label{5.21}
\end{eqnarray}
where
\begin{eqnarray}
&&
G_{ijk} =  2 \left[{}^{1 1}_{i j}\right]
\left[{}^2_k\right] +4 \left[ 12_i 1_j 2_k \right],
\nonumber
\\
&&
H_{ij} = 4 \left[1_i12_j^2\right]  + \left[{ }^2_i\right] \left[ {}^2_j\right] + 2 \left[12_i12_j\right],
\nonumber
\\
&&
\left [1^a2^b_i1^c2^d_j \cdots \right ] = \int\int \mu_i \left(x^ay^b\right) \mu_j \left(x^c y^d\right) \cdots dF(x)dF(y),
\nonumber
\end{eqnarray}
so that
\begin{eqnarray*}
&&
\left[1_i^a 1_j^b \cdots \right] = \left[ {}^{a b}_{i j}\cdots \right],
\\
&&
\left[ 12_i 1_j 2_k \right] = (i)_2 \mu_{i-2} A_j A_k - i \sum_{jk}^2 B_{ij} A_k
\end{eqnarray*}
for
\begin{eqnarray*}
A_j = \mu_{j+1} - j \mu_{j-1} \mu_2,
\
B_{ij} = \mu_{i+j-1} - j \mu_{j-1} \mu_i - \mu_{i-1} \mu_j.
\end{eqnarray*}
By  (\ref{5.699}),
\begin{eqnarray*}
\left[ {}^{111}_{ijk} \right ]
&=&
-ijk \mu_{i-1}\mu_{j-1}\mu_{k-1}\mu_3 + \sum^3 ij \mu_{i-1} \mu_{j-1} \left(\mu_{k+2} - \mu_k\mu_2\right)
\\
&&
- \sum^3 i\mu_{i-1} \left( \mu_{j+k+1} - \mu_{j+1}\mu_k-\mu_{k+1}\mu_j \right) + \mu_{i+j+k}  - \sum^3 \mu_i\mu_{j+k} + 2 \mu_i\mu_j\mu_k,
\\
\left[ {}^{2\;\;1}_{i \;\; j} \right]
&=&
-(i)_2 j \mu_{i-2}\mu_{j-1}\mu_3 + (i)_2\mu_{i-2} \left(\mu_{j+2} -\mu_j\mu_2\right)
\\
&&
+ 2ij \mu_{j-1} (\mu_{i+1} - \mu_{i-1}\mu_2)-2i \left( \mu_{i+j}-\mu_i\mu_j-\mu_{i-1}\mu_{j+1} \right),
\\
\left[ 1_i12_j^2 \right]
&=&
(j)_2 \left\{ \left( -3i \mu_{i-1} \mu_{j-1} +\mu_{i+j-2} - \mu_i\mu_{j-2} \right)
\mu_2 +2 \mu_{i+1} \mu_{j-1} \right\}
\\
&&
+ (j)_3 \left(i \mu_{i-1}\mu_{j-3} \mu_2^2 - \mu_{j-3}\mu_{i+1}\mu_2 \right),
\\
\left [12_i12_j \right ]
&=&
(i)_2(j)_2\mu_{i-2}\mu_{j-2}\mu_2^2 - 2 \sum^2 i (j)_2  \mu_i\mu_{j-2}\mu_2
\\
&&
+ 2ij \left( \mu_{i+j-2}\mu_2 - \mu_{i-1}\mu_{j-1}\mu_2 + \mu_i\mu_j \right).
\end{eqnarray*}
Also $[ {}^i_r]$ for $2 \leq i \leq 4$ and $\mu_i (1^21^2)$ are given by (\ref{4.4.1})-(\ref{4.4.6}).
\end{example}

\begin{example}
Consider Example 5.7 with $T(F) = \mu_r^p $.
Then
\begin{eqnarray*}
&&
T \left(1^2\right) / T (F) = p \left[{}^2_r\right] \mu_r^{-1} + (p)_2 \mu_r^{-2} \left[{}^{11}_{rr}\right],
\\
&&
T\left(1^3\right) / T (F) = p \mu_r^{-1} \left[{}^3_r\right]+3 (p)_2 \mu_r^{-2} \left[{}^{2\;\;1}_{r\;\;r}\right] + (p)_3 \mu_r^{-3} \left[{}^{111}_{rrr}\right],
\\
&&
T \left(1^2 1^2\right) /T (F) = p \mu_r^{-1} \mu_r \left(1^2 1^2\right) + (p)_2 \mu_r^{-2} H_{rr}
+ (p)_3 \mu_r^{-3} G_{rrr} + (p)_4 \mu_r^{-4} \left[ {}^{11}_{rr}\right]^2.
\end{eqnarray*}
\end{example}

\begin{example}
Consider Example 5.8 with $ T(F) = \mu_2^p$.
Set $ \beta_r = \mu_r \mu_2^{-r/2} $.
Then
\begin{eqnarray*}
&&
T \left(1^2\right) / T (F) = -2p + (p)_2 \left(\beta_4 - 1 \right),
\\
&&
T \left(1^3\right) / T(F) = -6 (p)_2 \left(\beta_4 -1\right) + (p)_3 \left(\beta_6 -3 \beta_4 +2\right),
\\
&&
T \left(1^2 1^2\right) /T(F) = 12(p)_2 -4(p)_3 \left(\beta_4 -1 +2 \beta_3^2\right) + (p)_4 \left(\beta_4 - 1\right)^2.
\end{eqnarray*}
So,
\begin{eqnarray*}
&&
-T_1 (F) /T (F) = C_1 /T (F) = -p + (p)_2 \left(\beta_4 - 1\right) /2,
\\
&&
C_2 / T (F) = (p)_2 \left(5/2 - \beta_4\right) + (p)_3 \left\{ \beta_6/6 - \beta_4 - \beta_3^2 + 5/6 \right\} + (p)_4 \left(\beta_4 -1\right)^2,
\\
&&
T_2 (F) /T(F) = p + (p)_2 \left(4-5 \beta_4/2\right) + (p)_3 \left(2 \beta_6 -9 \beta_4 + 7-6 \beta_3^2 \right)/6
\\
&&
\qquad \qquad
+ (p)_4 \left(\beta_4 - 1\right)^2 /8 = \sum_{i=1}^r (p)_i A_i \mbox{ say, }
\\
&&
S_2 (F) /T (F) = (p)_2 \left(7/2-2 \beta_4\right) + \sum_{i=3}^4 (p)_i A_i.
\end{eqnarray*}
For $p=2$ this gives $ T (F) = \mu_2^2$,
\begin{eqnarray}
&&
C_1 = \mu_4 -3 \mu_2^2,
\
T_1 (F) = -\mu_4 + 3 \mu_2^2,
\label{5.2299}
\\
&&
C_2 = -2\mu_4+5\mu_2^2,
\
T_2(F)    = -5 \mu_4+10 \mu_2^2,
\
S_2(F)    = -4 \mu_4+7 \mu_2^2.
\label{5.2399}
\end{eqnarray}
Note that $C_1$, $C_2$ agree with $\mu(2^2)$ of Sukhatme (1944, page 368).

The UE of $\mu_2^2$ has the form
\begin{eqnarray*}
l_{22} = \left(\sum_{i=0}^2 a_{i22} \widehat{(F)}n^{-i}\right) / \prod_{i=1}^3 (1-i/n).
\end{eqnarray*}
So, $\{ a_i=a_{i22} (F) \}$ are given by
\begin{eqnarray*}
&&
a_0 = T(F) = \mu_2^2,
\\
&&
a_1 = -6T (F) + T_1 (F) = -\mu_4-3 \mu_2^2,
\\
&&
a_2 = 11 T (F) -6T_1 (F) + T_2(F) = \mu_4 + 3 \mu_2^2.
\end{eqnarray*}
\end{example}

We now present a second method for finding an UE
of $\prod_i \mu_i^{p_i}$.
This method avoids computing $\{T_i(F)\}$,
but derives the UE of the vector
\begin{eqnarray}
{\bf T}(F)' = \left\{ \prod_i \mu_i^{p_i} : \sum p_i = p\right\},
\label{5.2499}
\end{eqnarray}
that is, for all products of a given degree $p$,
directly from their first few coefficients $\{ {\bf C}_i \}$.
Suppose ${\bf T} (F)$ has dimension $d = d_{p}$.
Then
\begin{eqnarray}
{\bf C}_i = {\bf A}_i {\bf T} (F),
\nonumber
\end{eqnarray}
where ${\bf A}_i$ is a $d \times d$ matrix of integers and ${\bf A}_0 = {\bf I}_d$, the identity matrix.
So,
\begin{eqnarray*}
{\bm \alpha} (n)^{-1} \ {\bf T} \widehat{(F)}
\end{eqnarray*}
is the UE of ${\bf T}(F)$, where
\begin{eqnarray*}
{\bm \alpha} (n) = \sum_0^{\infty} {\bf A}_i n^{-i}.
\end{eqnarray*}
But this is known to have the form
\begin{eqnarray}
{\bf T}_n \widehat{(F)} = \widehat {\bm \beta}_n / \prod_{i=1}^{p-1} (1-i/n),
\label{5.2699}
\end{eqnarray}
where
\begin{eqnarray*}
\widehat {\bm \beta}_n = \left\{ \sum_{i=0}^{[p/2]} {\bf B}_i n^{-i}\right \} {\bf T} \widehat{(F)},
\end{eqnarray*}
where ${\bf B}_i$ is a $d \times d$ matrix of integers with ${\bf B}_0 = {\bf I}_d$.
So,
\begin{eqnarray*}
\sum_{i=0}^{[p/2]} {\bf B}_i \varepsilon^i
&=&
\left\{ \prod_{i=1}^{p-1} \left(1-i \varepsilon\right) \right\} {\bm \alpha} \left(\varepsilon^{-1}\right)
\\
&=&
\left\{ 1-D_1 (p) \varepsilon + D_2 (p) \varepsilon^2 - \cdots \right\}
\\
&&
\times
\left\{ {\bf I}_d - {\bf A}_1 \varepsilon + \left(-{\bf A}_2 + {\bf A}_1^2 \right) \varepsilon^2 + \left(-{\bf A}_3 +
{\bf A}_1 {\bf A}_2 + {\bf A}_2{\bf A}_1-{\bf A}_1^3 \right) \varepsilon^3 + \cdots \right\},
\end{eqnarray*}
where $D_1(p) = (p)_2 /2$ and $D_2(p) = (p)_3 (p-1/3)/8$.
So, the UE (\ref{5.2699}) is given in terms of $\{A_i, i \leq p/2 \}$:
\begin{eqnarray*}
&&
{\bf B}_0 = {\bf I}_d,
\\
&&
{\bf B}_1 = -D_1(p){\bf I}_d-{\bf A}_1,
\\
&&
{\bf B}_2 = D_2(p){\bf I}_d+D_1(p){\bf A}_1-{\bf A}_2+{\bf A}_1^2,
\\
&&
{\bf B}_3 = -D_3(p){\bf I}_d-D_2(p){\bf A}_1-D_1(p)(-{\bf A}_2+{\bf A}_1^2) - {\bf A}_3+{\bf A}_1 {\bf A}_2 + {\bf A}_2 {\bf A}_1 -{\bf A}_1^3,
\end{eqnarray*}
and so on.
The method also applies to obtaining an UE for
\begin{eqnarray*}
{\bf T} (F)^\prime = \left\{ \mu_1^{p_1} \prod_{i=2}^{q} \ \mu_i^{p_i} : \ \sum p_i = p \right\},
\end{eqnarray*}
where ${\bm \mu} = {\bm \mu} (F)$.
A third method (for $p \leq 8$) due to Fisher (1929) is given in Section 12 of Stuart and Ord (1987).
Their Tables 11 and 10, pages 554-555 may be used to verify Examples 5.8 to 5.11 after some labor.

\begin{example}
Consider Example 5.7 with ${\bf T}(F) = (\mu_4, \mu_2^2)^{\prime}$.
So, (\ref{5.2499}) holds with $p=4$ and $d = [p/2] = 2$.

By (\ref{5.999}), (\ref{5.1099}), for $\mu_4$, $C_1 = -4 \mu_4 + 6 \mu_2^2$ and  $C_2 =6 \mu_4 - 15 \mu_2^2$,
in agreement with $\mu(4)$ on Sukhatme (1944, page 368).
So, by (\ref{5.2299}), (\ref{5.2399})
\begin{eqnarray*}
{\bf A}_1 = \left( \begin{array}{ccc} -4
      & 6 \\ 1 & -3 \end{array} \right ) \mbox{and } {\bf A}_2 =
            \left ( \begin{array}{ccc} 6
     & -15 \\ -2 & 5
\end{array} \right).
\end{eqnarray*}
So,
\begin{eqnarray*}
{\bf B}_1 = -6 {\bf I}_2 - {\bf A}_1 = \left(
\begin{array}{ccc}-2&-6\\-1&-3 \end{array} \right),
\
{\bf B}_2 = 11 {\bf I}_2 + 6 {\bf A}_1 - {\bf A}_2 + {\bf A}_1^2 = \left(
\begin{array}{ccc}3&9\\ 1&3\end{array}\right).
\end{eqnarray*}
So, UEs of  $\mu_4$ and  $\mu_2^2$ are $\mu_{4n} \widehat{(F)}$ and $ \mu_{22n} \widehat{(F)}$, where
\begin{eqnarray*}
\mu_{4n}(F) = \left\{ \mu_4 + \left(-2 \mu_4 - 6\mu_2^2\right) n^{-1} + \left( 3 \mu_4 + 9\mu_2^2 \right) n^{-2} \right\} / \prod_{i=1}^3 (1-i/n),
\end{eqnarray*}
and
\begin{eqnarray*}
\mu_{22n}(F) = \left\{ \mu_2^2+ \left(-\mu_4-3\mu_2^2\right) n^{-1} +
\left(\mu_4 +3\mu_2^2\right) n^{-2} \right\}/ \prod_{i=1}^3(1-i/n).
\end{eqnarray*}
Table 5.2 gives the relative bias of $S_{np}(\widehat{F})$ as estimated from
two runs of 60,000 simulations for $p \leq 2$ and $F$ normal and exponential.
For $n = 100$ and $p=2$, the estimates are poor: see Example C.3.
For $p=3$ the bias is zero.
\end{example}

\noindent
[Table 5.2 about here.]

\begin{example}
Consider Example 5.7 with ${\bf T} (F) = (\mu_5, \mu_3\mu_2)^{\prime}$.
So, (\ref{5.2499}) holds with $p = 5$ and $d= [p/2] =2$.

By (\ref{5.999}), (\ref{5.1099}) for $\mu_5$, $C_1=-5\mu_5+10\mu_3\mu_2$ and $C_2=10\mu_5-50\mu_3\mu_2$,
in agreement with $\mu(5)$ of Sukhatme (1944, page 368).
By (\ref{5.1999})-(\ref{5.21}), for  $\mu_2\mu_3$,
\begin{eqnarray*}
T{\left(1^2\right)} = 2 \mu_5 - 16\mu_3\mu_2,
\
T\left(1^3\right) = -24\mu_5 + 72\mu_3\mu_2,
\
T\left(1^21^2\right) = 96 \mu_3\mu_2,
\end{eqnarray*}
giving $C_1=\mu_5-8\mu_3\mu_2$ and $C_2=-4\mu_5+ 24 \mu_3\mu_2$.
So,
\begin{eqnarray*}
{\bf A}_1 = \left( \begin{array}{ccc} -5 & 10 \\ 1  &  -8 \end{array}
\right )  \mbox{ and } {\bf A}_2 = \left ( \begin{array}{ccc}10
&  -50 \\ -4 & 24 \end{array}  \right ).
\end{eqnarray*}
So,
\begin{eqnarray*}
{\bf B}_1 = -10 {\bf I}_2 - {\bf A}_1 = \left(
\begin{array}{ccc}-5&-10 \\ -1&-2\end{array}\right )
\end{eqnarray*}
and
\begin{eqnarray*}
{\bf B}_2 = 35 {\bf I}_2 + 10 {\bf A}_1 - {\bf A}_2 + {\bf A}_1^2 =\left(
\begin{array}{ccc} 10&20 \\ 1 & 5 \end{array} \right).
\end{eqnarray*}
That is, UEs of $\mu_5$ and $\mu_3 \mu_2$ are $\mu_{5n} \widehat{(F)}$, and $\mu_{32n} \widehat{(F)}$, where
\begin{eqnarray*}
\mu_{5n}(F) = \left\{ \mu_5+ \left(-5 \mu_5 -10 \mu_3 \mu_2\right) n^{-1} +
\left(10 \mu_5 + 20 \mu_3 \mu_2\right)n^{-2} \right\} / \prod^4_{i=1} (1 - i/n)
\end{eqnarray*}
and
\begin{eqnarray*}
\mu_{32n}(F)= \left\{\mu_3 \mu_2+ \left(-\mu_5-2 \mu_3 \mu_2\right)n^{-1}+\left(\mu_5+ 5 \mu_3\mu_2\right)n^{-2} \right\} / \prod^4_{i=1} (1-i/n).
\end{eqnarray*}
\end{example}

\begin{example}
Suppose $ k=s_1=1$  and $T(F) = g(\mu_2)$.
Set $g^r = g^{(r)} (\mu_2)$, and $\beta_r = \mu_r\mu_2^{-r/2}$.
Then
\begin{eqnarray*}
\mu_{x} = \mu_{F}(x) = x - \mu,
\
\mu_{2x} =\mu_{2F}{(x)} =\mu_{x}^2-\mu_2,
\
\mu_{2xy} = \mu_{2F}  (x,y) = -2 \mu_{x}\mu_{y}
\end{eqnarray*}
by (\ref{5.699}).
By (\ref{A8}),
\begin{eqnarray*}
|2| = T\left(1^2\right) = g^2 \mu_{22} (1,1)+g^1 \mu_2 \left(1^2\right),
\end{eqnarray*}
where
\begin{eqnarray*}
&&
\mu_{22}(1,1) = \int \mu_{2x}^2=\int \mu_{2x}^2 dF(x)= \mu_4-\mu_2^2,
\\
&&
\mu_2 \left(1^2\right) = \int \mu_{2xx}= -2 \mu_2 \mbox{ by (\ref{4.4.1})}.
\end{eqnarray*}
Similarly, by (\ref{A9}) to (\ref{A11}) and (\ref{A14}),
\begin{eqnarray*}
T\left(1^3\right)
&=&
g^3\mu_{222}(1,1,1)+3g^2\mu_{22}(1,1^2)+g^1\mu_2\left(1^3\right),
\\
T\left(1^4\right)
&=&
g^4 \mu_{2222} (1,1,1,1)+6g^3 \mu_{222}\left(1,1,1^2\right)
+g^2\left\{4\mu_{22} \left(1,1^3\right) + 3\mu_{22} \left(1^2,1^2\right) \right\}
\\
&&
+g^1 \mu_2 \left(1^4\right),
\\
T\left(1^21^2\right)
&=&
g^4\mu_{22}(1,1)^2+g^3 \left\{ 2 \mu_{22}(1,1)\mu_2 \left(1^2\right) + 4 \mu_{222} (ab, a, b) \right\} + g^1\mu_2 \left(a^2b^2\right)
\\
&&
+g^2 \left\{4 \mu_{22} \left(a,ab^2\right) + \mu_2 \left(1^2\right)^2 + 2\mu_{22}(ab,ab)\right\} \mbox{ at } a=b=1
\\
&=&
\sum_{i=2}^4 g^i a_i \mbox{ say, }
\\
T\left(1^21^3\right)
&=&
g^3 A_3 +g^4 A_4 +g^5 A_5,
\end{eqnarray*}
and by (\ref{A15})
\begin{eqnarray*}
T\left(1^21^21^2\right) = \sum_{i=3}^6 g^i B_i,
\end{eqnarray*}
where
\begin{eqnarray*}
&&
\mu_{222} (1,1,1) = \int \mu_{2x}^3 = \mu_6-3\mu_4\mu_2 +2\mu_2^3,
\\
&&
\mu_{22} \left(1,1^2\right) = \int \mu_{2x} \mu_{2xx}= -2 \left(\mu_4-\mu_2^2\right),
\\
&&
\mu_2 \left(1^3\right) = \int \mu_{2xxx} = 0,
\\
&&
\mu_{2222} (1,1,1,1) = \int \mu_{2x}^4 = \mu_8 - 4 \mu_6  \mu_2 +6\mu_4 \mu_{2}^2 -3 \mu_2^4,
\\
&&
\mu_{222} \left(1,1,1^2\right) = \int \mu_{2x}^2 \mu_{2xx}=-2 \left(\mu_6-2 \mu_4\mu_2 +\mu_2^3\right),
\\
&&
\mu_{22} \left(1,1^3\right) =  \mu_2 \left(1^4\right)=0,
\\
&&
\mu_{22} \left(1^2,1^2\right) = \int \mu_{2xx}^2 =4\mu_4,
\\
&&
\mu_{22} \left(a,ab^2\right)_{a=b=1} = \int \mu_{2x} \mu_{2xyy} = 0,
\\
&&
\mu_{222} (ab,a,b)_{a=b=1} = \int \int \mu_{2xy} \mu_{2x} \mu_{2y} = -2 \mu_3^2,
\\
&&
\mu_{22}(ab,ab)_{a=b=1} = \int \mu_{2xy}^2 = 4\mu_2^2,
\\
&&
\mu_2 \left(a^2 b^2\right)_{a=b=1} = \int \mu_{2xxyy} =0,
\\
&&
a_2 = 1 2 \mu_2^2,
\
a_3 = -4 \left( \mu_4 \mu_2 - \mu_2^3 + 2\mu_3^2 \right),
\
a_4 = \left(\mu_4 - \mu_2^2\right)^2,
\end{eqnarray*}
and
\begin{eqnarray*}
A_3
&=&
6 \mu_{222} \left(a,ab,b^2\right) + 3 \mu_2\left(a^2\right)\mu_{22}\left(b,b^2\right)+
6\mu_{222}(b,ab,ab) \mbox{ at } a=b=1
\\
&=&
3 \int\int \left\{2\mu_{2x} \mu_{2xy} \mu_{2yy} + \mu_{2y} \mu_{2yy}\mu_{2xx} +2 \mu_{2y} \mu_{2xy}^2 \right\}
\\
&=&
3 \int \int \left\{ 8 \left(\mu_x^2 - \mu_2\right) \mu_{x} \mu_{y}^3 + 12\left(\mu_{y}^2-\mu_2\right) \mu_{x}^2 \mu_{y}^2  \right\}
\\
&=&
12 \left\{ 2 \mu_3^2 +3 \left(\mu_2 \mu_4 -\mu_2^3\right)   \right\}
\\
&=&
12 \mu_2^3 \left\{2 \beta_3^2 +3 \beta_4-3 \right\},
\\
A_4
&=&
\int \int \left\{ \mu_{2xx} \mu_{2y}^3 +6 \mu_{2xy}\mu_{2y}^2 +3 \mu_{2yy}  \mu_{2y} \mu_{2x}^2 \right\}
\\
&=&
-2 \int \int \left\{ \mu_{x}^2 \left(\mu_{y}^2 -\mu_2\right)^3 +
6\mu_{x} \mu_{y} \left( \mu_{x}^2-\mu_2 \right) \left( \mu_{y}^2-\mu_2 \right)^2 +
3 \mu_{y}^2 \left(\mu_{x}^2 - \mu_2\right)^2 \left(\mu_{y}^2-\mu_2\right)\right\}
\\
&=&
-2\left\{ \mu_2 \left(\mu_6-3 \mu_4 \mu_2 + 2\mu_2^3\right) +
6\mu_3 \left( \mu_5- 2\mu_3 \mu_2 \right)+3 \left(\mu_4-\mu_2^2\right)^2\right\}
\\
&=&
-2 \mu_2^4 \left\{ \beta_6 -3 \beta_4+2 +6 \beta_3\left(\beta_5 -2\beta_3\right)+3 \left(\beta_4-1\right)^2 \right\},
\\
A_5
&=&
\int \int \mu_{2x}^2 \mu_{2y}^3 = \int \left(\mu_x^2-\mu_2\right)^2 \int \left(\mu_{y}^2-\mu_2\right)^3
\\
&=&
\left(\mu_4-\mu_2^2\right) \left(\mu_6 -3 \mu_4 \mu_2 +2 \mu_2^3\right)
\\
&=&
\mu_2^5 \left(\beta_4 -1\right) \left(\beta_6 - 3\beta_4 + 2\right),
\\
B_3
&=&
B_3^{ijk} \mbox{ at } \left\{ a=b=c=1, S=\mu\right\}
\\
&=&
\int \int\int \left\{  \mu_{2xx}\mu_{2yy}\mu_{2zz} + 6 \mu_{2xx}\mu_{2yz}^2+8 \mu_{2xy}\mu_{2yz}\mu_{2zx}\right\}
\\
&=&
-120\mu_2^3,
\\
B_4
&=&
B_2^{ijkl} \mbox{ at } \left\{ a=b=c=1, S=\mu_2 \right\}
\\
&=&
3 \int\int\int \left\{ \mu_{2x}^2\mu_{2yy}\mu_{2zz}+2 \mu_{2x}^2\mu_{2yz}^2+
4\mu_{2x}\mu_{2y}\mu_{2xy}\mu_{2zz} +8 \mu_{2x}\mu_{2y}\mu_{2xz}\mu_{2yz}\right\}
\\
&=&
36 \left\{ \left(\mu_4-\mu_2^2\right)\mu_2^2 + 4 \mu_3^2\mu_2\right\}
\\
&=&
36 \mu_2^4 \left\{\beta_4-1 +4 \beta_3^2\right\},
\\
B_5
&=&
3 \int\int\int \left\{ \mu_{2xx}\mu_{2y}^2 + \mu_{2xy}\mu_{2x}\mu_{2y} \right\} \mu_{2z}^2
\\
&=&
-6 \left\{ \mu_2 \left(\mu_4-\mu_2^2\right) + \mu_3^2\right\} \left(\mu_4-\mu_2^2\right)
\\
&=&
-6 \mu_2^5 \left\{\beta_4-1+\beta_3^2\right\} \left(\beta_4-1\right),
\\
B_6
&=&
\int\int\int \mu_{2x}^2 \mu_{2y}^2 \mu_{2z}^2 = \left(\mu_4-\mu_2^2\right)^3=\mu_2^6\left(\beta_4-1\right)^3.
\end{eqnarray*}
So,
\begin{eqnarray*}
C_1
&=&
-g^1\mu_2+g^2\left(\mu_4-\mu_2^2\right)/2,
\\
C_2
&=&
g^2\left(5\mu_2^2/2-\mu_4\right)+g^3\left(\mu_6/6- \mu_3^2 - \mu_4\mu_2 + 5 \mu_2^3/6\right)+g^4\left(\mu_4-\mu_2^2\right)^2/8,
\\
C_3
&=&
g^2\mu_4/2+g^3\left(-\mu_6/2+4\mu_4\mu_2+2\mu_3^2-6\mu_2^3\right)
\\
&&
+g^4\left(\mu_8/24-\mu_6\mu_2/3-\mu_5\mu_3-\mu_4^2/2+5\mu_4\mu_2^2/2 +5\mu_3^2\mu_2-41\mu_2^4/24\right)
\\
&&
+g^5\left(\mu_4-\mu_2^2\right) \left( 2\mu_6-9\mu_4\mu_2-3\mu_3^2+7\mu_2^3\right)/24+g^6\left(\mu_4-\mu_2^2\right)^3/48,
\\
T_1(F)
&=&
S_1(F) = -g^2\left(\mu_4-\mu^2_2\right)/2 + g^1\mu_2,
\\
T_2(F)
&=&
g^4 \left(\mu_4-\mu^2_2\right)^2/8 + g^3 \left(\mu_6/3 - \mu_3^2 - 3 \mu_4\mu_2/2 + 7 \mu^3_2/6\right) + g^2\left(-5\mu_4/2 + 4\mu^2_2\right) +g^1\mu_2,
\\
T_3(F)
&=&
\sum^6_{i=1} g^iT_{3i},
\\
S_2(F)
&=&
g^4 \left(\mu_4-\mu^2_2\right)^2/8 + g^3\left(\mu_6/3 - \mu_3^2 - 3 \mu_4\mu_2/2 + 7 \mu^3_2/6\right) + g^2\left(-2\mu_4 + 7\mu^2_2/2\right),
\\
S_3(F)
&=&
\sum^6_{i=2} g^iS_{3i},
\end{eqnarray*}
where
\begin{eqnarray*}
&&
S_{32} = -3\mu_4 + 9\mu^2_2 /2,
\\
&&
S_{33} = 3\mu_6 - 27 \mu_4 \mu_2/2 -7 \mu_3^2+ 13 \mu^3_2,
\\
&&
S_{34} = -\mu_8/4 + 4\mu_6\mu_2/3 +2\mu_5\mu_3 + 11\mu^2_4 /8- 6\mu_4\mu^2_2 - 7 \mu^2_3 \mu_2     + 85\mu^4_2/24,
\\
&&
S_{35} =  \left(\mu_4-\mu_2^2\right) \left(-4\mu_6+15\mu_4\mu_2+3\mu_3^2 -11\mu_2^3 \right)/24,
\\
&&
S_{36} =  -B_6/48,
\\
&&
T_{31} = \mu_2,
\\
&&
T_{32} = -19\mu_4/2 + 31\mu^2_2 /2,
\\
&&
T_{33} = 4\mu_6 - 18\mu_4\mu_2  + 33 \mu_2^3/2 - 10 \mu_3^2,
\\
&&
T_{34} = -\mu_8/4 + 4\mu_6\mu_2/3 +2\mu_5\mu_3 + 7\mu^2_4 /4  - 27\mu_4\mu^2_2/4 - 7 \mu^2_3 \mu_2     + 47\mu^4_2/12,
\\
&&
T_{35} = \left(\mu_4-\mu_2^2\right) \left(-4\mu_6+15\mu_4\mu_2+3\mu_3^2 -11\mu_2^3 \right)/24,
\\
&&
T_{36} =  -B_6/48.
\end{eqnarray*}
\end{example}

\begin{example}
Consider Example 5.12 with $T(F)=\mu_2^{q}$.
Then
\begin{eqnarray*}
g^i
&=&
(q)_i \mu_2^{q-i},
\\
T\left(1^2\right) / \mu_2^{q}
&=&
(q)_2 \left(\beta_4-1\right) - 2q,
\\
T\left(1^3\right) / \mu_2^q
&=&
(q)_3 \left(\beta_6 -3 \beta_4 +2\right) -6(q)_2 \left(\beta_4 -1\right),
\\
T\left(1^4\right) / \mu_2^{q}
&=&
(q)_4 \left(\beta_8-4\beta_6+6\beta_4-3\right)-12(q)_3\left(\beta_6 -2\beta_4+1\right) + 12 (q)_2 \beta_4,
\\
T\left(1^2 1^2\right) / \mu_2^q
&=&
(q)_4 \left(\beta_4 -1\right)^2 -4 (q)_3 \left(\beta_4 -1 + 2 \mu_3^2\right)+12(q)_2,
\\
T \left(1^21^3\right) / \mu_2^{q}
&=&
12(q)_3 \left(2 \beta_3^2 +3 \beta_4-3\right)
\\
&&
-2(q)_4 \left\{\beta_6-3\beta_4+2+6\beta_3 \left(\beta_5-2\beta_3\right)  + 3 \left(\beta_4-1\right)^2\right\}
\\
&&
+ (q)_5 \left(\beta_4 -1\right) \left(\beta_6 -3 \beta_4 +2\right),
\\
T\left(1^2 1^2 1^2\right) / \mu_2^q
&=&
-120(q)_3+36(q)_4\left(\beta_4-1+4\beta_3^2\right)
\\
&&
-6(q)_5 \left(\beta_4-1 +\beta_3^2\right) \left(\beta_4 -1\right) + (q)_6 \left(\beta_4-1\right)^3.
\end{eqnarray*}
So, $t_i=T_i (F) / T (F)$ and $s_i=S_i (F) / T (F)$ are given by
\begin{eqnarray*}
&&
t_1 = s_1=-(q)_2 \left(\beta_4-1\right)/2 +q,
\\
&&
t_2 = (q)_4 \left(\beta_4 - 1\right)^2/8 + (q)_3 \left( \beta_6/3 - 3 \beta_4 /2 + 7/6\right) + (q)_2 \left(-5 \beta_4 / 2 + 4\right) + q,
\\
&&
s_2 = (q)_4 \left(\beta_4-1\right)^2/8 d+(q)_3 \left(\beta_6/3 - \beta_3^2 -3 \beta_4/2 + 7/6\right)  + (q)_2 \left(-2 \beta_4 +7/2\right),
\\
&&
t_3 = \sum_{i=1}^6 (q)_i t_{3i},
\
s_3 = \sum_{i=2}^6 (q)_i s_{3i}
\end{eqnarray*}
for
\begin{eqnarray*}
&&
t_{31} = 1,
\\
&&
t_{32} = \left(31-19 \beta_4\right) /2,
\\
&&
t_{33} = 4\beta_6- 18 \beta_4- 10 \beta_3^2 + 33 /2,
\\
&&
t_{34} =  \left\{-3 \beta_8+ 16 \beta_6 +24 \beta_5\beta_3-84 \beta_3^2 + 21 \beta_4^2 - 81 \beta_4 + 47 \right\}/12,
\\
&&
t_{35} = s_{35}= \left(\beta_4-1\right) \left(-4 \beta_6 +15\beta_4-11+3\beta_3^2\right) /24,
\\
&&
t_{36} = s_{36}= -\left(\beta_4 -1\right)^3 /48,
\\
&&
s_{32} = -3 \beta_4 +9 /2,
\\
&&
s_{33} = 3\beta_6 -27 \beta_4 /2 +13-7 \beta_3^2,
\\
&&
s_{34} = \left\{-6 \beta_8 +32 \beta_6 -138 \beta_4 + 33 \beta_4^2 + 85
\right\} / 24 -6 \beta_4 - 7 \beta_3^2+2\beta_3\beta_5.
\end{eqnarray*}
\end{example}

\begin{example}
Consider Example 5.13 with $T(F) = \mu_2$, so $ET(\widehat{F})= (1-n^{-1})T (F)$.
As a check $q=1$ above gives
$T(1^2) = -2 \mu_2$, $T(1^3) = T(1^4) = T(1^21^2) = T(1^21^3)=T (1^21^21^2)=0$, so $t_1 = t_2= t_3=1$, $s_1=1$, $s_2=s_3=0$.
\end{example}

\noindent
[Table 5.3 about here.]

\begin{example}
Consider Example 5.13 with $T(F) = \mu_2^{1/2} = \sigma(F)$ say.
Putting $q=1/2$ gives $t_1 =s_1= (\beta_4 +3) /8$,
so an estimate of $\sigma (F)$ of bias $O(n^{-2})$ is
\begin{eqnarray*}
\sigma \left(\widehat{F}\right) \left\{ 1+n^{-1} \left(\beta_4 \left(\widehat{F}\right)+3\right) /8 \right\},
\end{eqnarray*}
where $\beta_4 (F)=\beta_4=\mu_4\mu_2^{-2}$.
To reduce the bias further use
\begin{eqnarray*}
s_2
&=&
\left(16 \beta_6 + 22 \beta_4  + 164 -15 \beta_4^2\right) /128,
\\
s_3
&=&
\bigg( 240 \beta_8 +432 \beta_6-2503 \beta_4 + 2817-165\beta_4^2
\\
&&
+4764 \beta_3^2 + 315\beta_4^3-560\beta_4\beta_6+420 \beta_4\beta_3^2-1920\beta_3\beta_5\bigg) /1024.
\end{eqnarray*}
Table 5.3 gives the relative bias of $S_{np}(\widehat{F})$ estimated from
simulations for $p \leq 2$ and $F$ normal and exponential.
For $n = 100$  and $p=2$, the estimates are poor: see Example C.4.

The usual estimator of $\sigma(F)$ is the sample standard deviation,
s.d. $= \{ n\mu_2(\widehat{F})/(n-1) \}^{1/2}$,
with mean $\sigma \{ 1 - t^*_1 n^{-1} + O(n^{-2}) \}$, where $t^*_1 = t_1 - 1/2$.
So, bias $\{$s.d.$\}$/bias $\{ \sigma( \widehat{F})\} = \lambda_1 + O(n^{-1})$,
where $\lambda_1 = (\beta_4 - 1)/(\beta_4+3)$.

For the normal, exponential and gamma ($\gamma$),
$\beta_4 = 3$, 9 and 3 $+$ 6 $\gamma^{-1}$, so
$\lambda_1 = 1/3$, $2/3$ and $(5\gamma + 12)/(6\gamma + 12)$
and the s.d. improves on $\sigma(\widehat{F})$, although both are first
order estimates, that is, both have bias $O(n^{-1})$.

To see how $S_{n2}(\widehat{F})$ improves on the s.d., note that
bias $\{ S_{n2}(\widehat{F}) \}/$ bias $\{$ s.d. $\} =  \lambda_2 n^{-1} + O(n^{-2})$, where $\lambda_2 = s_2/t^*_1$.
For the normal, exponential and gamma $(\gamma)$,
\begin{eqnarray*}
\beta_6 = 15, 265 \mbox{ and } 120 \gamma^{-2} + 130\gamma^{-1} + 15,
\end{eqnarray*}
so
\begin{eqnarray*}
&&
s_2 = 65/64, 767/32 \mbox{ and } N(\gamma)/64,
\\
&&
\lambda_2 = 65/16 \approx 4.06, 767/32 \approx 24.1 \mbox{ and }
N(\lambda) \left(2.5 + 6\lambda^{-1}\right)^{-1}/64,
\end{eqnarray*}
where $N(\gamma) = 690 \gamma^{-2} + 788 \gamma^{-1} + 65$.
\end{example}

\begin{example}
Suppose $k=s_1=1$, $T(F)= \mu /\sigma= \mu \mu_2^{-1/2}=g (\mu, \mu_2)=\beta$  say.
Again set $\beta_r = \mu_r \mu_2^{-r/2}$.
Then the partial derivatives of $g$ are $g_1 = \mu_2^{-1/2}$, $g_{11} = 0$, $g_2 = -\mu \mu_2^{-3/2}/2$,
$g_{12} = -\mu_2^{-3/2} /2$, $g_{22} = 3 \mu \mu_2^{-5/2} /4$, $g_{122} = 3 \mu_2^{-5/2}/4$,
$g_{222} = -15 \mu \mu_2^{-7/2} /8$, and so on.
Set $U_1(F) = \mu$, $U_2(F) = \mu_2$.
Then defining $U_{ij} \cdots (1^I, 1^J, \ldots)$ as in (\ref{A12})-(\ref{A13}),
\begin{eqnarray*}
&&
U_{11} (1,1) = \int U_{1x}^2=\int \mu_{x}^2=\mu_2,
\\
&&
U_{12} (1,1) = \int U_{1x} U_{2x}=\int \mu_x \mu_{2x}=\mu_3,
\\
&&
U_{22} (1,1) = \int U_{2x}^2 = \mu_4-\mu_2^2.
\end{eqnarray*}
So, by (\ref{A20}),
\begin{eqnarray*}
T \left(1^2\right) = \beta_3+ \beta \left(3 \beta_4+1\right) /4.
\end{eqnarray*}
Also
\begin{eqnarray*}
&&
U_{122} (1,1,1) = \int \mu_{x} \mu_{2x}^2 = \mu_5-2 \mu_2\mu_3,
\\
&&
U_{222} (1,1,1) = \int \mu_{2x}^3=\mu_6-3 \mu_4 \mu_2 +2 \mu_2^3,
\\
&&
U_{12} \left(1,1^2\right) = \int \mu_{x} \mu_{2xx}=-2 \mu_3,
\\
&&
U_{21} \left(1,1^2\right) = \int \mu_{2x} \mu_{xx}=0,
\\
&&
U_{22} \left(1,1^2\right) = \int \mu_{2x} \mu_{2xx}=-2 \left(\mu_4-\mu_2^2\right),
\\
&&
U_1\left(1^3\right) = \int \mu_{xxx},
\\
&&
U_2\left(1^3\right) = \int \mu_{2xxx}=0.
\end{eqnarray*}
So, by (\ref{A21})
\begin{eqnarray*}
T\left(1^3\right) /3=\left(3 \beta_5-2\beta_3\right)/4 +\beta \left(-5 \beta_6+11 \beta_4 -6\right) /8.
\end{eqnarray*}
Similarly, at $(1,1,1^2)$, $U_{221} =0$,
\begin{eqnarray*}
&&
U_{122} = -2\left(\mu_5-\mu_3 \mu_2\right),
\
U_{222}=-2 \left(\mu_6-2\mu_4\mu_2+\mu_2^3\right),
\\
&&
U_{ij}\left(1, 1^3\right) = U_i (1^4)=0,
\
U_{12} \left(1^2, 1^2\right) = 0,
\
U_2 \left(1^2,1^2\right)=4 \mu_4,
\end{eqnarray*}
so by (\ref{A22}),
\begin{eqnarray*}
T\left(1^4\right)=3 \left(-5 \beta_7+3 \beta_5-3 \beta_3\right) /2 +3 \beta \left(35 \beta_8-132 \beta_6 + 242 \beta_4 - 97\right) /16.
\end{eqnarray*}
Also at $a=b=1$,
\begin{eqnarray*}
&&
U_{12}(ab,ab) = U_1\left(a^2b^2\right) = U_2\left(a^2b^2\right) = 0,
\
U_{22}(ab,ab) =\int \mu_{2xy}^2=4\mu_2^2,
\\
&&
U_{122} (ab,a,b) = U_{22} \left(a,ab^2\right) = U_{12} \left(a,ab^2\right) = U_{21} \left(a, ab^2\right) = 0,
\\
&&
U_{222} (ab,a,b) = \int \int \mu_{2xy} \mu_{2x} \mu_{2y} = -2 \mu_3^2.
\end{eqnarray*}
So, by (\ref{A23})
\begin{eqnarray*}
T\left(1^21^2\right)
&=&
4 g_{1222} \mu_3 \left(\mu_4-\mu_2^2\right) +g_{2222}\left(\mu_4-\mu_2^2\right)^2 -
4g_{122} \mu_3\mu_2
\\
&&
-4 g_{222} \left\{ \left(\mu_4-\mu_2^2\right) \mu_2+2 \mu_3^2\right\} +12g_{22}\mu_2^2
\\
&=&
-3\left(5 \beta_4 - 43\right) \beta_3 /8 + 3 \beta \left(35 \beta_4^2+ 90 \beta_4 + 320 \beta_3^2 - 77\right) /16.
\end{eqnarray*}
So,
\begin{eqnarray*}
&&
S_1(F) = T_1 (F)= -\beta_3 /2 - \beta \left(3\beta_4 +1\right) /8,
\\
&&
S_2(F) = \left(48 \beta_5-15 \beta_4\beta_3 -23 \beta_3\right) /64 +
\beta \left(-80 \beta_6 +446 \beta_4 - 327 +105 \beta_4^2 + 960 \beta_3^2\right)/128.
\end{eqnarray*}
Note that $T(1^21^3)$, $T(1^21^21^2)$ and $S_3 (F)$ may be calculated similarly using (\ref{A7}).
\end{example}

\begin{note}
In the one sample example above ${\bm \mu}$ is the mean of ${\bf X} \sim F$.
In many cases ${\bf X}_i = {\bf h}({\bf Y}_i)$, where ${\bf h} : R^t\rightarrow R^s$
is a given transformation and ${\bf Y}_1,\ldots, {\bf Y}_n \sim G$ on $R^t$ is the original sample.
So, ${\bm \mu} (F)=\int {\bf x} dF({\bf x}) = \int {\bf h}({\bf y}) dG({\bf y})$.
Equivalently, we may replace ${\bm \mu} (F) = \int {\bf x} dF ({\bf x})$ by
${\bm \mu} (F)=\int {\bf h} ({\bf x}) dF({\bf x})$, so that ${\bm \mu}_{\bf x} = {\bf h} ({\bf x}) - {\bm \mu}$.
Similarly, if $s=1$ replace $\mu_r (F)=\int (x-\mu)^r dF(x)$ by $\int (h(x)-\mu)^r dF(x)$ so that
(\ref{5.699}) holds with $h_i=h_{x_i} = h(x_i) - \mu$.
A similar remark holds for several samples.
\end{note}

The next four examples apply this idea to return times and exceedances.

\begin{example}
Take $ k=1$, $h({\bf x})=I({\bf x}\leq {\bf a})$ for some ${\bf a}$ in $R^s$, and $T(F)=\mu^{-1}$.
Since $\mu=F({\bf a})$, $T (F)$ is the return period of the event $\{ {\bf X} \leq {\bf a} \}$, where ${\bf X} \sim F$.
But the case $T(F)=\mu^{-1}$ was dealt with in Example 5.3 in terms of $\mu_r$.
In this instance $\mu_r = \mu_r (Bi(1,p))$, where $p=F({\bf a})$, so
$\mu_2=pq$, where $q=1-p$, $\mu_3=pq(1-2p)$ and $\mu_4=pq(1-3p q)$.
So, by Examples 5.6, 5.7 and Note 4.3 an estimate of the return period $p^{-1}$ of
bias $O(n^{-4})$ is $\widetilde{S_{n4}}[\widehat{p}]=S_{n4}[\widehat{p}]$ if $\widehat{p}>l$
or $l^{-1}$ if $\widehat{p}\leq l$, where $0<l<p$,
\begin{eqnarray*}
S_{n4}[p]  =  p^{-1}+\sum^3_{i=1} S_i[p]/(n-1)_i,
\end{eqnarray*}
and $S_i[p] = S_i(F)$ is given by $S_1[p] = p^{-1}-p^{-2}$,
$S_2[p] = -p^{-1}+p^{-3}$, $S_3[p] = 2p^{-1}+p^{-2}-2p^{-3}-p^{-4}$.
\end{example}

The same formula with $p=1-F({\bf a})$ and $\widehat{p}=1-\widehat{F}({\bf a})$ gives an estimate
of bias $O(n^{-4})$ for the return time of the event $\{ {\bf X} > {\bf a} \}$.
Similarly, for the event $\{ {\bf x} \in A \}$ with $p=F(A)$ and $\widehat{p}=\widehat{F} (A)$.
Similarly, we can apply Example 5.4 to obtain estimates of bias $O(n^{-p})$ for any smooth function
$g(p_1,\ldots, p_k)$ given independent $n_i \widehat{p_i} \sim Bi(n_i, p_i)$, $1\leq i \leq k$.
This problem can also be solved by the parametric method of Withers (1987).

\begin{example}
Suppose $k=1$, ${\bf X} \sim F $ on $R^t$ and $T(F)=E r ({\bf X}) \mid ({\bf X} \in A)$,
where $A \subset R^t$ is a measurable set, $F(A)>0$ and $r : R^t \rightarrow R$ is a given function.
Then $T(F) = \mu_1/\mu_2=\mu_1(F)/\mu_2(F)$, where $\mu_i(F)=\int h_i ({\bf x}) dF({\bf x})$, $h_1({\bf x})=r({\bf x}) I ({\bf x} \in A)$
and $h_2({\bf x}) = I({\bf x} \in A)$.
So, $\{T_i, S_i, 1 \leq i \leq 3\}$ are given in Example 5.2 in terms of the moments of (\ref{5.199})
in which $x_{j_i}$ now needs to be replaced by $h_{j_i}({\bf x})$.
Set
\begin{eqnarray*}
p=F(A),
\
q=1-p,
\
I_i=\int_A \left( r({\bf x}) - \mu_1 \right)^i dF({\bf x}).
\end{eqnarray*}
So, $\mu[2^j] = \mu_i (Bi(1,p))$ is given for $2 \leq j \leq 4$ in Example 5.17 and
\begin{eqnarray*}
\mu \left[1^i 2^j\right] = I_iq^j + \left(-\mu_1\right)^i (-p)^j q.
\end{eqnarray*}
Using $I_1=0$ simplification yields
\begin{eqnarray*}
S_{n4} (F) = \mu_1 p^{-1} \left\{ 1-q^2p^{-1}/(n-1)+q^3p^{-2}/(n-1)_2 +q^3p^{-3}(2p-1)/(n-1)_3\right\}.
\end{eqnarray*}
Unlike Example 5.17, one does not need to know a lower bound for $p$, since
$\mu_{1} =0$ if $p=0$; so, if $\widehat{p}=0$ one interprets $S_{n4}(\widehat{F})$ as an
arbitrary constant.
This shows, surprisingly that the bias reduction problem
for $T(F)=\mu_1 /p$ can be treated as a parametric problem, the parameters
being $(\mu_1, p)$.
The more general problem of $T(F)=g(\mu_1,p)$ does
{\it not} reduce to a finite parameter problem as it involves $\{\int_A r^i dF, i \geq 1\}$.
\end{example}

\begin{example}
The conditional distribution of exceedances is
\begin{eqnarray}
F_u(x) = P\left(X-u < x \mid X-u > 0\right)=  \left\{ F(x+u)- F(u)\right\}/\left\{1-F(u)\right\}
\label{5.2899}
\end{eqnarray}
for $x \geq 0$.
This is $\mu_1/\mu_2$ with $A=\{y:y>u\}=(u,\infty)$, $B-\{y:x+u>y>u\} = (u,x+u)$ and $r(y) = I(y \in B)$.
So, Example 5.18 applies with $\mu_1= F(x+u)-F(u)$, $\mu_2=1-F(u)$.
\end{example}

\begin{example}
The mean conditional exceedance is
\begin{eqnarray*}
\mu \left(F_u \right) = \int x dF_u(x)=\mu_1/\mu_2
\end{eqnarray*}
for
\begin{eqnarray*}
\mu_1 = \int(x-u)_{+} dF(x),
\
\mu_2 =1-F(u),
\end{eqnarray*}
where
\begin{eqnarray*}
x_{+} = \left\{ \begin{array}{ll}
x,  & \mbox{if $x>0$,}\\
0,  & \mbox{if $x\leq 0$.}
\end{array}
\right.
\end{eqnarray*}
So, $r(y)=(y-u)_{+}$ and Example 5.18 applies.
\end{example}

The {\it central} moments of $F_u$ of (\ref{5.2899}) are {\it not}
covered by Example 5.18 and are probably best dealt with by writing them as
functions of the noncentral moments and applying Example 5.1 with
$\mu=\{\int (x-u)_{+}^i dF(x), i\geq 0\}$.
A more direct approach is given by the following example.

\begin{example}
Suppose $T(F) =S(F_{u})$ for $F_u$ of (\ref{5.2899}).
Set $C^y(F) = F_{u}(y)$.
Then
\begin{eqnarray*}
C^y \left((1-\epsilon)F+\epsilon \delta_x\right)  =  F_u(y)
+\epsilon C^y_F(x) + O\left(\epsilon^2\right),
\end{eqnarray*}
and
\begin{eqnarray*}
T\left((1-\epsilon)F+\epsilon \delta_x \right)= S\left(F_u(\cdot)+\epsilon C^{\cdot}_F(x) +O \left(\epsilon^2\right)\right) =
S(F)+\epsilon \int S_{F_u}(y)  C^y_F (x) dy +O \left(\epsilon^2\right),
\end{eqnarray*}
where $C^y_F(x) = \mu_2^{-1}I(u<x<u+y)-\mu_1\mu_2^{-2} I(u<x)$.
So,
\begin{eqnarray}
T_F(x) = \int S_{F_u}(y) C^y_F(x) dy =\mu_2^{-1}S_{F_{u}}(x-u).
\label{5.2999}
\end{eqnarray}
Higher derivatives can be calculated from (\ref{5.2999}).
\end{example}

Now let us apply the previous note with $s=1$, $t=r$, $h({\bf y})={\bf a}'{\bf y}$, where ${\bf a}$ lies in $R^r$.
Set ${\bm \mu} = E {\bf Y}$.
Then the joint central moment $\mu_{1 \cdots r}=E({\bf Y}-{\bm \mu})_1 \cdots ({\bf Y}-{\bm \mu})_r$
is the coefficient of $a_1 \cdots a_r/r!$ in $\mu_r({\bf a}' {\bf Y})$,
so the same relation is true of their derivatives.
The same is also true of the cumulants.
This device allows us to derive results
for multivariate moments and cumulants from their univariate analogs.

For example, from Example 5.6, for a univariate random variable,
$\mu_{2}(x)=(x-\mu)^2-\mu_2$ and $\mu_{2}(x_1,x_2)=-2(x_1-\mu)(x_2-\mu)$.
So, for a bivariate random variable,
$\mu_{12}({\bf x})=({\bf x}-{\bm \mu})_1({\bf x}-{\bm \mu})_2-\mu_{12}$ and
$\mu_{12}({\bf x}_1, {\bf x}_2)=-2({\bf x}_1 - {\bm \mu})_1 ({\bf x}_2 - {\bm \mu})_2$.

We illustrate this device further with the problems of estimating
multivariate moments and the correlation of a bivariate distribution
and its square.

\begin{example}
Suppose $k=1$, $s=2$ and $T(F)=\mu_{12}$.
From Example 5.6 and the previous remark, an UE of $\mu_{12}$ is $\mu_{12}/(1-n^{-1})$ at $F=\widehat{F}$.
\end{example}

Similarly, we have
\begin{example}
Suppose $k=1$, $s=3$ and $T(F)=\mu_{123}$.
An UE of
$\mu_{123}$ is  $\mu_{123}/\{(1-n^{-1})(1-2n^{-2})\}$ at $F=\widehat{F}$.
\end{example}

\begin{example}
Suppose $k=1$, $s=2$, and $T(F)=\mu_{12}\{\mu_{11}\mu_{22}\}^{-1/2}$, the correlation of a bivariate sample.
So, (\ref{A1}) of Appendix A holds with ${\bf S} (F)=(\mu_{12}, \mu_{11}, \mu_{22})$ and $g({\bf S})=S_1(S_2S_3)^{-1/2}$.
We shall apply (\ref{A8}).
Set $\nu_{ij \cdots}=\mu_{ij \cdots}(\mu_{ii}\mu_{jj} \cdots)^{-1/2}$.
So, $T(F)=\nu_{12}$.
Now $S_1(1^2)=\int S_{1 {\bf xx}} =-2\mu_{12}$, $S_2(1^2)=\int S_{2 {\bf xx}}=-2\mu_{11}$ and $S_3(1^2)=\int S_{3 {\bf xx}}=-2\mu_{22}$.
Also $g_1=(\mu_{11}\mu_{22})^{-1/2}$, $g_2=-\nu_{12}/\mu_{11}$, $g_3=-\nu_{12}/\mu_{22}$.
So, $g_iS_i(1^2)=T(F)(-2+1+1)=0$.
Similarly, $S_{1{\bf x}}=({\bf x}-{\bm \mu})_1({\bf x}-{\bm \mu})_2-\mu_{12}$, so
$ S_{11}(1,1)=\int S_{1{\bf x}}^2=\mu_{1122}-\mu_{12}^2$, and similarly
$ S_{12}(1,1)=\mu_{1112}-\mu_{11}\mu_{12}$,
$S_{13}(1,1)=\mu_{1222}-\mu_{12}\mu_{22}$,
$S_{22}(1,1)=\mu_{1111}-\mu_{11}^2$,
$S_{33}(1,1)=\mu_{2222}-\mu_{22}^2$, and $S_{23}(1,1)=\mu_{1122}-\mu_{11}\mu_{22}$.
So, an estimate of bias $O(n^{-2})$ is
$T(F)-T(1^2)/(2n)$ or $T(F)-T(1^2)/(2n-2)$ at $F=\widehat{F}$, where by (\ref{A8}),
$T(1^2)=\nu_{12}(3\nu_{1111}+3\nu_{2222}+2\nu_{1122})/4-\nu_{1112}-\nu_{1222}$.
\end{example}

\begin{example}
Suppose $k=1$, $s=2$ and $T(F)=\mu_{12}^2\{\mu_{11}\mu_{22}\}^{-1}= \nu_{12}^2$, the square of the correlation of a bivariate sample.
Again (\ref{A1}) holds with ${\bf S} (F)=(\mu_{12}, \mu_{11}, \mu_{22})$ but now $g({\bf S}) = S_1^2(S_2S_3)^{-1}$, so $g_1=2T(F)S_1^{-1}$,
$g_2=-T(F)S_2^{-1}$,
$g_3=-T(F)S_3^{-1}$,
$g_{ii}=2T(F)S_i^{-2}$,
$g_{12}=-2T(F)(S_1S_2)^{-1}$,
$g_{13}=-2T(F)(S_1S_3)^{-1}$, and $g_{23}=T(F)(S_2S_3)^{-1}$.
Again $g_iS_i(1^2)=T(F)(-4+2+2)=0$.
So, an estimate of bias $O(n^{-2})$ is
$T(F)-T(1^2)/(2n)$ or $T(F)-T(1^2)/(2n-2)$ at $F=\widehat{F}$, where by (\ref{A8}),
$T(1^2)=2\nu_{12}^2(\nu_{1111}+\nu_{2222}+2\nu_{1122}-2\nu_{1112}-2\nu_{1222})$.
\end{example}

\section{Estimating covariances of estimates}
\setcounter{equation}{0}

In this section, we give an estimate of bias $O(n^{-3})$ for ${\bf V}_n (F)$,
the covariance of ${\bf T} \widehat{(F)}$, where now ${\bf T} (F)$ is a $q \times 1$
{\it vector} with {\it components}
$ \{ T^\alpha (F), 1 \leq \alpha \leq q \}$.
After Example 6.1, we estimate the covariance of more general estimates of ${\bf T}(F)$.

From the formulas for $ \{K_i^{ab} \}$ on Withers (1982, pages 66 and 67),
\begin{eqnarray}
V_n^{\alpha \beta} (F) =  covar  \left( T^\alpha \widehat{(F)},
T^\beta  \widehat{(F)} \right) = \sum_{i=1}^{\infty} n^{-i } K_i^{\alpha \beta} (F),
\label{6.199}
\end{eqnarray}
where
\begin{eqnarray}
K_1^{\alpha \beta} (F)
&=&
t_i^\alpha t_j^\beta  k^{ij} = \sum \lambda_a \int \int T^\alpha_F \left( \begin{array}{c} a \\ x
  \end{array} \right ) T_F^\beta  \left ( \begin{array}ca \\ y
  \end{array} \right ) d \kappa_a (xy)
\nonumber
\\
&=&
\sum \lambda_a T^{\alpha \beta} (a,a),
\label{6.299}
\\
K_2^{\alpha \beta} (F)
&=&
\sum^2 t_{ij}^\alpha t_k^\beta  k^{ijk} /2 + \left( \sum^2 t^\alpha_{ijk} t^\beta _{l} +  t^\alpha_{ik} t^\beta _{jl}\right) k^{ij} k^{kl} /2
\nonumber
\\
&=&
\sum \lambda_a \sum^2 \int^3 T^\alpha_F \left( \begin{array}{cc} a&a\\ x&y \end{array} \right )T^\beta_F  \left(
\begin{array}{cc} a\\z \end{array} \right ) d \kappa_a (xyz) /2
\nonumber
\\
&&
+ \sum \lambda_a \lambda_b \int^4 \left \{ \sum^2 T^\alpha_F \left(
  \begin{array}{ccc} a&a&b \\ w&x&y \end{array} \right ) T^\beta_F\left(
  \begin{array}cb\\z \end{array} \right ) \right.
\nonumber
\\
&&
\left. + T^\alpha_F \left(
  \begin{array}{cc} a&b \\ w&x \end{array} \right ) T^\beta_F
  \left( \begin{array}{cc} a&b\\y&z \end{array} \right ) \right \}
  d \kappa_a (wx) d \kappa_b (yz) /2
\nonumber
\\
&=&
\sum \lambda_a \sum^2 T^{\alpha\beta} \left(a^2, a\right)/2
\nonumber
\\
&&
+ \sum\lambda_a \lambda_a \left\{ \sum^2 T^{\alpha\beta} \left(a^2b,b\right) +
T^{\alpha\beta} (ab, ab) \right\} /2,
\label{6.399}
\\
\sum^2 f_{\alpha\beta}
&=&
f_{\alpha\beta} + f_{\beta\alpha},
\nonumber
\\
T^{\alpha\beta}(a,a)
&=&
\int T^\alpha_F \left( \begin{array}c a \\ x
\end{array} \right ) T^\beta _F \left( \begin{array}c a \\ x
\end{array} \right )dF_a (x),
\nonumber
\\
T^{\alpha\beta} \left(a^2, a\right)
&=&
\int T^\alpha_F \left(
\begin{array}{cc} a&a\\ x&x \end{array} \right )T^\beta_F  \left(
\begin{array}c a \\ x \end{array} \right ) dF_a (x),
\label{6.599}
\\
T^{\alpha\beta} \left(a^2b, b\right)
&=&
\int \int T^\alpha_F \left(
\begin{array}{ccc} a&a&b \\ x&x&y \end{array} \right ) T^\beta_F\left(
\begin{array}cb \\ y \end{array} \right) d F_a (x) dF_b (y),
\label{6.699}
\end{eqnarray}
and
\begin{eqnarray}
T^{\alpha\beta} (ab,ab) = \int\int T^\alpha_F \left(
\begin{array}{cc}a&b\\ y&x \end{array} \right) T^\beta_F  \left(
\begin{array}{cc}a&b\\ x&y \end{array} \right) d F_a (x) dF_b (y).
\label{6.799}
\end{eqnarray}
Also, setting $V^{\alpha\beta} (F) $ = $K^{\alpha \beta}_1 (F)$ and differentiating, we have
\begin{eqnarray*}
V^{\alpha\beta}_F \left( \begin{array}ca \\
x \end{array} \right ) / \lambda_a
 = T^\alpha_F \left( \begin{array}ca \\
   x \end{array}
   \right ) T_F^\beta  \left( \begin{array}ca\\
x \end{array} \right) -
T^{\alpha\beta} (a,a) + \sum^2 \int T_F^\alpha \left(
\begin{array}{cc}a&a\\
  y & x \end{array} \right ) T_F^\beta  \left(
\begin{array}c a\\y \end{array} \right )d F_a (y),
\end{eqnarray*}
and
\begin{eqnarray*}
V^{\alpha\beta}_F \left( \begin{array}{cc}
   a&a \\ x&x \end{array} \right ) /  \lambda_a
&=&
\sum^2 \bigg[
  \left\{T_F^\alpha  \left(
  \begin{array}{cc} a&a\\x&x \end{array} \right) - T_F^\alpha
  \left(\begin{array}c a \\ x \end{array} \right) \right\} T_F^\beta
  \left( \begin{array}ca \\ x \end{array} \right)
\\
&&
+ T_F^\alpha \left(
  \begin{array}{cc} a&a\\x&x \end{array}\right) T_F^\beta  \left(
  \begin{array}c a\\x \end{array} \right)
\\
&&
- \int T_F^\alpha \left( \begin{array}{cc} a&a\\x&y \end{array} \right )T_F^\beta
\left( \begin{array}c a\\y \end{array} \right) d  F_a (y)
\\
&&
+ \int T_F^\alpha \left( \begin{array}{cc} a&a\\x&y
  \end{array} \right) T_F^\beta  \left( \begin{array}{cc} a&a\\x&y
  \end{array} \right) d F_a (y)
\\
&&
+ \int \{ T_F^\alpha \left( \begin{array}{ccc} a&a&a\\x&x&y
  \end{array} \right) - T_F^\alpha \left( \begin{array}{cc}
 a&a\\x&y\end{array} \right) \} T_F^\beta  \left( \begin{array}c
  a\\y \end{array} \right) d F_a (y) \bigg],
\end{eqnarray*}
so that
\begin{eqnarray}
C_1 \left(V^{\alpha\beta}, F\right)
&=&
\sum \lambda_a V^{\alpha\beta}(a^2)
\nonumber
\\
&=&
\sum \lambda_a^2 \left[ \sum^2 \left\{ T^{\alpha \beta} \left(a^2,a\right) +
T^{\alpha\beta} \left( a^{2} b, b\right) /2 \right\} +
2T^{\alpha\beta} (a b, a b) - T^{\alpha \beta} (a,a) \right]_{b=a}.
\nonumber
\end{eqnarray}
So, $n^{-1} K_1^{\alpha\beta} \widehat{(F)}$ given by (\ref{6.299}) estimates
$V_n^{\alpha\beta} (F)$ with bias $O(n^{-2})$ and
$n^{-1}K_1^{\alpha\beta} \widehat{(F)}+n^{-2} L^{\alpha \beta} \widehat{(F)}$
estimates $V^{\alpha\beta}_n (F)$ with bias $O(n^{-3})$, where
\begin{eqnarray}
L^{\alpha\beta} (F)
&=&
K_2^{\alpha\beta} (F) - C_1 \left( V^{\alpha\beta}, F \right)
\nonumber
\\
&=&
\sum \left(\lambda_a-\lambda_a^2\right) \sum^2 T^{\alpha\beta} \left( a^2, a \right)/2
\nonumber
\\
&&
+\sum\lambda_a \lambda_b \left\{ \sum^2 T^{\alpha\beta} \left(a^2 b,b\right) + T^{\alpha\beta} (ab,ab) \right\}/2
\nonumber
\\
&&
-\sum\lambda_a^2 \left\{\sum^2 T^{\alpha\beta} \left(a^{2}b, b\right)/2+
2 T^{\alpha\beta}(a b,a b) - T^{\alpha\beta} (a,a) \right\}_{b=a}.
\nonumber
\end{eqnarray}
If $k = 1$ this reduces to
\begin{eqnarray}
L^{\alpha\beta} (F) = T^{\alpha\beta} (a,a) - 3T^{\alpha\beta} (a b, a b) /2
\label{6.1099}
\end{eqnarray}
at  $a=b=1$, so that
\begin{eqnarray}
(n-1)^{-1} T^{\alpha\beta} (a,a) - 3n^{-2} T^{\alpha\beta} (ab, ab)/2
\label{6.1199}
\end{eqnarray}
at  $\{ F=\widehat{F}, \ a=b=1 \}$ estimates $V^{\alpha\beta}_n (F)$ with bias $O (n^{-3})$, where at $a=b=1$,
\begin{eqnarray*}
T^{\alpha\beta} (a,a) = \int T_F^\alpha (x) T_F^\beta  (x) dF (x),
\end{eqnarray*}
and
\begin{eqnarray*}
T^{\alpha\beta} (ab,ab) = \int \int T_F^\alpha (xy) T^\beta_F  (xy) dF(x) dF (y).
\end{eqnarray*}
One may prefer to use $n^{-1} - n^{-2}$ instead of $(n-1)^{-1}$ in (\ref{6.1199}).
Remarkably, unlike the case $ k > 1$, the estimate (\ref{6.1199}) does not depend on $T^{\alpha \beta} (a^2, a)$ or $T^{\alpha\beta} (a^2 b, b)$ at $a=b=1$.

We now show how to estimate
\begin{eqnarray}
{\bf W}_n (F) = covar \ {\bf T}_{(n)} \widehat{(F)},
\label{6.1299}
\end{eqnarray}
where
\begin{eqnarray*}
{\bf T}_{(n)} = \sum_{i=0}^{\infty} n^{-i} {\bf T}_i
\end{eqnarray*}
is $q \times 1$ and ${\bf T}_0 = {\bf T}$.
Clearly, ${\bf T}_{(n)} \widehat{(F)}$ estimates ${\bf T}(F)$.
Now
\begin{eqnarray*}
{\bf W}_n (F) = \sum_{i,j \geq 0} n^{-i-j} {\bf W}_n \left( {\bf T}_i, {\bf T}_j \right),
\end{eqnarray*}
where
\begin{eqnarray*}
{\bf  W}_n \left( {\bf T}_i, {\bf T}_j \right) = covar \left( {\bf T}_i \widehat{(F)}, {\bf T}_j \widehat{(F)} \right)
\end{eqnarray*}
has $(\alpha,\beta)$ element
\begin{eqnarray*}
W_n^{\alpha\beta} \left( {\bf T}_i, {\bf T}_j \right) = {\bf W}_n \left( T^\alpha_i, T_j^\beta \right) = V_n^{12} (F)
\end{eqnarray*}
of (\ref{6.199}) with $(T^1, T^2) = (T^\alpha_i, T^\beta _j)$.
So,
\begin{eqnarray*}
W_n^{\alpha\beta} (F) = \sum_{l=1}^{\infty} n^{-l}K_{l}^{\alpha\beta} [F],
\end{eqnarray*}
where
\begin{eqnarray*}
K_{l}^{\alpha\beta} [F] = \sum_{i+j+k=l} K_k \left(T^\alpha_i, T^\beta _j\right),
\end{eqnarray*}
and
\begin{eqnarray*}
K_k \left(T^1, T^2 \right) = K_k^{12} (F) \mbox{ of } (\ref{6.199}).
\end{eqnarray*}
So,
\begin{eqnarray*}
K_1^{\alpha\beta} [F] = K_1 \left(T^\alpha, T^\beta \right) = K_1^{\alpha\beta} (F)
\end{eqnarray*}
of (\ref{6.299}), and
\begin{eqnarray*}
K_2^{\alpha\beta} [F] = K_2^{\alpha\beta} (F) + \triangle^{\alpha\beta},
\end{eqnarray*}
where
\begin{eqnarray*}
\triangle^{\alpha\beta}= \sum^2 K_1 \left(T^\alpha, T^\beta _1\right),
\nonumber
\end{eqnarray*}
and
\begin{eqnarray*}
K_1 \left(T^\alpha, T_1^\beta \right) = K_1^{\alpha\beta}(F)
\end{eqnarray*}
of (\ref{6.299}) at $T^\beta  = T_1^\beta$.

So, $n^{-1} K_1^{\alpha\beta} \widehat{(F)}$ and $n^{-1} K_1^{\alpha\beta} \widehat{(F)} + n^{-2} L^{\alpha\beta} \widehat{[F]}$
estimate $W^{\alpha \beta}_n (F)$ with bias $O(n^{-2})$ and $O(n^{-3})$, respectively, where
\begin{eqnarray}
L^{\alpha\beta} [F] = K_2^{\alpha\beta} [F]- C_1 \left(V^{\alpha\beta}, F\right) = L^{\alpha\beta} (F) + \triangle^{\alpha\beta}.
\label{6.1699}
\end{eqnarray}
Alternatively, for $k = 1$, the sum of (\ref{6.1199}) and $n^{-2} \triangle^{\alpha\beta}$ at $F = \widehat{F}$
estimates $W^{\alpha \beta}_n (F)$ with bias $O(n^{-3})$.
Now for $p \geq 2$, $T_{np}$ of (\ref{1.3}) has the form ${\bf T}_{(n)}$ of (\ref{6.1299}) with $T_1$ given by (\ref{4.1}), so that
\begin{eqnarray*}
T_{1F}^\beta \left( \begin{array} c a \\ x \end{array} \right ) = -\lambda_a \left\{ T_F^\beta \left( \begin{array}c
a^2\\ x^2 \end{array} \right) - T^\beta \left(a^2\right) + \int T_F^\beta\left( \begin{array}{cc}a^2a \\y^2x \end{array} \right)
dF_a (y) \right\} /2,
\end{eqnarray*}
and so
\begin{eqnarray}
K_1 \left( T^\alpha, T^\beta_1 \right)
&=&
-\sum \lambda_a^2 \left\{ T^{\beta\alpha} \left(a^2, a\right) + T^{\beta\alpha} \left(a^3, a\right) \right\} /2,
\nonumber
\\
\triangle^{\alpha\beta}
&=&
-\sum \lambda_a^2 \sum^2 \left\{ T^{\alpha\beta} \left(a^2,a\right) + T^{\alpha\beta} \left(a^{2} b, b\right) \right\}_{b=a} /2,
\label{6.1899}
\\
K_2^{\alpha\beta} [F]
&=&
\sum \left(\lambda_a - \lambda_a^2\right) \sum^2 T^{\alpha\beta} \left(a^{2}b, b\right)/2 -
\sum \lambda_a^2 \sum^2 T^{\alpha\beta} \left(a^2 b, b\right)_{b=a}/2
\nonumber
\\
&&
+ \sum \lambda_a\lambda_b \left\{ \sum^2 T^{\alpha\beta} \left(a^2 b, b\right) + T^{\alpha\beta} (ab,ab) \right\}/2,
\nonumber
\\
L^{\alpha\beta}[F]
&=&
\sum \left( \lambda_a/2 - \lambda_a^2 \right) \sum^2 T^{\alpha\beta} \left(a^2,a\right)
\nonumber
\\
&&
+\sum \lambda_a\lambda_b\left\{ \sum^2 T^{\alpha\beta} \left(a^2 b, b \right) +
T^{\alpha\beta} (ab,ab) \right\} /2
\nonumber
\\
&&
-\sum \lambda_a^2 \left\{ \sum^2 T^{\alpha\beta} \left(a^{2} b, b\right) + 2T^{\alpha\beta} (ab, ab) -
T^{\alpha\beta} (a,a) \right\}_{b=a}.
\nonumber
\end{eqnarray}
For $k = 1$, at $ a=b=1$, this gives
\begin{eqnarray}
&&
\triangle^{\alpha\beta} = -\sum^2\left\{ T^{\alpha\beta} \left(a^2, a\right) + T^{\alpha\beta} \left(a^2b, b\right) \right\} /2,
\nonumber
\\
&&
K_2^{\alpha\beta} [F] = T^{\alpha\beta} (ab, ab) /2,
\label{6.2199}
\\
&&
covar  \ \left( T^\alpha_{np} \widehat{(F)}, T^\beta_{np}(\widehat{F}) \right) = n^{-1} T^{\alpha\beta} (a, a) +  n^{-2}T^{\alpha\beta} (ab, ab) /2 + O\left(n^{-3}\right)
\nonumber
\end{eqnarray}
which, remarkably, does not depend on $T(a^2, a )$ or $T(a^2 b, b)$ to this accuracy - whereas $L^{\alpha\beta} [F]$ does.

\begin{example}
Consider again Example 5.1, that is $k=1$, ${\bf T} (F) = {\bf g} ({\bm \mu})$, where now ${\bf g}$ may be a vector $ \{ g^\alpha \}$.
By (\ref{A16})-(\ref{A19}) at $a=b=1$
\begin{eqnarray}
&&
K_1^{\alpha\beta} (F) = T^{\alpha\beta} (a, a) =g_i^\alpha g_j^\beta \mu [ij],
\nonumber
\\
&&
T^{\alpha\beta} (ab, ab) = g_{ij}^\alpha \ g_{kl}^\beta  \ \mu[ik] \ \mu [jl],
\nonumber
\\
&&
T^{\alpha\beta} \left(a^2, a\right) = g_{ij}^\alpha  \  g_k^\beta \ \mu [ijk],
\nonumber
\\
&&
T^{\alpha\beta} \left(a^2b, b\right) = g_{ijk}^\alpha\ g_l^\beta  \mu [ij] \mu [kl],
\nonumber
\end{eqnarray}
and $K_2^{\alpha\beta} (F)$, $L^{\alpha\beta} (F)$, $K_2^{\alpha\beta} [F]$, $L^{\alpha\beta}[F]$ are given by (\ref{6.399}), (\ref{6.1099}),
(\ref{6.1699}), (\ref{6.1899}), (\ref{6.2199}).
Note that $L^{\alpha\beta}$ depends only on the first and
second moments of $F$, even though $K_2^{\alpha\beta}$ depends on the third moments!
\end{example}

\begin{example}
Consider Example 6.1 with $g({\bm \mu}) = {\bm \alpha}^\prime {\bm \mu} / {\bm \beta}^\prime {\bm \mu} = N/D$, say, -- that is, Example 5.2.
Since $q=1$ we drop suffixes $\alpha$, $\beta$.
Define $\mu [\cdot]$ and $\delta_i$ as in (\ref{5.199}) and (\ref{5.399}).
Then at $a=b=1$
\begin{eqnarray*}
&&
K_1 (F) = T(a,a) =D^{-2} \mu_2 [\delta\delta],
\\
&&
T(ab, ab) = 2\mu_2 [\delta \beta]^2 + 2\mu_2 [\delta\delta]\mu_2 [\beta\beta],
\\
&&
T\left(a^2, a\right) = -2D^{-3} \mu_3 [\delta\delta\beta],
\\
&&
T\left(a^2 b, b\right) = 2D^{-4} \left\{ 2 \mu_2 [\delta \beta]^2 + \mu_2 [\delta \delta] \mu_2 [\beta \beta] \right\},
\end{eqnarray*}
where $\mu_2 [\delta\beta] = \delta_i \beta_j \mu [ij]$ and $\mu_3 [\alpha\beta\gamma] = \alpha_i\beta_j\gamma_k \mu [ijk]$.
In particular, for $g ({\bm \mu}) = \mu_1 / \mu_2$, at $a=b=1$ setting $\gamma_{ij \cdots} = \mu (ij \cdots) \mu_i^{-1} \mu_j^{-1} \cdots$, we have
\begin{eqnarray}
&&
K_1 (F) = T(a,a) = \left(\mu_1 / \mu_2\right)^2 \left(\gamma_{11} - 2 \gamma_{12} + \gamma_{22} \right),
\label{5.26}
\\
&&
T (ab, ab) = 2 \left(\mu_1 / \mu_2\right)^2 \left(\gamma_{11} \gamma_{22} - 4 \gamma_{12} \gamma_{22} + 2\gamma_{22}^2 \right),
\nonumber
\\
&&
T \left(a^2, a\right) = -2\left(\mu_1 / \mu_2\right)^2 \left(\gamma_{112} - 2\gamma_{122}+\gamma_{222}\right),
\\
&&
T \left(a^2b, b\right) = 2 \left(\mu_1 / \mu_2\right)^2 \left(2 \gamma_{12}^2 - 5 \gamma_{12} \gamma_{22} + 3 \gamma_{22}^2 +\gamma_{11} \gamma_{22} \right).
\nonumber
\end{eqnarray}
Note that (\ref{5.26}) is in agreement with equation (10.17) of Kendall and Stuart (1977).
\end{example}

\begin{example}
Consider Example 6.1 with $g (\mu) = N^p$, where $N = {\bm \alpha}^{\prime} {\bm \mu}$, that is, we consider Example 5.3.
In the notation there, with $a=b=1$
\begin{eqnarray*}
&&
K_1 (F) = T(a, a) = p^2 N^{2p} \alpha_{(2)},
\\
&&
T(ab, ab) = p^2 (p-1)^2 N^{2p} \alpha_{(2)}^2,
\\
&&
T \left(a^2, a\right) = p^2 (p-1)N^{2p} \alpha_{(3)},
\\
&&
T \left(a^2 b, b\right) = (p)_3 p N^{2p} \alpha_{(2)}^2.
\end{eqnarray*}
In particular, for $s=1$ and $g (\mu) = \mu^p$, with $a=b=1$
\begin{eqnarray}
&&
T(a, a) = p^2 \mu^{2p-2}\mu_2,
\
T(ab, ab)=p^2 (p-1)^2 \mu^{2p-4} \mu_2^2,
\nonumber
\\
&&
T \left(a^2, a\right) = p^2 (p-1) \mu^{2p-3} \mu_3,
\
T \left(a^2 b, b\right) = (p)_3 p \mu^{2p-4} \mu_2^2.
\nonumber
\end{eqnarray}
For example,  $ var \{ \widehat{\mu}^{-1} \}$ or (if Note 4.3 needs to be
applied), $var \{ \widehat{\mu}^{-1} I \left( | \widehat{\mu} |> l \right) \}$,
where $l > 0$ is a known lower bound for $|\mu|$, can be estimated by
\begin{eqnarray*}
\widehat{T}_{n2} = (n-1)^{-1} \widehat{\mu}^{-4} \widehat{\mu}_2 - 6n^{-2}
\widehat{\mu}^{-6} \widehat{\mu}^2_2
\end{eqnarray*}
or by
\begin{eqnarray*}
\widehat{T}_{n2} I \left( |\widehat{\mu}| > l \right)
\end{eqnarray*}
with bias $O(n^{-3})$, where $( \widehat{\mu}, \widehat{\mu}_2)$ is $( \mu, \mu_2)$ at $F = \widehat{F}$.
Alternatively, replacing $n^{-2}$ in $\widehat{T}_{n2}$ by $(n - 1)^{-2}$ and setting $s^2 = \widehat{\mu}_2 \ n/(n-1)$,
the UE of $\mu_2$, we obtain
\begin{eqnarray*}
T^{\star}_{n2} = n^{-1} \widehat{\mu}^{-4} s^2 - 6n^{-2} \widehat{\mu}^{-6} s^4,
\
T^{\star}_{n2} I \left( |\widehat{\mu}| > l \right)
\end{eqnarray*}
as estimates with bias $O (n^{-3})$.
\end{example}

\section{Estimating the covariance of an estimate of bias}
\setcounter{equation}{0}

The emphasis of this paper has been to reduce bias, not estimate it.
However, a number of papers have given methods for estimating the
variance of an estimate of bias for the case $k = 1$.
See, for example, Efron (1981) and Davison {\it et al.} (1986).
These papers provide bootstrap and jackknife methods of an order of magnitude less
efficient computationally than the Taylor series method (also called
the delta method or the infinitesimal jackknife when $p=2$) used here.

Suppose then ${\bf T}(F)$ is a $q \ \times \ 1$ functional.
Note that ${\bf T} \widehat{(F)}$ has bias $n^{-1} {\bf B} (F) /2 + O (n^{-2})$, where ${\bf B} (F)= |2| = \sum \lambda_a T (a^2)$.
Its estimate $ n^{-1} {\bf B} \widehat{(F)} /2 $ has covariance $n^{-2} {\bf V} (F) /4 + O (n^{-3})$, where
\begin{eqnarray*}
V^{\alpha\beta} (F)
&=&
\sum \lambda_a \int  B_{F}^\alpha \left( \begin{array}c a \\ x \end{array} \right )
B_{F}^\beta   \left( \begin{array}c a \\ x \end{array} \right )   dF_a (x)
\\
&=&
\sum \lambda_a^3 \bigg\{ \int T^\alpha \left(
\begin{array}{cc} aa\\xx \end{array} \right) T^\beta  \left(
\begin{array}{cc} aa \\ xx\end{array} \right) - T^\alpha \left(a^2\right) T^\beta  \left(a^2\right)
\\
&&
+ \sum^2 \int \int T^\alpha \left( \begin{array}{ccc} aaa \\
xxy\end{array} \right) T^\beta  \left( \begin{array}{cc} aa \\
yy\end{array} \right) + \int\int\int T^\alpha \left(
  \begin{array}{ccc} aaa \\
   xxz \end{array} \right) T^\beta  \left(
  \begin{array} {ccc} aaa \\ yyz \end{array} \right)  \bigg\}
\end{eqnarray*}
and $dF_a (x)$, $d F_a (y)$, $dF_a (z)$ are implicit in the integrals.
Finally, $ n^{-2} {\bf V} \widehat{(F)} /4 $ estimates $ covar \left\{ n^{-1} {\bf B} \widehat{(F)} /2  \right \}$ with bias $O(n^{-3})$.

The same is true if we replace ${\bf B} \widehat{(F)}$ by ${\bf B}_{np} \widehat{(F)}$.
If desired, one could apply Section 6 to reduce this bias to $O( n^{-4})$.

\begin{note}
In equation (2.6) of Davison {\it et al.} (1986) and the following
line a factor $ 1/2 $ should be inserted.
So, the
usual bootstrap and the usual jackknife estimates of bias as
well as our estimate $n^{-1} {\bf B} (F) /2 $, all have bias $O (n^{-2})$.
\end{note}

\section*{Appendix A}
\renewcommand{\theequation}{$\mbox{A.\arabic{equation}}$}
\setcounter{equation}{0}

Here, we note and illustrate the following chain rule for the partial
derivatives of
\begin{eqnarray}
T(F) = g( {\bf S}(F)),
\label{A1}
\end{eqnarray}
where ${\bf S} (F)$ is $q \times 1$ and $g: R^{q} \rightarrow R$.

First, suppose $k = 1$, that is, $F$ is a single d.f.
Given $r \geq 1$, let ${\bf s} ({\bf y}) : R^r \rightarrow R^{q}$ be an
arbitrary function.
Set $ \partial_i = \partial/\partial y_i$.
Then
\begin{eqnarray}
T_{F} \left( {\bf x}_1 \cdots {\bf x}_r\right) = \partial_1 \cdots \partial_r g ({\bf s}({\bf y})),
\label{A2}
\end{eqnarray}
evaluated with ${\bf s} ({\bf y})$ replaced by ${\bf S} (F)$,
and $\partial_1 \cdots \partial_r {\bf s} ({\bf y})$ replaced by ${\bf S}_{F} ({\bf x}_1 \cdots {\bf x}_r)$.
So, setting
\begin{eqnarray*}
&&
T_{1 \cdots r} = T_{F} \left( {\bf x}_1, \ldots, {\bf x}_r \right),
\\
&&
S_{i 1 \cdots r} = S_{iF} \left( {\bf x}_1, \ldots, {\bf x}_r \right),
\\
&&
g_{ij \cdots} = \partial_i \partial_j \cdots g({\bf s})
\end{eqnarray*}
with $\partial_i = \partial / \partial s_i$ at ${\bf s} = {\bf S}(F)$, we have
\begin{eqnarray}
&&
T_1 = g_i S_{i1},
\
T_{12} = g_{ij} S_{i1} S_{j2} +g_i S_{i12},
\label{A3}
\\
&&
T_{123} = g_{ijk} S_{i1} S_{j2} S_{k3} + g_{ij} \sum^3 S_{i12} S_{j3} + g_i S_{i123},
\label{A4}
\\
&&
T_{1234} = g_{ijkl} S_{i1} S_{j2} S_{k3} S_{l4} + g_{ijk} \sum^6  S_{i1} S_{j2} S_{k34}
\nonumber
\\
&&
\qquad \qquad
+g_{ij} \left( \sum^4 S_{i1} S_{j234} + \sum^3 S_{i12} S_{j34} \right) + g_i S_{i1234},
\label{A5}
\end{eqnarray}
where summation over repeated suffixes $i,j, \ldots$ is implicit, and
by the multivariate version of Faa de Bruno's chain
rule given in Withers (1984), for $r \geq 1$,
\begin{eqnarray}
T_{1 \cdots r} = \sum_{k=1}^r g_{i_1 \cdots i_k} \left( {\bf S}(F) \right) \sum_{\bf n} \sum^{m({\bf n})} S_{i_1 \pi_1} \cdots S_{i_k \pi_k},
\label{A6}
\end{eqnarray}
where $\sum^{m({\bf n})}$ sums over all $m({\bf n}) = r!/ \prod_{i=1}^r (i!^{n_i} n_i!)$
partitions $(\pi_1 \cdots \pi_k)$ of $1 \cdots r$ giving
distinct terms with $n_i $ of the $\pi$'s of length $i$, and $\sum_{\bf n}$ sums
over $ \{ {\bf n}  \in \ N^r, \sum n_i = k, \sum in_i = r \}$.
For example,
\begin{eqnarray*}
\sum^3 S_{i12} S_{j34} = S_{i12}S_{j34} + S_{i13} S_{j24} + S_{i14} S_{j23}.
\end{eqnarray*}
The reader can derive $T_{123}$ from $T_{12}$ using equation (2.6) of Withers (1983) to appreciate the labor-saving this rule gives.

By [4c] of Comtet (1974) the general term can be written in terms of
the multivariate exponential Bell polynomials, $\{ B_{rk} ({\bf S})_{i_1 \cdots i_k} \}$:
\begin{eqnarray}
T_{1 \cdots r} = \sum_{k=1}^r g_{i_1 \cdots i_k} B_{rk} ({\bf S})_{i_1 \cdots i_k}.
\label{A7}
\end{eqnarray}
This is a much easier form to use than (\ref{A6}) as these
polynomials are immediately derived from the univariate
polynomials $B_{r_k} ({\bf S})$ tabled on pages 307-308 of Comtet (1974).
For example, the table gives
\begin{eqnarray*}
&&
B_{41} ({\bf S}) = S_4,
\\
&&
B_{42} ({\bf S}) = 4S_1 S_3 + 3 S_2^2,
\\
&&
B_{43} ({\bf S}) = 6 S_1^2 S_2,
\\
&&
B_{44} ({\bf S}) = S_1^4,
\end{eqnarray*}
so
\begin{eqnarray*}
&&
B_{41}({\bf S})_{i_1} = S_{i_11234},
\\
&&
B_{42}({\bf S})_{i_1i_2} = \sum^4 S_{i_11} S_{i_2234} + \sum^3 S_{i_11_2} S_{i_234},
\\
&&
B_{43}({\bf S})_{i_1i_2i_3} = \sum^6 S_{i_11}S_{i_22}S_{i_334},
\\
&&
B_{44} ({\bf S})_{i_1 \cdots i_4} = S_{i_11} \cdots S_{i_44},
\end{eqnarray*}
and (\ref{A7}) for $r \leq 4$ reduces to (\ref{A3})-(\ref{A5}).

Now suppose $F$ consists of $k$ d.f.s:  the only
change is to replace $({\bf x}_1 \cdots {\bf x}_r)$ by $\left(
\begin{array}{ccc} a_1& \cdots & a_r \\   {\bf x}_1 & & {\bf x}_r
\end{array} \right)$ wherever it occurs.
So, in the notation of (\ref{3.1}), (\ref{A3})-(\ref{A5}) imply
\begin{eqnarray}
T\left(a^2\right)
&=&
g_{ij} S_{ij} (a,a) + g_i S_i \left(a^2\right),
\label{A8}
\\
T\left(a^3\right)
&=&
g_{ijk} S_{ijk} (a,a,a) + 3g_{ij} S_{ij} \left(a,a^2\right) +g_i S_i \left(a^3\right),
\label{A9}
\\
T\left(a^4\right)
&=&
g_{ijkl} S_{ijkl} (a,a,a,a) + 6g_{ijk} S_{ijk} \left(a,a,a^2\right)
\nonumber
\\
&&
+ g_{ij} \left\{ 4S_{ij} \left(a,a^3\right) + 3S_{ij}\left(a^2,a^2\right) \right\} + g_iS_i (a^4),
\label{A10}
\\
T\left(a^2 b^2\right)
&=&
g_{ijkl} S_{ij} (a,a) S_{kl} (b,b)
\nonumber
\\
&&
+ g_{ijk} \left\{ S_{ij} (a,a) S_k \left(b^2\right) + S_{ij} (b, b) S_k \left(a^2\right) + 4 S_{ijk} (ab, a, b) \right\}
\nonumber
\\
&&
+g_{ij} \left\{ 2S_{ij} \left(a, ab^2\right) + 2S_{ij} \left(b, a^2 b\right) + S_i \left( a^2 \right)
S_j \left(b^2\right) + 2S_{ij}(ab,ab) \right\}
\nonumber
\\
&&
+ g_iS_i \left(a^2 b^2\right),
\label{A11}
\end{eqnarray}
where
\begin{eqnarray}
&&
S_{ij \cdots} \left( a^I,a^J, \ldots \right) =
\int S_{i_F}\left({}^{a^I}_{x^I} \right ) S_{j_F} \left( {}^{a^J}_{x^J} \right ) \cdots dF_a (x),
\label{A12}
\\
&&
\left( \begin{array}c a^I \\ x^I \end{array} \right)
= \begin{array}{ccc} a &\cdots& a\\ x && x \end{array}
\mbox{ with } I \mbox{ columns},
\\
&&
S_{ij} \left( a^Ib^J, \ldots, a^{K}b^{L}, \ldots \right) =
\int\cdots\int S_{iF}
\left( {}^{a^I}_{x^I} {}^{b^J}_{y^J} \cdots\right) S_{jF}
   \left( {}^{a^K}_{x^K} {}^{b^L}_{y^L} \cdots \right)  dF_a(x)dF_b(y),
\nonumber
\\
&&
\label{A13}
\end{eqnarray}
and so on.
Similarly, from (\ref{A7}) at $r=5$ we obtain
\begin{eqnarray}
T \left(a^2b^3\right) = \sum_{k=1}^5 g_{i_1 \cdots i_k} A^{i_1 \cdots i_k},
\label{A14}
\end{eqnarray}
where
\begin{eqnarray*}
&&
A^i = S_i \left(a^2 b^3\right),
\nonumber
\\
&&
A^{ij} = 2S_{ij} \left(a, ab^3 \right) +3 S_{ij} \left(b,a^2 b^2\right)+ S_i\left(a^2\right) S_j\left(b^3\right) + 6S_{ij} \left(ab, ab^2\right)+3S_{ij} \left(b^2,a^2b\right),
\nonumber
\\
&&
A^{ijk} = S_{ij} (a,a)S_k \left(b^3\right) + 3S_{ijk} \left(b,b,a^2b\right) + 6S_{ijk} \left(a,b,ab^2\right) +
6S_{ijk} \left(a, ab, b^2\right)
\nonumber
\\
&&
\qquad \qquad
+ 3S_{ik} \left(b,b^2\right) S_j \left(a^2\right) + 6S_{ijk}(b,ab,ab),
\nonumber
\\
&&
A^{ijkl} = S_i \left(a^2\right) S_{jkl}(b,b,b) + 6S_{ijkl} (ab,a,b,b) + 3S_{il}(b^2,b)S_{jk}(a,a),
\nonumber
\\
&&
A^{i_1 \cdots i_5} = S_{i_1 i_2} (a,a) S_{i_3 i_4 i_5}(b,b,b),
\nonumber
\end{eqnarray*}
and from (\ref{A7}) at $r=6$ we obtain
\begin{eqnarray}
T \left(a^2b^2c^2\right) = \sum_{k=1}^6g_{i_1 \cdots i_k}B^{i_1 \cdots i_k},
\label{A15}
\end{eqnarray}
where
\begin{eqnarray*}
&&
B^i = S_i \left(a^2b^2c^2\right),
\\
&&
B^{ij} = B_1^{ij} + B_2^{ij} +B_3^{ij},
\\
&&
B_1^{ij} = 2 \sum^3 S_{ij} \left( a,ab^2c^2 \right),
\\
&&
B_2^{ij} = \sum^3 S_i \left(a^2\right) S_j \left(b^2c^2 \right) +4 \sum^3 S_{ij} \left( ab, abc^2 \right),
\\
&&
B_3^{ij} = 2\sum^3 S_{ij} \left(a^2b, bc^2 \right) +4 S_{ij} (abc, abc),
\\
&&
B^{ijk}  = B_1^{ijk} +B_2^{ijk}+B_3^{ijk},
\\
&&
B_1^{ijk} = \sum^3 S_{ij} (a,a)S_k \left(b^2c^2\right) + 4 \sum^3 S_{ijk} \left(a, b, abc^2 \right),
\\
&&
B_2^{ijk} = 2\sum^6 S_{ik} \left(a,ac^2\right) S_j\left(b^2\right) + 4 \sum^6 S_{ijk} \left(a,ab,bc^2\right) + 8 \sum^3 S_{ijk} (a,bc,abc),
\\
&&
B_3^{ijk} = S_i(a^2)S_j \left(b^2\right) S_k\left(c^2\right) + 2\sum^3 S_i \left(a^2\right) S_{jk}(bc,bc) + 8S_{ijk} (ab,bc,ca),
\\
&&
B^{ijkl} = B_1^{ijkl} +B_2^{ijkl},
\\
&&
B_1^{ijkl} = 2\sum^6S_{ij} \left(a^2b,b\right) S_{kl}(c,c) + 8S_{ijkl}(abc,a,b,c),
\\
&&
B_2^{ijkl} = \sum^3 \big\{ S_{ij} (a,a)S_k \left(b^2\right) S_l\left(c^2\right) + 2 S_{ij} (a,a)S_{kl}(bc,bc)
\\
&&
\qquad \qquad
+ 4 S_{ijk} (a,b,ab)S_l(c^2) +8S_{ijkl} (a,b,ac,bc) \big\},
\\
&&
B^{i_1 \cdots i_5} = \sum^3 \left\{ S_{i_1}(a^2)S_{i_2i_3}(b,b) S_{i_4 i_5}(c,c) + S_{i_1 i_2 i_3} (ab,a,b) S_{i_4 i_5}   (c,c)\right\},
\\
&&
B^{i_1 \cdots i_6} = S_{i_1 i_2} (a,a) S_{i_3 i_4}(b,b) S_{i_5 i_6} (c,c),
\end{eqnarray*}
and $\sum^m$ is interpreted in the obvious manner by permuting $a,b,c$.
For example,
\begin{eqnarray*}
\sum^3S_{ij} \left(a,ab^2c^2\right) = S_{ij} \left(a,ab^2c^2\right) + S_{ij}\left(b,bc^2a^2\right) + S_{ij}\left(c,ca^2b^2\right).
\end{eqnarray*}
Similarly, if we now allow ${\bf T}$ and ${\bf g}$ to be $r$-vectors with
components $\{T^\alpha \}$ and $ \{g^\alpha  \}$, then by (\ref{A3}), $T^{\alpha\beta} (a,a)$ of (\ref{6.299}) is given by
\begin{eqnarray}
T^{\alpha \beta \cdots} (a, a, \ldots) = g_i^\alpha g_i^\beta \cdots S_{ij \cdots} (a, a, \ldots)
\label{A16}
\end{eqnarray}
and $T^{\alpha\beta} (ab,ab)$ of (\ref{6.799}) satisfies
\begin{eqnarray}
T^{\alpha\beta} (ab,ab) = g_{ij}^\alpha  g_{kl}^\beta  S_{ik}(a,a) S_{jl} (b,b)+
\sum_{\alpha\beta}^2 g_i^\alpha g_{jk}^\beta  S_{ijk} (ab,a,b) +g^\alpha _ig_j^\beta  S_{ij} (ab,ab),
\label{A17}
\end{eqnarray}
where
\begin{eqnarray*}
S_{ijk} (ab,a,b) = \int\int S_{iF}
\left(\begin{array}c ab\\xy \end{array} \right) S_{jF} \left(
\begin{array}c a\\x\end{array} \right) S_{kF}
\left( \begin{array}c b\\y\end{array} \right) dF_a (x) dF_b (y).
\end{eqnarray*}
Similarly, (\ref{6.599}), (\ref{6.699}) yield
\begin{eqnarray}
T^{\alpha\beta} \left( a^2, a \right) = \left \{ g_{ij}^\alpha  S_{ijk} (a,a,a)+g_i^\alpha S_{ik} \left(a^2, a\right) \right\} g_k^\beta,
\label{A18}
\end{eqnarray}
and
\begin{eqnarray}
T^{\alpha\beta} \left(a^2b, b\right)
&=&
\bigg\{ g_{ijk}^\alpha S_{ij}(a,a) S_{kl} (b,b) + g_{ij}^\alpha  \left[ S_i \left(a^2\right) S_{jl} (b,b) +
2 S_{ijl} (ab,a,b) \right]
\nonumber
\\
&&
+ g_i^\alpha S_{il} \left(a^2b, b\right)  \bigg\} g_{l}^\beta .
\label{A19}
\end{eqnarray}
Similarly,
\begin{eqnarray*}
T^{\alpha \beta \delta} (ab,a,b)= \left\{ g_{ij}^\alpha
S_{ijkl} (a,b,a,b) + g_i^\alpha S_{ikl} (ab, a,b) \right\} g_k^\beta g_l^\delta.
\end{eqnarray*}
We now consider the case, where ${\bf S}(F)$ is bivariate, that is $q = 2$.
Since $ S_{ij} (a^I,a^J)=S_{ji} (a^J,a^I)$, (\ref{A8})-(\ref{A11}) can be written as
\begin{eqnarray}
T\left(a^2\right)
&=&
\left\{g_{11}S_{11}+2 g_{12}S_{12}+g_{22}S_{22}\right\} (a,a)+\left\{g_1S_1+g_2S_2\right\} \left(a^2\right),
\label{A20}
\\
T\left(a^3\right)
&=&
\left\{g_{111}S_{111}+3g_{112}S_{112}+3 g_{122}S_{122}+g_{222}S_{222}\right\} (a,a,a)
\nonumber
\\
&&
+3 \left\{ g_{11}S_{11}+g_{12} \left(S_{12}+S_{21}\right)+g_{22}S_{22}\right\} \left(a,a^2\right)
\nonumber
\\
&&
+\left\{ g_1S_1+g_2S_2\right\} \left(a^3\right),
\label{A21}
\\
T\left(a^4\right)
&=&
\big\{ g_{1111}S_{1111}+4 g_{1112}S_{1112}+6g_{1122}S_{1122}+4 g_{1222}S_{1222}
\nonumber
\\
&&
+g_{2222}S_{2222}\big\}(a,a,a,a)+6 \big\{ g_{111}S_{111}+g_{112} S_{112}+2g_{121}S_{121}
\nonumber
\\
&&
+g_{221}S_{221}+2 g_{122}S_{122}+g_{222}S_{222}\big\}  \left(a,a,a^2\right)
\nonumber
\\
&&
+ 4 \left\{ g_{11}S_{11} + g_{12}(S_{12}+S_{21})+g_{22}S_{22}\right\}\left(a, a^3\right)
\nonumber
\\
&&
+ 3\left\{ g_{11}S_{11}+2g_{12} S_{12}+g_{22}S_{22}\right\} \left(a^2, a^2\right)
\nonumber
\\
&&
+\left\{g_1S_1+g_2S_2\right\} \left(a^4\right),
\label{A22}
\\
T\left(a^2b^2\right)
&=&
\big\{g_{1111}S_{11}S_{11}+2g_{1112}S_{11}S_{12}+g_{1122}S_{11}S_{22}
\nonumber
\\
&&
+ 2g_{1211}S_{12}S_{11}+4g_{1212}S_{12}S_{12}+2g_{1222}S_{12}S_{22}
\nonumber
\\
&&
+g_{2211}S_{22}S_{11}+2g_{2212}S_{22}S_{12}+g_{2222}S_{22}S_{22}\big\}  (a,a)(b,b)
\nonumber
\\
&&
+\big\{g_{111}S_{11}S_1+2g_{121}S_{12}S_1+g_{221}S_{22}S_1
\nonumber
\\
&&
+ g_{112}S_{11}S_2+2g_{122}S_{12}S_2+g_{222}S_{22}S_2 \big\}
\left\{(a,a) \left(b^2\right) + (b,b) \left(a^2\right) \right\}
\nonumber
\\
&&
+ 4\left\{ g_{111} S_{111} +3 g_{112} S_{1} 12 +3 g_{122} S_{122} + g_{222} S_{222} \right\} (ab,a,b)
\nonumber
\\
&&
+2 \left\{ g_{11}S_{11}+g_{12} \left(S_{12}+S_{21}\right) +g_{22}S_{22} \right\} \left\{ \left(a,ab^2\right) + \left(b,a^2 b\right) \right\}
\nonumber
\\
&&
+  \left\{g_{11}S_1S_1+g_{12} \left(S_1S_2+S_2S_1\right) + g_{22}S_2S_2  \right\} \left(a^2\right) \left(b^2\right)
\nonumber
\\
&&
+2 \left\{g_{11}S_{11}+2g_{12}S_{12}+g_{22}S_{22}\right\} (ab,ab) + \left\{g_1S_1+g_2S_2\right\} \left(a^2b^2\right).
\label{A23}
\end{eqnarray}
The convention here is that
\begin{eqnarray*}
&&
\left( g_{\pi_1} S_{\pi_2} + \cdots \right) \left( a^I,\ldots \right) = g_{\pi_1} S_{\pi_2} \left(a^I, \ldots\right),
\\
&&
\left( g_{\pi_1} S_{\pi_2} S_{\pi_3}+ \cdots \right) \left( a^I, \ldots \right) \left(b^J, \ldots \right) =
g_{\pi_1}S_{\pi_2} \left(a^I, \ldots \right) S_{\pi_3} \left(b^J, \ldots \right).
\end{eqnarray*}
Similarly, for $q=2$, splitting the third term in (\ref{A14}), $g_{ijk} A^{ijk}$,
into the six components corresponding to $A^{ijk}$, the first is
\begin{eqnarray*}
g_{ijk}S_{ijk} = \left\{g_{11k}S_{11k} +2g_{12k}S_{12k} +g_{22k}S_{22k}\right\}
\end{eqnarray*}
at $(a,a,b^3)$ and similarly for the second and sixth components.
Similarly, for the three components of the fourth term, the first being
\begin{eqnarray*}
g_{i \cdots l} S_{i \cdots l} = \left\{\sum_{j=1}^2g_{ijjj} S_{ijjj} + 3
g_{i112}S_{i112} + g_{i122} S_{i122}\right\}
\end{eqnarray*}
at $(a^2,b,b,b)$, and for the fifth term
\begin{eqnarray*}
g_{i_1 \cdots i_5} S_{i_1 \cdots i_5}
&=&
\left( g_{11-} S_{11-} + 2g_{12-} S_{12-} + g_{22-} S_{22-} \right)
\\
&&
\times
\left( g_{-111} S_{-111} + 3g_{-112} S_{-112}
+ 3g_{-122} S_{-122} + g_{-222}  S_{-222} \right)
\end{eqnarray*}
at $(a,a,b,b,b)$, where  $g_{\pi -} S_{\pi -} g_{-\pi^\prime}S_{-\pi^\prime}$ is interpreted as $g_{\pi\pi^\prime} S_{\pi \pi^\prime}$.

Similarly, for $q = 2$, the  term $B^{i}_3$ in (\ref{A15}) has the component
\begin{eqnarray*}
4g_{ij} S_{ij} = 4\sum^2_{i=1} g_{ii} S_{ii} + 8g_{12}S_{12}
\end{eqnarray*}
at $(abc,abc)$.
The sixth component is
\begin{eqnarray*}
&&
\left( g_{11-}  S_{11-}  + 2g_{12-} S_{12-} + g_{22-} S_{22-} \right) \left( g_{-11-} S_{-11-}  + 2g_{-12-} S_{-12-}  + g_{-22-} S_{-22-}\right)
\\
&&
\times \left( g_{-11} S_{-11} + 2g_{-12} S_{-12} + g_{-22} S_{-22} \right)
\end{eqnarray*}
at $(a,a,b,b,c,c)$, where $g_{\pi_1-}  S_{\pi_1-} g_{-\pi_2-} S_{-\pi_2-} g_{-\pi_3} S_{-\pi_3}$ interpreted as
$g_{\pi_1 \pi_2 \pi_3} S_{\pi_1 \pi_2 \pi_3}$, and so on.

\section*{Appendix B}
\renewcommand{\theequation}{$\mbox{B.\arabic{equation}}$}
\setcounter{equation}{0}

The nonparametric analogs of the terms for $t_2$ and equation (D.1) of
Withers (1987) needed for $T_2$ and $T_3$ - apart from those given in (\ref{3.3.1})-(\ref{3.3.3}) are as follows.
Summation over $a$, $b$, $c$ is implicit, where they occur.
These terms are listed both for the purpose of
checking and for application to other problems.
Note that $T_2$ requires
\begin{eqnarray*}
\left| \begin{array}c 22\\10 \end{array}
\right| = |3| \mbox{ and } \left| \begin{array}c 22\\20\end{array}
\right| = -2 \lambda_a^2 |2|_a
\end{eqnarray*}
and that $T_3$ requires
\begin{eqnarray*}
&&
\left| \begin{array}c 23\\10 \end{array}
\right| = \left| \begin{array}{cc}222\\ 110\end{array} \right|
= \lambda^3_a \left\{T \left(a^4\right) - T\left(a^2a^2\right) \right\},
\\
&&
\left| \begin{array}c 23\\20 \end{array} \right| = -2
\lambda_a^3T \left(a^3\right),
\\
&&
\left| \begin{array}{ccc}2&2&2\\ && \\ 1&0&0\end{array}\right| = |23|,
\
\left| \begin{array}{ccc} 2&2&2 \\ && \\ 2&0&0
\end{array}\right|_2 = -2\lambda_a^2\lambda_b T \left(a^2b^2\right),
\\
&&
\left| \begin{array}{ccc}2&2&2\\ && \\0&2&0 \end{array} \right|_1 =
-2\lambda_a^3 T \left(a^2a^2\right),
\
\left| \begin{array}{ccc}2&2&2\\1&2&0\end{array}\right|_i = \left| \begin{array}{ccc} 2&2&2 \\ &&
\\ 2&1&0 \end{array}\right| -2\lambda_a^3T \left(a^3\right) \mbox{ for } 1 \leq i \leq 3,
\\
&&
\left| \begin{array}{ccc} 2&2&2 \\ &&  \\ 2&2&0\end{array}\right| = 4 \lambda_a^3 T\left(a^2\right),
\
\left|
\begin{array}c32\\10\end{array}\right| = \lambda_a^3
\left\{ T \left(a^4\right) - 3T\left(a^2a^2\right) \right\},
\\
&&
\left| \begin{array}c 32\\20\end{array}\right| = -6|3|.
\end{eqnarray*}
Also,
\begin{eqnarray*}
\left|\begin{array}c 23\\30\end{array} \right| = \left|
\begin{array}{ccc} 2&2&2 \\ & & \\0&3&0 \end{array}\right| = \left|
\begin{array}{ccc} 2&2&2\\ & & \\0&4&0\end{array} \right|=0
\end{eqnarray*}
since  $\kappa_a (x_1 x_2)$, being quadratic in $F_a$, has functional derivatives higher than two equal to zero.
To illustrate the proof,
\begin{eqnarray*}
\left| \begin{array}c222\\210 \end{array} \right|_1
&=&
\kappa_{y_1z_1}^{x_1x_2} \kappa_{z_2}^{y_1y_2}   \kappa^{z_1z_2}  t_{x_1x_2y_2}
\\
&=&
\int^6 \lambda_a\lambda_b\lambda_c d_{\bf x} U_F \left(  \begin{array}{cc} b&c\\ y_1&z_1 \end{array} \right) d_{\bf y} V_F  \left(\begin{array}c c\\
z_2 \end{array} \right) d\kappa_c  \left(z_1z_2\right) T_F \left( \begin{array}{ccc} a&a&b\\
x_1 & x_2 & y_2
\end{array} \right),
\end{eqnarray*}
where $U(F) = \kappa^{x_1x_2} = \kappa_a (x_1  x_2)$ and $V(F) = \kappa^{y_1y_2} = \kappa_b(y_1y_2)$.
Note that
\begin{eqnarray*}
V_F \left( \begin{array}c c\\
z_2 \end{array} \right) = 0
\end{eqnarray*}
unless $c=b$ and
\begin{eqnarray*}
U_F \left( \begin{array}{cc} b&c\\
y_1&z_1
\end{array} \right) = 0
\end{eqnarray*}
unless $b=c=a$.
Also
\begin{eqnarray*}
U_F\left( \begin{array}{cc} a&a\\
y_1&z_1\end{array}
\right) = - \sum_{x_1x_2}^2 \triangle_{y_1} \left(x_1\right) \triangle_{z_1} \left(x_2\right),
\end{eqnarray*}
and
\begin{eqnarray*}
V_F \left(\begin{array}c a\\
z\end{array} \right) =
\triangle_{z} \left(y_1 \wedge y_2\right) - \sum^2_{y_1y_2} \triangle_{z} \left(y_1\right) F_a \left(y_2\right),
\end{eqnarray*}
where $\triangle_{y} (x) = (F_a (x))_{y} = I(y \leq x) - F_a (x)$.
Integrate first with respect to ${\bf x} =(x_1x_2)$: since columns in $T_F
\left( \begin{array}{ccc} \cdot & \cdot & \cdot \\ \cdot & \cdot & \cdot \end{array}\right)$
are interchangeable we may replace
$\sum_{x_1x_2}^2$ by 2.
Since
\begin{eqnarray}
\int T_F \left( \begin{array}{ccc} a&a&a \\ x_1&x_2&y_2
\end{array} \right) dF_a (x_i) =  0
\label{B1}
\end{eqnarray}
for $i = 1, 2$, and
\begin{eqnarray*}
d_{\bf x} \left\{ I \left(y_1 \leq x_1\right) I\left(z_1 \leq x_2\right)
\right \} = \delta \left(x_1-y_1\right) \delta \left(x_2-z_1\right) dx_1 dx_2
\end{eqnarray*}
with $\delta$ the Dirac delta function,
\begin{eqnarray*}
\int^2 d_{\bf x} U_F \left( \begin{array}{cc}
a&a\\ y_1 & z_1\end{array}\right) T_F \left(\begin{array}{ccc}
a&a&a\\x_1&x_2&y_2 \end{array}\right) = -2T_F\left(
\begin{array}{ccc} a&a&a\\y_1&z_1&y_2 \end{array} \right).
\end{eqnarray*}
So,
\begin{eqnarray*}
\left| \begin{array}c222\\210\end{array} \right|_1
= -2 \lambda_a^3 \int^4 d\kappa_a (z_1z_2) T_F
\left(\begin{array}{ccc} a&a&a\\y_1&y_2&z_1 \end{array} \right)
d_{\bf y} V_F \left( \begin{array}ca\\z_2\end{array} \right).
\end{eqnarray*}
Integrate with respect to ${\bf y} = (y_1y_2)$: (\ref{B1}) implies the contribution
from the last two out of the three terms in
$V_F \left(\begin{array}ca\\z \end{array}\right)$ is zero.
Also,
\begin{eqnarray*}
\triangle_{z} (y_1 \wedge y_2 )= I \left(z\leq y_1\right) I\left(z \leq y_2\right) - F_a \left(y_1\wedge y_2\right),
\end{eqnarray*}
so
\begin{eqnarray*}
d_{\bf y} \triangle_{z} \left(y_1 \wedge y_2\right) = \delta \left(y_1 - z\right)
\delta \left(y_2-z\right) dy_1dy_2 - \delta \left(y_1-y_2\right)dy_2 dF_a \left(y_2\right).
\end{eqnarray*}
So,
\begin{eqnarray*}
\int^2 T_F \left(\begin{array}{ccc}
a&a&a\\y_1&y_2&z_1\end{array} \right) d_{\bf y} V_F
\left(\begin{array}ca\\z_2\end{array} \right) = T_F
\left(\begin{array}{ccc} a&a&a\\z_2&z_2&z_1\end{array}\right) -
\int T_F \left( \begin{array}{ccc}a&a&a\\y_1&y_1&z_1\end{array}
\right) dF_a \left(y_1\right).
\end{eqnarray*}
Now integrate with respect to ${\bf z} = (z_1z_2)$:
by (\ref{B1}) the second out of two terms from $d \kappa_a(z_1z_2)$ contributes zero.
So, putting
\begin{eqnarray*}
L = \int dF_a \left(z_2\right) T_F \left(\begin{array}{ccc}
a&a&a\\y_1&y_1&z_2 \end{array}\right) = 0,
\end{eqnarray*}
we obtain
\begin{eqnarray*}
\left|\begin{array}c 222\\210\end{array}\right|_1
&=&
-2 \lambda_a^3 \int^2 dF_a \left(z_1 \wedge z_2\right) \left\{ T_F
\left(\begin{array}{ccc} a&a&a\\z_2&z_2&z_1\end{array}\right) -
\int T_F \left(\begin{array}{ccc}
a&a&a\\y_1&y_1&z_1\end{array}\right) dF_a \left(y_1\right) \right\}
\\
&=&
-2 \lambda_a^3 \left\{ \int T_F \left(\begin{array}{ccc}a&a&a\\z&z&z\end{array} \right) dF_a
(z) - \int dF_a\left(y_1\right)L \right\}= -2 \lambda_a^3 T \left(a^3\right).
\end{eqnarray*}

\section*{Appendix C}
\renewcommand{\theequation}{$\mbox{C.\arabic{equation}}$}
\setcounter{equation}{0}

Here, we show how to estimate $N$, the number of simulated samples
needed to estimate the bias to within a given relative error $\epsilon$.

Note that $T_{np}(\widehat{F})$ has bias $-n^{-p}T_p(F) + O(n^{-p-1})$ and that
$S_{np}(\widehat{F})$ has bias $-(n-1)^{-1}_p S_p(F) + O(n^{-p-1})= -n^{-p} S_p(F) + O( n^{-p-1})$.
Suppose we estimate the bias of $Y = S_{np}(\widehat{F})$ by $Z = \overline{Y} - T(F)$,
where $\overline{Y} = N^{-1} \sum^N_{j=1} Y_j$, $Y_j = S_{np}(\widehat{F}_j)$ and $\widehat{F}_j$ is the empirical distribution
of the $j$th simulated sample.
Then $EZ = ES_{np}(\widehat{F}) - T(F)$ is the true bias of $Y$ and we can write
$Z = EZ + (v_n/N)^{1/2} \{ \mathcal{N}(0,1) + o_p(1))\}$ as $N\to \infty$,
where $ v_n = var Y_1 = V_Tn^{-1} + O(n^{-2})$ as $ n \to \infty$,
and $ V_T=V_T(F) =  \sum \lambda_aT(a,a)$ with $ T(a,a) = \int T_F({}^a_x)^2 dF_a(x)$.
So, if $S_p=S_p(F) \neq 0$, the relative error in the estimate of bias,
\begin{eqnarray*}
\mbox{ (bias estimate - bias)/bias }
&\approx&
-\left(v_n/N\right)^{1/2} \mathcal{N}(0,1) n^pS_p(F)
\\
&\approx&
-V_T(F)^{1/2} S_p(F)^{-1} n^{p-1/2} N^{-1/2} \mathcal{N} (0,1)
\end{eqnarray*}
is bounded by a given number $\epsilon$ with probability greater than $0.975 + O_p(n^{-1/2})$  if
\begin{eqnarray*}
2 V_T(F)^{1/2} S_p(F)^{-1} n^{p-1/2} N^{-1/2} \leq \epsilon,
\end{eqnarray*}
that is, if
\begin{eqnarray*}
N \geq N_{\epsilon p n}=\epsilon^{-2} n^{2p-1} \phi_p,
\end{eqnarray*}
where $\phi_p = 4 V_T(F) S_p(F)^{-2}$.
This implies that for $\epsilon=0.1$ and $n$ large, say $n=100$, it is not practical
to carry out enough simulations to give meaningful estimates of bias unless $p=1$.
This is reflected by the poor estimates of bias in the tables
for the case $p=2$ obtained for $n=100$ using $N=10,000$.

Consider the following one sample examples.
Set $\beta_r = \mu_r\mu_2^{-2/2}$.
For $F = {\cal N}(0,1)$, $ \mu_4=3$, $\mu_6=15$, $\mu_8=105$
and for $ F = \exp(1)$, $\mu_2=1$, $\mu_3=2$, $\mu_4=9$, $\mu_5=44$, $\mu_6=305$, $ \mu_8=14,833$.

\noindent
{\bf Example C.1}
{\it Consider $T(F) = \mu_2$.
Then $V_T=\mu_4-\mu_2^2$, $S_1=\mu_2$, $\phi_1 = 4(\beta_4-1)$.
So, for a normal sample $ \phi_1 = 8$
and $\widehat{\mu}_2 = \mu_2(\widehat{F})$ needs
\begin{eqnarray*}
N \geq N_{ \epsilon 1 n} = 8\epsilon^{-2}n =  \left\{ \begin{array}{ll}
80,000 n \mbox { simulations}  & \mbox{for } \epsilon = 0.01,\\
800 n \mbox { simulations} & \mbox{for } \epsilon = 0.1.
\end{array} \right.
\end{eqnarray*}
For an exponential sample $ \phi_1 = 32$, so one needs four times as many simulations.
Since $S_2(F) = 0$, $\phi_2$ is not defined.}

\noindent
{\bf Example C.2}
{\it Consider $T(F) = \mu_2^2$.
Then $ V_T = 4\mu_2^2(\mu_4-\mu_2^2)$ and by Example 5.8,
$S_1 = -\mu_4+\mu_2^2$, $S_2 = -4\mu_4+7\mu_2^2$ so for a unit normal, $V_T = 8$, $S_1 = -2$, $\phi_1 = 8$, $S_2 = -29$, $\phi_2 = 0.1522$ so
$N_{0.1 1 n} = 800n$ and $N_{0.1 2 n} = 152n^3$ and for $ \exp(1)$, $V_T = 14,048$, $S_1 = 30$, $\phi_1 = 62.44$, $S_2 = 87$, $\phi_2 = 7.424$,
so $N_{0.1 1 n} = 6,244n$ and $N_{0.1 2 n} = 74.24n^3$.}

\noindent
{\bf Example C.3}
{\it Consider $T(F) = \mu_4$.
Then $V_T = \mu_8-\mu_4^2-8\mu_5\mu_3$, and by Example 5.6 or 5.10,
$S_1 = 2(2\mu_4-3\mu_2^2)$, $S_2 = 3(4\mu_4-7\mu_2^2)$, so
for a unit normal, $V_T = 96$, $S_1 = 6$, $\phi_1 = 32/3$, $S_2=15$, $\phi_2 = 128/75$, so
$ N_{0.1 1 n} = 1067n$ and $N_{0.1 2 n} = 171n^3$ and for $ \exp(1)$, $V_T = 14,048$, $S_1 = 30$,
$\phi_1 = 62.44$, $S_2 = 87$, $\phi_2 = 7.424$, so $ N_{0.1 1 n} = 6,244n$ and $N_{0.1 2 n} = 74.24n^3$.}

\noindent
{\bf Example C.4}
{\it Consider $T(F) =  \sigma = \mu^{1/2}_2$.
Then $V_T = \mu_2(\beta_4-1)/4$, so by Example 5.15, for a unit normal, $V_T = 1/2$, $S_1 = 3/4$, $\phi_1 = 32/9$, $S_2 = 1/32$, $\phi_2 = 2048$,
so  $N_{0.1 1 n} = 356n$ and $N_{0.1 2 n} = 204,800n^3$ and for $ \exp(1)$, $V_T = 2$, $S_1 = 3/2$, $\phi_1 = 32/9$,
$S_2 = 213/8 = 26.625$, $\phi_2 = 0.01129$, so $N_{0.1 1 n} = 356n$ and $N_{0.1 2 n} = 1.129n^3$.}

\section*{Appendix D}
\renewcommand{\theequation}{$\mbox{D.\arabic{equation}}$}
\setcounter{equation}{0}

Here, we list the non-zero derivatives $\mu_{r \cdot 12 \cdots p}=\mu_{rF}(x_1 \cdots x_p)$ for $2 \leq p \leq r \leq 6 $.
They are obtained from (\ref{5.699}) in terms of $h_i=\mu_{x_i}$, where  $\mu_x=x-\mu$, the first derivative of $\mu$:

\begin{eqnarray*}
&&
\mu_{2 \cdot 1}=h_1^2-\mu_2,
\\
&&
\mu_{2 \cdot 12}=-2h_1h_2,
\\
&&
\mu_{3 \cdot 1}=h_1^3-\mu_3 -3h_1\mu_2,
\\
&&
\mu_{3 \cdot 12}=-3\left( h_1^2-\mu_2 \right)h_2-3h_1\left( h_2^2-\mu_2 \right),
\\
&&
\mu_{3 \cdot 123}=12h_1h_2h_3,
\\
&&
\mu_{4 \cdot 1}=h_1^4-\mu_4-4h_1\mu_3,
\\
&&
\mu_{4 \cdot 12}= 12h_1h_2\mu_2-4\left( h_1^3-\mu_3 \right)h_2-4h_1\left( h_2^3-\mu_3 \right),
\\
&&
\mu_{4 \cdot 123}=12\left( h_1^2-\mu_2 \right)h_2h_3+12h_1\left( h_2^2-\mu_2 \right)h_3+12h_1h_2\left( h_3^2-\mu_2 \right),
\\
&&
\mu_{4 \cdot 1234}=-72h_1h_2h_3h_4,
\\
&&
\mu_{5 \cdot 1}=h_1^5-\mu_5-5h_1\mu_4,
\\
&&
\mu_{5 \cdot 12}=20h_1h_2\mu_3-5\left( h_1^4-\mu_4 \right)h_2-5h_1\left( h_2^4-\mu_4 \right),
\\
&&
\mu_{5 \cdot 123}= -60h_1h_2h_3\mu_2+20\left( h_1^3-\mu_3 \right)h_2h_3+20h_1\left( h_2^3-\mu_3 \right)h_3
+20h_1h_2\left( h_3^3-\mu_3 \right),
\\
&&
\mu_{5 \cdot 1234}= -60\left( h_1^2-\mu_2 \right)h_2h_3h_4-60h_1\left( h_2^2-\mu_2 \right)h_3h_4-60h_1h_2\left( h_3^2-\mu_2 \right)h_4
\\
&&
\qquad \qquad
-60h_1h_2h_3\left( h_4^2-\mu_2 \right),
\\
&&
\mu_{5 \cdot 12345}=480h_1h_2h_3h_4h_5,
\\
&&
\mu_{6 \cdot 1}=h_1^6-\mu_6-6h_1\mu_5,
\\
&&
\mu_{6 \cdot 12}= 30h_1h_2\mu_4-6\left( h_1^5-\mu_5 \right)h_2-6h_1\left( h_2^5-\mu_5 \right),
\\
&&
\mu_{6 \cdot 123}= -120h_1h_2h_3\mu_3 +30\left( h_1^4-\mu_4 \right)h_2h_3 +30h_1\left( h_2^4-\mu_4 \right)h_3+30h_1h_2\left( h_3^4-\mu_4 \right),
\\
&&
\mu_{6 \cdot 1234}/120= 3h_1h_2h_3h_4\mu_2-\left( h_1^3-\mu_3 \right)h_2h_3h_4
 -h_1\left( h_2^3-\mu_3 \right)h_3h_4
\\
&&
\qquad \qquad
-h_1h_2\left( h_3^3-\mu_3 \right)h_4 -h_1h_2h_3\left( h_4^3-\mu_3 \right),
\\
&&
\mu_{6 \cdot 12345}/360=\left( h_1^2-\mu_{2} \right)h_2h_3h_4h_5+h_1\left( h_2^2-\mu_2 \right)h_3h_4h_5
+h_1h_2\left( h_3^2-\mu_2 \right)h_4h_5
\\
&&
\qquad \qquad
+h_1h_2h_3\left( h_4^2-\mu_2 \right)h_5
+h_1h_2h_3h_4\left( h_5^2-\mu_2 \right).
\end{eqnarray*}
Note that
\begin{eqnarray*}
\mu_{r \cdot 12 \cdots r}=(-1)^{r-1}(r-1)r!\prod_{j=1}^{r}h_j,
\end{eqnarray*}
and
\begin{eqnarray*}
\mu_{r \cdot 12 \cdots r-1}=(-1)^{r}(r!/2)\sum^{r-1} \left(h_1^2-\mu_{2}\right) h_2 \cdots h_{r-1},
\end{eqnarray*}
where $\sum^{r-1}$ sums over all $r-1$ like terms.

\newpage

\begin{center}
{\bf Table 5.1.}~~Relative bias of $\widetilde{S}_{np} \widehat{(F)}$\\
for $T(F) = \mu^{-1}$ estimated from two runs of 5000 simulations.\\
\begin{tabular}{||l|l||rr|rr||}
\hline
&& $n=10$  & &   $n=100$  &\\
&& $p=1$ & $p=2$ & $p=1$ & $p=2$\\
\hline
Norm $(1/2, 1)$& Run 1 & 0.0773 & -0.0242 & 0.0089 & 0.0013 \\
& Run2 & 0.0916 & -0.0092 & 0.0087 & 0.0011 \\
\hline
Norm $(1,1)$& Run1 & -0.0780 & -0.0105 & -0.0149 & -0.0094  \\
            & Run2 & -0.0660 & -0.0040 & -0.0141 & -0.0087  \\
\hline
Norm $(2,1)$ &Run1 & 0.0208 & -0.0048 & -0.0046 & -0.0070  \\
       & Run2 & 0.0202 & -0.0056 & -0.0056 & -0.0078  \\
\hline
Exp $(1)$ & Run1 & 0.1096 & 0.0120 & 0.0052 &-0.0045   \\
          & Run2 & 0.1062 & 0.0184 & 0.0062 & -0.0035  \\
\hline
\end{tabular}
\end{center}

\begin{center}
{\bf Table 5.2.}~~Relative bias of $S_{np} \widehat{(F)}$\\
for $T(F) = \mu_4$ estimated from two runs of 60,000 simulations.
\begin{tabular}{||l|l||rr|rr|rr||}
\hline
&& $n=5$ & &  $n=10$  & &   $n=100$  &\\
&& $p=1$ & $p=2$ & $p=1$ & $p=2$ & $p=1$ & $p=2$\\
\hline
Norm $(0,1)$& Run1 & -0.3584 & -0.1988 & -0.1934 & -0.0543 & -0.0174 & 0.0021
\\
& Run2& -0.3572 & -0.1947 & -0.1871 & -0.0460 & -0.0206 & 0.0012
\\
\hline
Exp $(1)$ & Run1 & -0.4957 & -0.2861 & -0.2831 & -0.0754 & -0.0380 & -0.0063 \\
          & Run2 & -0.4943 & -0.2851 & -0.2964 & -0.0923 & -0.0399 & -0.0082 \\
\hline
\end{tabular}
\end{center}

\begin{center}
{\bf Table 5.3.}~~Relative bias of $S_{np} \widehat{(F)}$ for $T(F) = \sigma$.\\
\begin{tabular}{||l|l||rr|rr|rr||}
\hline
&& $n=5$ & &  $n=10$  & &   $n=100$  &\\
&& $p=1$ & $p=2$ & $p=1$ & $p=2$ & $p=1$ & $p=2$\\
\hline
Normal $(0,1)$& Run1 & -0.1578 & -0.0265 & -0.0764 & -0.0082 & 0.0281 & -0.0045
\\
    & Run2& -0.1592 & -0.0277 & -0.0745 & -0.0080 & 0.0003 & 0.0031
\\
\hline
Exp $(1)$ & Run1 & -0.2278 & -0.1019 & -0.1251 & -0.0422 & -0.0158 & -0.0029\\
          & Run2 & -0.2331 & -0.1084 & -0.1206 & -0.0422 & -0.0176 & -0.0004\\
\hline
Number of & &&&&&&\\
simulations/run & & 10,000 && 30,000 && 30,000 &\\
\hline
\end{tabular}
\end{center}


\begin{thebibliography}{999}



\bibitem{}
Cabrera, J. and Fernholz, L. T. (1999).
Target estimation for bias and mean square error reduction.
{\it Annals of Statistics}, {\bf 27}, 1080-1104.


\bibitem{}
Cabrera, J. and Fernholz, L. T. (2004).
Multivariate targeting with applications to ellipse estimation.
Contemporary data analysis: theory and methods.
{\it Journal of Statistical Planning and Inference}, {\bf 122}, 79-94.


\bibitem{}
Comtet, L. (1974).
{\it Advanced Combinatorics}.
Reidel, Dordrecht.


\bibitem{}
Davison, A. C., Hinkley, D. V. and Schechtman, E. (1986).
Efficient bootstrap simulation.
{\it Biometrika}, {\bf 73}, 555-566.

\bibitem{}
Efron, B. (1981).
Nonparametric estimates of standard errors: The jackknife, the bootstrap and other methods.
{\it Biometrika}, {\bf 68}, 589-599.


\bibitem{}
Efron, B. (1982).
{\it The Jackknife, the Bootstrap and Other Resampling Plans}.
SIAM, Philadelphia.

\bibitem{}
Fernholz, L. T. (2001).
On multivariate higher order von Mises expansions.
{\it Metrika}, {\bf 53}, 123-140.


\bibitem{}
Fisher, R. A. (1929).
Moments and product moments of sampling distributions.
{\it Proceedings of the London Mathematical Society}, 2, {\bf 30}, 199-238.

\bibitem{}
Gray, H. L. and Schucany, W. R. (1972).
{\it The Generalized Jackknife Statistic}.
Marcel Dekker, New York.


\bibitem{}
Hall, P. (1992).
{\it The Bootstrap and Edgeworth Expansion}.
Springer, New York.

\bibitem{}
Jaeckel, L. A. (1972).
The infinitesimal jackknife.
{\it Bell Laboratories Technical Report MM 72-1215-11}, June 30, 1972, New Jersey.


\bibitem{}
James, G. S. (1958).
On moments and cumulants of systems of statistics.
{\it Sankhy\=a}, {\bf 20}, 1-30.

\bibitem{}
Kendall, M. G. and Stuart, A. (1977).
{\it The Advanced Theory of Statistics}, volume 1, fourth edition.
Griffin, London.


\bibitem{}
Quenouille, M. H. (1956).
Notes on bias in estimation.
{\it Biometrika}, {\bf 43}, 353-360.

\bibitem{}
Schucany, W. R., Gray, H. and Owen, D. B. (1971).
On bias reduction in estimation.
{\it Journal of the American Statistical Association}, {\bf 66}, 524-533.


\bibitem{}
Sen, P. K. (1988).
Functional jackknifing: rationality and general asymptotics.
{\it Annals of Statistics}, {\bf 16}, 450-469.


\bibitem{}
Stuart, A. and Ord, J. K. (1987).
{\it Kendall's Advanced Theory of Statistics}, volume 1, fifth edition.
Griffin, London.

\bibitem{}
Sukhatme, P. V. (1944).
Moments and product moments of moment-statistics
for samples of the finite and infinite populations.
{\it Sanhky\=a}, {\bf 6}, 363-382.

\bibitem{}
Tiit, E. (1988).
Unbiased and $n^{-k}$-estimations of entire rational functions of moments.
In: {\it Statistical and Probabilistic Models}, Report 798, Tartu State University, Estonia, pp. 3-17.

\bibitem{}
Wishart, J. (1952).
Moment coefficients of the $k$-statistics in samples from finite population.
{\it Biometrika}, {\bf 39}, 1-13.


\bibitem{}
Withers, C. S. (1982).
The distribution and quantiles of a function of
parameter estimates.
{\it Annals of the Institute of Statistical Mathematics}, A, {\bf 34}, 55-68.

\bibitem{}
Withers, C. S. (1983).
Expansions for the distribution and quantiles
of a regular functional of the empirical distribution with applications to
nonparametric confidence intervals.
{\it Annals of Statistics}, {\bf 11}, 577-587.


\bibitem{}
Withers, C. S. (1984).
A chain rule for differentiation
with applications to multivariate Hermite polynomials.
{\it Bulletin of the Australian Mathematical Society}, {\bf 30}, 247-250.

\bibitem{}
Withers, C. S. (1987).
Bias reduction by Taylor series.
{\it Communications in Statistics---Theory and Methods}, {\bf 16}, 2369-2384.


\bibitem{}
Withers, C. S. (1988).
Nonparametric confidence intervals for
functions of several distributions.
{\it Annals of the Institute of Statistical Mathematics}, {\bf 40}, 727-746.

\bibitem{}
Withers, C. S. and Nadarajah, S. (2008).
Analytic bias reduction for $k$ sample functionals.
{\it Sankhy\=a}, A, {\bf 70}, 186-222.


\end{thebibliography}
\end{document}